\documentclass[11pt]{article} \pdfoutput=1 

\newcommand\manualformatting{}

\usepackage{authblk} 
\usepackage{fullpage}

\usepackage{bbold} 
\usepackage{bbm}
\usepackage{float}
\usepackage{euscript}
\usepackage{graphicx} 
\usepackage{extarrows} 

\usepackage{mathtools} 
\usepackage{amsmath}
\usepackage{amssymb}
\usepackage{amsthm}
\usepackage{enumerate, latexsym, mathrsfs, mdwlist}
\usepackage{booktabs}

\usepackage[color,all,2cell,arrow,matrix]{xy} \UseAllTwocells \SilentMatrices

\usepackage{array} 

\usepackage[dvipsnames]{xcolor}
\definecolor{Ocean}{HTML}{016064}

\newcommand{\eqn}[1]{\begin{equation}#1\end{equation}}
\newcommand{\eqnn}[1]{\begin{equation*}#1\end{equation*}}
\newcommand{\agn}[1]{\begin{align*}#1\end{align*}}
\newcommand{\agnn}[1]{\begin{align}#1\end{align}}
\newcommand \pf			{\begin{proof}}
\newcommand \epf			{\end{proof}}
\newcommand \bnu			{\begin{enumerate}}
\newcommand \enu			{\end{enumerate}}
\newcommand \bit             	{\begin{itemize}}
\newcommand \eit             	{\end{itemize}}

\newcommand \diagram		{\xymatrix}          
\newcommand{\newdef}[1]{\textit{#1}}
\newcommand \defdtobe		{\coloneqq}
\renewcommand{\to}			{\longrightarrow}

\newcommand\sboxtimes		{{_\boxtimes}}
\newcommand\fib			{\mathrm{Fib}}
\newcommand\eend			{\mathrm{end}}
\newcommand \funu			{\mathrm{Fun}^{\mathrm{u}}}
\newcommand \repu			{\mathrm{Rep}^\ast}
\newcommand\Tube			{\mathrm{Tube}}
\newcommand\vectG		{\mathrm{Vec}_G}
\newcommand\Cu			{\mathfrak{C}}
\newcommand\Mu			{\mathfrak{M}}
\newcommand\Au			{\Aaa}

\newcommand \vect			{\mathrm{Vec}_\kb}
\newcommand \mone            	{\mathbf{1}} 
\newcommand \nat			{\mathrm{Nat}}
\newcommand \fun		      	{\mathrm{Fun}}
\newcommand \Hom          	{\mathrm{Hom}}
\newcommand \Id		       	{\mathrm{Id}}
\newcommand \id          		{\mathrm{id}}
\newcommand \op   		    	{\mathrm{op}}
\newcommand \cop   		{\mathrm{cop}}
\newcommand \rev			{\mathrm{rev}}
 
\newcommand \rhs 			{R.H.S.} 
\newcommand \os			{\otimes}
\newcommand \bos			{\mathbin{\overline{\os}}}
\newcommand \Oplus		{\bigoplus}
\newcommand \funto			{\Rightarrow} 
\newcommand \isom			{\overset{\sim}{\to}}
\newcommand \ladj			{\dashv}

\newcommand \inverse		{^{-1}}
\renewcommand \epsilon		{\varepsilon}
\renewcommand \leq			{\leqslant} 
\renewcommand \geq		{\geqslant}
\renewcommand \:			{\colon}

\newcommand \<			{\langle}
\renewcommand \>			{\rangle}

\newcommand \kb			{k}

\newcommand \Aa			{{\mathcal A}}
\newcommand \Bb			{{\mathcal B}}
\newcommand \Cc			{{\mathcal C}}
\newcommand \Dd			{{\mathcal D}}
\newcommand \Ee			{{\mathcal E}}

\newcommand \Gg			{{\mathcal G}}

\newcommand \Kk			{{\mathcal K}}
\newcommand \Ll			{{\mathcal L}}
\newcommand \Mm			{{\mathcal M}}
\newcommand \Nn			{{\mathcal N}}

\newcommand \Rr			{{\mathcal R}}

\newcommand \Zz			{{\mathcal Z}}

\def\Aaa					{\mathscr{A}}

\def\Fff					{\mathscr{F}}
\def\Ggg					{\mathscr{G}}

\def\Uuu					{\mathscr{U}}
\def\Vvv					{\mathscr{V}}

\newcommand \fa			{\mathbf{a}}

\newcommand \fe 			{\mathbf{e}}
\newcommand \ff			{\mathbf{f}}

\newcommand \fh			{\mathbf{h}}

\DeclareMathAlphabet{\mathcal}{OMS}{cmsy}{m}{n}	
\DeclareMathAlphabet{\mathsf}{OT1}{cmss}{m}{n}	

\newcommand \dmo			{\DeclareMathOperator}

\dmo \Coev 				{\mathrm{Coev}}
\dmo \Ev					{\mathrm{Ev}}
\dmo \coev       				{coev}
\dmo \ev          				{ev}
\dmo \im 		          		{im}
\dmo \Irr                            		{Irr}
\dmo \rep                          		{Rep}
\dmo \BiMod				{BiMod}
\dmo \LMOD 				{\Ll\Mm\mathrm{od}}
\dmo \Mor          			{Mor}
\dmo \End	          			{End}

\usepackage{varioref}    
\usepackage{hyperref}      
\usepackage{cleveref}   
\hypersetup{
   colorlinks=true,
   linkcolor=Ocean,
   filecolor=blue,      
   urlcolor=blue,
   citecolor=blue
}  

\theoremstyle{plain}

\newtheorem{mainthm}{Theorem}

\swapnumbers

\newtheorem{thm}{Theorem}[section]
\newtheorem{prp}[thm]{Proposition}
\newtheorem{lem}[thm]{Lemma}
\newtheorem{cor}[thm]{Corollary}

\theoremstyle{definition}

\newtheorem{dfn}[thm]{Definition}
\newtheorem{ntn}[thm]{Notation}
\newtheorem{expl}[thm]{Example}
\newtheorem{rmk}[thm]{Remark}

\Crefname{lem}{Lemma}{Lemmas}
\Crefname{expl}{Example}{Examples}
\Crefname{thm}{Theorem}{Theorems}
\Crefname{mainthm}{Theorem}{Theorems}

\usepackage{tikz}
\usetikzlibrary{matrix,arrows,decorations}

\title{On the Representation Categories of Weak Hopf Algebras Arising from Levin-Wen Models}

\author[1,2,3]{Ansi Bai\thanks{\href{mailto:crippledbai@163.com}{\tt crippledbai@163.com}}}
\author[1,2,3,4]{Zhi-Hao Zhang\thanks{\href{mailto:zhangzhihao@bimsa.cn}{\tt zhangzhihao@bimsa.cn}}}

\affil[1]{Beijing Institute of Mathematical Sciences and Applications, Beijing 101408, China}
\affil[2]{International Quantum Academy, Shenzhen 518048, China}
\affil[3]{Guangdong Provincial Key Laboratory of Quantum Science and Engineering, Southern University of Science and Technology, Shenzhen 518055, China}
\affil[4]{Wu Wen-Tsun Key Laboratory of Mathematics of Chinese Academy of Sciences,
School of Mathematical Sciences,
University of Science and Technology of China, Hefei 230026, China}

\begin{document}

\date{}
\maketitle

\begin{abstract}
	In their study of Levin-Wen models [Commun. Math. Phys. 313 (2012) 351-373], Kitaev and Kong proposed a weak Hopf algebra associated with a unitary fusion category $\Cc$ and a unitary left $\Cc$-module $\Mm$, and sketched a proof that its representation category is monoidally equivalent to the unitary $\Cc$-module functor category $\funu_\Cc(\Mm,\Mm)^\rev$. We give an independent proof of this result without the unitarity conditions. In particular, viewing $\Cc$ as a left $\Cc \boxtimes \Cc^\rev$-module, we obtain a quasi-triangular weak Hopf algebra whose representation category is braided equivalent to the Drinfeld center $\Zz(\Cc)$. 
In the appendix, we also compare this quasi-triangular weak Hopf algebra with the tube algebra $\Tube_\Cc$ of $\Cc$ when $\Cc$ is pivotal. These two algebras are Morita equivalent by the well-known equivalence $\rep(\Tube_\Cc)\cong\Zz(\Cc)$. However, we show that in general there is no weak Hopf algebra structure on $\Tube_\Cc$ such that the above equivalence is monoidal.
\end{abstract}

\setcounter{tocdepth}{2}
\tableofcontents

\section*{Introduction} 
\addcontentsline{toc}{section}{Introduction}

Given a unitary fusion category $\Cu$ and a finite unitary left $\Cu$-module $\Mu$, Kitaev and Kong \cite{Kitaev_Kong_2012} introduced an algebra (see also \cite[\S 6]{Morrison_Walker_2012}), denoted by $\Au_\Mu^\Cu$, to study the topological excitations on the boundaries of the Levin-Wen model \cite{Levin_Wen_2005}. The study relies on two observations they made about this algebra:
\bnu[(Obs.~I)]
	\item \label{item1:intro} $\Au_\Mu^\Cu$ is a $C^\ast$-weak Hopf algebra with the structure maps given in \cite[\S 4]{Kitaev_Kong_2012}.
	\item \label{item2:intro} There is an equivalence of unitary multi-tensor categories
	\eqn{\label{eq:you_need_to_ref}
		\repu(\Au_\Mu^\Cu)\isom \funu_{\Cu}(\Mu,\Mu)^\rev\,,
	}
	where $\repu(\Au_\Mu^\Cu)$ is the category of finite-dimensional $\ast$-representations over $\Au_\Mu^\Cu$, and $\funu_{\Cu}(\Mu,\Mu)$ is the category of unitary $\Cu$-module functors from $\Mu$ to itself.
\enu
The proofs of these two observations were known to the authors of \cite{Kitaev_Kong_2012}, but were not fully written out in the original article, apart from an outline of a proof of (Obs. \ref{item2:intro}). Later, a proof of (Obs. \ref{item1:intro}) based on graphical calculus appeared in a detailed study of $\Au_\Mu^\Cu$ \cite{Jia_Tan_Kaszlikowski_2024}; see also \cite[\S 3.2]{Cordova_Holfester_Ohmori_2024}. The proof of (Obs.~\ref{item2:intro}) remains intricate:  although the proof was known \cite{Kitaev_Kong_2012} (see also \cite{Barter_Bridgeman_Jones_2019a} for a partial proof), a detailed proof has yet to be published despite the passage of time. As a result, (Obs.~\ref{item2:intro}) has become folklore among physicists. 
This folklore, on the other hand, holds considerable importance due to its relevance both in Levin-Wen models and in various other contexts. For example, the original article \cite{Kitaev_Kong_2012} applied (Obs.~\ref{item2:intro}) to classify the topological excitations on the $_\Cu\Mu$-boundary of the Levin-Wen model as unitary $\Cu$-module endofunctors on $\Mu$. To be more specific, they first identified those excitations as $\ast$-representations over $\Au_\Mu^\Cu$ via physical principles, and then utilized (Obs.~\ref{item2:intro}) to complete the classification. The equivalence \eqref{eq:you_need_to_ref} also finds broader applications in a variety of contexts, such as designing algorithms for computing $F$-symbols \cite{Barter_Bridgeman_Wolf_2022}, studying $S^1$-parameterized families of general $\Cu$-symmetric gapped systems \cite{Inamura_Ohyama_2024}, and analyzing twisted boundary states and entanglement entropy in conformal field theory \cite{Choi_Rayhaun_Zheng_2024b, Choi_Rayhaun_Zheng_2024a}.
Additionally, a recent study \cite{Gagliano_Grigoletto_Ohmori_2025} assumes a higher-categorical version of the equivalence \eqref{eq:you_need_to_ref} to study phases in Yang-Mills theory. Some other aspects of the algebra $\Au_\Mu^\Cu$ are also addressed or employed in numerous studies, including but not limited to \cite{Kong_2012, Lan_Wen_2014, Barter_Bridgeman_Jones_2019a, Bridgeman_Barter_Jones_2019b, Bridgeman_Barter_2020a, Bridgeman_Barter_2020b, Bridgeman_Lootens_Verstraete_2023, Jia_Tan_Kaszlikowski_2024, Cordova_Holfester_Ohmori_2024}.

Considering the importance of (Obs. \ref{item2:intro}) and its widespread applications, a detailed proof would be a valuable addition to the literature. In this work, we partially address this need by proving (Obs. \ref{item2:intro}) without assuming the unitarity conditions. Along the way we also prove the non-unitary version of (Obs. \ref{item1:intro}). To be concrete, for a fusion category $\Cc$ and a finite semisimple left $\Cc$-module $\Mm$, we define a weak Hopf algebra $A_\Mm^\Cc$, which is the non-$C^\ast$-version of the algebra $\Au_\Mu^\Cu$. Then we show
\begin{mainthm}[\Cref{thm:main_equivalence}]\label{mainthm:1}
	There is a monoidal equivalence
	\eqn{\label{eq0:intro}
		\rep(A_\Mm^\Cc)\isom\fun_\Cc(\Mm,\Mm)\,.
	}
	Here $\rep(A_\Mm^\Cc)$ denotes the category of finite-dimensional left $A_\Mm^\Cc$-modules, and $\fun_\Cc(\Mm,\Mm)$ denotes the category of left $\Cc$-module functors from $\Mm$ to itself.
\end{mainthm}
We again note that a sketchy proof of the unitary version of \Cref{mainthm:1} is already provided in \cite[\S 4]{Kitaev_Kong_2012}; see \Cref{rmk:KK12_pf} for more discussion. However, the proof we provide here is slightly more conceptual and basis-independent. In particular, our formulation of $A_\Mm^\Cc$ reduces to the very concise form in \Cref{rmk:ih} if one employs the language of internal homs. These conceptual simplifications also enable proving certain generalization of \Cref{mainthm:1} \cite{Bai_Zhang} and potentially its higher categorical analogues. We also note that the equivalence (as categories) in \eqref{eq0:intro} could be derived from \cite[Proposition 10]{Barter_Bridgeman_Jones_2019a}, which is in turn based on \cite{Morrison_Walker_2012}. 
However, the complete derivation of \Cref{mainthm:1} from the results in \cite{Barter_Bridgeman_Jones_2019a} would require additional work, which, to the authors' awareness, has not been addressed in the literature. The proof presented in this article is independent of \cite{Barter_Bridgeman_Jones_2019a}.

Our second main result examines a special case of the equivalence \eqref{eq0:intro}, as follows. Note that $\Cc$ can be viewed as a left $\Cc\boxtimes\Cc^\rev$-module via treating $\Cc$ as the regular $\Cc$-$\Cc$-bimodule. Applying \Cref{mainthm:1}, we obtain a monoidal equivalence:
\eqn{
\label{eq1:intro}
	\rep(A_\Cc^{\Cc\boxtimes\Cc^\rev})\isom \fun_{\Cc\boxtimes\Cc^\rev}(\Cc,\Cc)\isom\Zz(\Cc)\,,
}
where $\Zz(\Cc)$ denotes the Drinfeld center of $\Cc$. It is well-established that $\Zz(\Cc)$ is a braided monoidal category, and that braidings on the representation category of a weak Hopf algebra are in 1:1 correspondence to quasi-triangular structures on the algebra. Our second main result is the explicit expression of the quasi-triangular structure on $A_\Cc^{\Cc\boxtimes\Cc^\rev}$ corresponding to the braiding on $\Zz(\Cc)$:
\begin{mainthm}[\Cref{thm:qt_main}]\label{mainthm:2}
	The quasi-triangular structure $\Rr$ on the weak Hopf algebra $A_\Cc^{\Cc\boxtimes\Cc^\rev}$, which corresponds to the braiding on $\Zz(\Cc)$ via the equivalence \eqref{eq1:intro}, is given by \eqref{eq:qt:1} or equivalently \eqref{eq:qt:2}. In particular, \eqref{eq1:intro} becomes a braided monoidal equivalence when $A_\Cc^{\Cc\boxtimes\Cc^\rev}$ is equipped with $\Rr$.
\end{mainthm}
In particular, \Cref{mainthm:2} offers a way to realize $\Zz(\Cc)$ as the representation category of certain quasi-triangular weak Hopf algebra. We give the explicit form of $(A_\Cc^{\Cc\boxtimes\Cc^\rev},\Rr)$ when $\Cc=\vectG^\omega$.

When $\Cc$ is a pivotal fusion category, it is well-known that there is an equivalence of categories
\eqn{
\label{eq2:intro}
	\rep(\Tube_\Cc)\isom\Zz(\Cc)\,,
}
where $\Tube_\Cc$ is Ocneanu's tube algebra associated with $\Cc$ \cite{Ocneanu_1994,Izumi_2000, Mueger_2003}. This means that $A_\Cc^{\Cc\boxtimes\Cc^\rev}$ and $\Tube_\Cc$ are Morita equivalent. It is also well-known by physicists that the algebra $A_\Cc^\Mm$, and in particular $A_\Cc^{\Cc\boxtimes\Cc^\rev}$ is related to $\Tube_\Cc$ \cite{Bridgeman_Barter_2020b, Jia_Tan_Kaszlikowski_2024}. These facts motivate a precise comparison between $A_\Cc^{\Cc\boxtimes\Cc^\rev}$ and $\Tube_\Cc$, which we undertake in an appendix by highlighting the following points:
\bnu
	\item The equivalence \eqref{eq1:intro} does not require a pivotal structure on $\Cc$, whereas the existence of \eqref{eq2:intro}, to the best of the authors' knowledge, does. 
	\item \label{item2:difference} There exists an infinite family of algebras that are Morita equivalent to $\Tube_\Cc$, and $\Tube_\Cc$ is the smallest one \cite{Mueger_2003}. When $\Cc$ is pivotal, $A_\Cc^{\Cc\boxtimes\Cc^\rev}$ lies within this family. 
	\item \label{item3:difference} In general, $\Tube_\Cc$ does not carry a weak Hopf algebra structure such that the induced monoidal structure on $\rep(\Tube_\Cc)$ renders \eqref{eq2:intro} a monoidal equivalence.
\enu
In particular, point \ref{item3:difference} presents a sharp contrast between $A_\Cc^{\Cc\boxtimes\Cc^\rev}$ and $\Tube_\Cc$. There are many nice works on tube algebras and the equivalence \eqref{eq2:intro}, or their variants \cite{Ocneanu_1994, Izumi_2000, Mueger_2003, Morrison_Walker_2012, Ghosh_Jones_2016, Popa_Shlyakhtenko_Vaes_2018, Hoek_2019, Liu_Ming_Wang_Wu_2023, Lan_2024}. In light of points \ref{item2:difference} and \ref{item3:difference}, we hope that our work, building on \cite{Kitaev_Kong_2012}, opens a new direction of studying the coalgebraic aspects of the algebras Morita equivalent to tube algebras (or their variants). This differs from existing works, although \cite{Neshveyev_Yamashita_2018} explores some coalgebraic structures of tube algebras in a different context.

The main tool used in our proof of \Cref{mainthm:1} is the reconstruction theorem for (finite-dimensional) weak Hopf algebras, also known as the Tananka-Krein duality for weak Hopf algebras \cite{Hayashi_1999, Szlachanyi_2000, Szlachanyi_2004}. This theorem asserts that a weak Hopf algebra $A^\Fff$ can be constructed from a finite multi-tensor category $\Dd$ together with a faithful exact separable Frobenius functor $\Fff\:\Dd\to\vect$ from $\Dd$ to the category of finite-dimensional vector spaces. Furthermore, there exists an equivalence $\rep(A^\Fff)\cong\Dd$ of monoidal categories. The strategy of our proof of \Cref{mainthm:1} is to recognize $A_\Mm^\Cc$ as the weak Hopf algebra constructed from certain faithful exact separable Frobenius functor $\fun_\Cc(\Mm,\Mm)\to\vect$. We emphasize that based on \cite{Hayashi_1999, Kitaev_Kong_2012}, we have obtained an explicit \emph{presentation} of $A_\Mm^\Cc$: once $\Cc$ and $\Mm$ are known, a basis for $A_\Mm^\Cc$ can be written, and its weak Hopf algebra structure can be expressed in terms of this basis. In particular, by considering a special case of our reconstruction process, for any given fusion category $\Cc$, one obtains an explicit presentation of a weak Hopf algebra whose representation category is equivalent to $\Cc$\footnote{This could be an answer to a mathoverflow question concerning reconstructions: \url{https://mathoverflow.net/questions/453975/how-does-the-tannaka-duality-work-for-weak-hopf-algebras-and-fusion-categories}.}; this is illustrated in \Cref{subsec:right_regular}. We hope that our explicit presentation of $A_\Mm^\Cc$ could serve as a non-unitary complement to \cite{Kitaev_Kong_2012}, offering a useful tool for physicists working with (non-unitary) fusion categories and their modules. 

We remark that the weak Hopf algebra $A_\Mm^\Cc$ fits into a broader class of algebraic structures in Levin-Wen models. This broader class was outlined in \cite[\S 6]{Kitaev_Kong_2012}, and \cite[\S VI]{Lan_Wen_2014}, based on \cite{Kong, Kong_2012}.
To incorporate this larger class of algebras, one needs to suitably generalize the concept of weak Hopf algebras. The present note serves as basis of the authors' future investigation into these generalized weak Hopf algebra structures (see \Cref{rmk:on_multi_object} for more discussion).

\Cref{sec:reconstruct} is devoted to the reconstruction theorem for weak Hopf algebras. \Cref{sec:2} and \Cref{sec:qt} are devoted to the proof of \Cref{mainthm:1} and \Cref{mainthm:2}, respectively. We compare $A_\Cc^{\Cc\boxtimes\Cc^\rev}$ with $\Tube_\Cc$ in \Cref{sec:comparison}.

Throughout this paper, we fix an algebraic closed field $\kb$ of characteristic $0$. All vector spaces, algebras and modules over $\kb$ are assumed to be finite-dimensional, although we will emphasize it whenever necessary. Algebras over $\kb$ are assumed to be associative with unit, and algebra homomorphisms are assumed to preserve units. Similar assumptions apply to coalgebras. For an algebra $A$ and a left $A$-module $(M,\rho)$, we use the notation $\rho(a\os m)=a.m$ for $a\in A$, $m\in M$; the notation is similar for right modules. We use the term ``$A$-representation'' as an synonym for a left $A$-module. Functors and equivalences between $\kb$-linear categories are implicitly assumed to be $\kb$-linear. 

For general facts on monoidal categories, as well as on fusion categories over $\kb$ and their module categories, we refer the reader to \cite{Etingof_Gelaki_Nikshych_Ostrik_2015}.

\section*{Acknowledgements}
This work has been reported in a talk by the first author at the 2024-2025 BIMSA TQFT and Higher Symmetries Seminar. The authors thank Liang Kong for proposing the project and for his suggestions on an early draft, which improved the article. The authors also thank Zhian Jia, Tian Lan, Chenqi Meng, Xiaoqi Sun, Sheng Tan, Yilong Wang, Wenjie Xi, Xinping Yang, Hao Y. Zhang, and Jiaheng Zhao for their stimulating discussions, which greatly accelerated the emergence of this work. The two authors are supported by NSFC (Grant No. 11971219) and by Guangdong Provincial Key Laboratory (Grant No. 2019B121203002) and by Guangdong Basic and Applied Basic Research Foundation (Grant No. 2020B1515120100). AB is additionally supported by the start-up funding from Hao Zheng at Beijing Institute of Mathematical Sciences and Applications. ZHZ is supported by grants from Beijing Institute of Mathematical Sciences and Applications and also by Wu Wen-Tsun Key Laboratory of Mathematics at USTC of Chinese Academy of Sciences. 

\section{The reconstruction theorem for weak Hopf algebras}\label{sec:reconstruct}
In this section, we introduce our main tool, the reconstruction theorem for weak Hopf algebras. The Reconstruction Theorem establishes a bijection between the following two sets:
\bnu
	\item the set $\mathtt{X}$ of isomorphism classes of finite-dimensional weak Hopf algebras.
	\item \label{item2:sec:reconstruct} the set $\mathtt{Y}$ of equivalence classes of pairs $(\Cc,\Fff)$ where $\Cc$ is a finite multi-tensor category and $\Fff\:\Cc\to\vect$ is a faithful exact separable Frobenius functor.
\enu

This theorem was due to Szlachányi \cite{Szlachanyi_2000, Szlachanyi_2004}. Other versions of the Reconstruction Theorem include \cite{Hayashi_1999, Bohm_Caenepeel_Janssen_2011, McCurdy_2012, Vercruysse_2013}; for modern viewpoints, we refer the reader to \cite{Bruguieres_Virelizier_2007} and \cite{Bohm_Lack_Street_2011}. Our version, which is formulated as the bijection mentioned above, differs slightly from those found in the references, for which we choose to include a proof for completeness.

We emphasize that our main examples of weak Hopf algebras can be reconstructed using the procedure dictated by Hayashi \cite{Hayashi_1999}, which is a special case of the Reconstruction Theorem by \cite{Szlachanyi_2000, Szlachanyi_2004}. In particular, they're all \emph{face algebras} in the sense of \cite{Hayashi_1999}. However, we choose to work with the language of weak Hopf algebras and the general reconstruction theorem \cite{Szlachanyi_2000, Szlachanyi_2004}.

In \Cref{subsec:sf}, we recall the notion of separable Frobenius algebras and separable Frobenius functors. In \Cref{subsec:wha_to_wff}, we recall definitions and properties of weak Hopf algebras, and construct a representative of an element in $\mathtt{Y}$ out of a weak Hopf algebra. In \Cref{subsec:wff_to_wha}, we construct a weak Hopf algebra out of a representative of an element in $\mathtt{Y}$, and show that the two constructions induce a bijection between $\mathtt{X}$ and $\mathtt{Y}$.

\subsection{Separable Frobenius algebras and separable Frobenius functors}\label{subsec:sf}
For vector spaces $V$ and $W$ and a vector $u\in V\os W$, we sometimes use the notation $u=u^{(1)}\os u^{(2)}$ for convenience.
\begin{dfn}
	A \newdef{Frobenius algebra over $\kb$} is an algebra $(A,m\:A\os A\to A,1\in A)$ equipped with a coalgebra structure $(A,s\:A\to A\os A,\delta\:A\to\kb)$ such that $s$ is an $A$-$A$-bimodule map, i.e., satisfies
	\[
		s(x)(1\os y)=s(xy)=(x\os 1)s(y),\;\forall x,y\in A\,.
	\]
	
	A \newdef{separable Frobenius algebra} is a Frobenius algebra $(A,m,1,s,\delta)$ such that $m\circ s=\id_A$.
\end{dfn}
\begin{rmk}
	Many authors use the term ``special Frobenius algebras'' for what are called separable Frobenius algebras here.
\end{rmk}
In a separable Frobenius algebra $A=(A,m,1,s,\delta)$, we call $p\defdtobe s(1)$ the \newdef{(canonical) separability idempotent} of $A$. 
\begin{lem}\label{lem:can_idp}
	The separability idempotent $p$ satisfies the following conditions:
\bnu
	\item \label{item1:lem:can_idp} For any $x\in A$, there is $xp^{(1)}\os p^{(2)}=p^{(1)}\os p^{(2)}x$.
	\item \label{item2:lem:can_idp} $p^{(1)}p^{(2)}=1$.
	\item \label{item3:lem:can_idp} $p$ is an idempotent in $A\os A^\op$, i.e., we have $p^{(1)}p^{(1')}\os p^{(2')}p^{(2)}=p^{(1)}\os p^{(2)}$.
\enu
\end{lem}
\begin{rmk}
	A separable Frobenius structure on an algebra is uniquely determined by the separability idempotent and the counit.
\end{rmk}
Given any $k$-algebra $A$, the category $\BiMod(A|A)$ of $A$-$A$-bimodules and bimodule maps is a monoidal category under relative tensor product over $A$. When $A$ is a separable Frobenius algebra, this relative tensor product can be explicitly computed.
\begin{cor}\label{cor:rel_over_sep}
	Let $A$ be a separable Frobenius algebra with separability idempotent $p$. Let $(V,\rho)$ be a right $A$-module and $(W,\lambda)$ be a left $A$-module. Then the relative tensor product $V\os_AW$ is given by the retract of the idempotent
	\[
		V\os W\to V\os W, \quad (v\os w)\longmapsto v.p^{(1)}\os p^{(2)}.w\,.
	\]
	\pf
		We denote the idempotent by $e$. Note that $e$ coequalizes the diagram 
		\eqn{\label{eq:cor:rel_over_sep}
			\diagram{V\os A\os W \ar@<0.5ex>[r]^-{\rho\os\id} \ar@<-0.5ex>[r]_-{\id\os\lambda}& V\os W}
		}
		by an application of \Cref{lem:can_idp}.\ref{item1:lem:can_idp} . Suppose $e$ has a retraction $r$ with the associated section $i$. Then $r$ also coequalizes the diagram, since $i$ is an injection. It suffices to show that for any vector space $Q$ and map $q\:V\os W\to Q$ which coequalizes \eqref{eq:cor:rel_over_sep}, $q\circ i$ is the unique map satisfying $q\circ i\circ r=q$. This follows from the fact that $r\circ i=\id$ and \Cref{lem:can_idp}.\ref{item2:lem:can_idp}.	\epf
\end{cor}

The forgetful functor $\BiMod(A|A)\to\vect$, which sends each bimodule to its underlying vector space, naturally carries the structure of \emph{separable Frobenius functor}, a key ingredient in the reconstruction theorem for weak Hopf algebras. We now introduce this notion.
\begin{dfn}[\textnormal{\cite[Definition 1.7]{Szlachanyi_2000}}]\label{dfn:sf}
	Let $\Cc=(\Cc,\os,\mone)$ and $\Dd=(\Dd,\os',\mone')$ be monoidal categories, which we may assume to be strict. A \newdef{Frobenius (monoidal) functor from $\Cc$ to $\Dd$} consists of the following data:
	\bit
		\item A functor $F\:\Cc\to\Dd$.
		\item A lax monoidal functor structure $(F,F_2,F_0)$ on $F$. This means a family of morphisms
		\[
			\{\diagram@C=2.75pc{F(X)\os' F(Y)\ar[r]^-{{F_2}_{X,Y}} & F(X\os Y)}\}_{X,Y\in\Cc}
		\]
		natural in $X$ and $Y$ and a morphism 
		\[
			\diagram@C=2.5pc{\mone' \ar[r]^-{F_0} & F(\mone)}
		\]
		satisfying that for any $X,Y,Z\in\Cc$, the following diagrams are commutative:
		\[
			\diagram@=2.75pc{
				F(X)\os' F(Y)\os' F(Z) \ar[d]_{1\os' {F_2}_{Y,Z}} \ar[r]^-{{F_2}_{X,Y}\os' 1} & F(X\os Y) \os' F(Z) \ar[d]^{{F_2}_{X\os Y,Z}}\\
				F(X)\os' F(Y\os Z) \ar[r]_-{{F_2}_{X,Y\os Z}} & F(X\os Y\os Z)
			}
		\]
		\[
			\diagram{\mone'\os' F(X)  \ar[d]_1 \ar[r]^-{F_0\os' 1} & F(\mone)\os' F(X) \ar[d]^{{F_2}_{\mone,X}} \\
				F(X) & F(\mone\os X) \ar[l]^-1
			}
			\qquad
			\diagram{
				F(X)\os'\mone' \ar[d]_-1 \ar[r]^-{1\os' F_0} & F(X) \os' F(\mone) \ar[d]^{{F_2}_{X,\mone}} \\
				F(X) & F(X\os \mone) \ar[l]^-1
			}\,.
		\]
		\item An oplax monoidal structure $(F,F_{-2},F_{-0})$ on $F$. This means a family of morphisms
		\[
			\{\diagram@C=3pc{F(X\os Y)\ar[r]^-{{F_{-2}}_{X,Y}} & F(X)\os' F(Y)}\}_{X,Y\in\Cc}
		\]
		natural in $X$ and $Y$ and a morphism 
		\[
			\diagram@C=2.5pc{F(\mone)\ar[r]^-{F_{-0}}  & \mone'}
		\]
		such that $(F,F_{-2},F_{-0})$ form a lax monoidal functor structure on the functor $F\:\Cc^\op\to\Dd^\op$.
	\eit
	They're required to render the following two diagrams commutative for any $X,Y,Z\in\Cc$:
		\eqn{\label{diagram:frob_1}
			\diagram@=2.6pc{
				F(X\os Y)\os' F(Z) \ar[d]_{{F_{-2}}_{X,Y}\os' 1} \ar[r]^-{{F_2}_{X\os Y,Z}} & F(X\os Y\os Z) \ar[d]^{{F_{-2}}_{X,Y\os Z}} \\
				F(X)\os' F(Y)\os' F(Z) \ar[r]_-{1\os' {F_2}_{Y,Z}} & F(X)\os' F(Y\os Z)
			}
		}
		\eqn{\label{diagram:frob_2}
			\diagram@=2.6pc{
				F(X)\os' F(Y\os Z)\ar[d]_{1\os' {F_{-2}}_{Y,Z}} \ar[r]^-{{F_2}_{X,Y\os Z}} & F(X\os Y\os Z) \ar[d]^{{F_{-2}}_{X\os Y,Z}}\\
				F(X)\os' F(Y)\os' F(Z) \ar[r]_-{{F_2}_{X,Y}\os' 1} & F(X\os Y)\os' F(Z)
			}\,.
		}
		
	A Frobenius functor $(F,F_2,F_0,F_{-2},F_{-0})\:\Cc\to\Dd$ is said to be \newdef{separable} if the following \newdef{separability condition} holds: for any $X,Y\in\Cc$, there is
		\eqn{\label{eq:s:dfn:sf}
			{F_2}_{X,Y}\circ {F_{-2}}_{X,Y}=\id_{F(X\os Y)}\,.
		}
\end{dfn}
\begin{dfn}\label{dfn:iso_spf}
	Let  $F=(F,F_2,F_0,F_{-2},F_{-0})$ and $G=(G,G_2,G_0,G_{-2},G_{-0})$ be two separable Frobenius functors from $\Cc$ to $\Dd$. A \newdef{natural isomorphism of separable Frobenius functors $F\funto G$} is a natural isomorphism $\xi\:F\funto G$ such that the equalities
	\agn{
		\xi_{X\os Y}\circ {F_2}_{X,Y}&={G_2}_{X,Y}\circ (\xi_X\os'\xi_Y) \qquad \xi_\mone \circ F_0 =G_0 \\
		(\xi_X\os'\xi_Y)\circ {F_{-2}}_{X,Y}&={G_{-2}}_{X,Y}\circ \xi_{X\os Y} \qquad  \qquad \manualformatting\hspace{0.8pc}
		 F_{-0} =G_{-0}\circ \xi_{\mone}
	}
	hold for any $X,Y\in\Cc$. We say that $F$ and $G$ are \newdef{isomorphic} if there exists a natural isomorphism of separable Frobenius functors $F\funto G$.
\end{dfn}
\begin{expl}
	Any strong monoidal functor is a separable Frobenius functor.
\end{expl}
\begin{expl}
	The composition of two separable Frobenius functors has a natural structure of separable Frobenius functor.
\end{expl}
\begin{expl}[\textnormal{\cite[Lemma 6.4]{Szlachanyi_2004}}]\label{expl:bimod_forget}
	Let $A$ be a separable Frobenius algebra with separability idempotent $p$. Then the forgetful functor $\Uuu\:\BiMod(A|A)\to\vect$ is naturally a separable Frobenius functor  :
	\bnu
		\item For $V,W\in\BiMod(A|A)$, by \Cref{cor:rel_over_sep}, the space $\Uuu(V\os_AW)$ is given by the retract of the idempotent $e_{V,W}\:V\os W\to V\os W,\;v\os w\longmapsto v.p^{(1)}\os p^{(2)}.w$. Then we define ${\Uuu_2}_{V,W}$ to be the retraction of $e_{V,W}$.
		\item $\Uuu_0\:\kb\to A$ is defined as the unit of $A$.
		\item For $V,W\in\BiMod(A|A)$, ${\Uuu_{-2}}_{V,W}$ is defined as the section of $e_{V,W}$ associated with ${\Uuu_2}_{V,W}$.
		\item $\Uuu_{-0}\:A\to\kb$ is defined as the counit of $A$. 
	\enu
\end{expl}
\begin{rmk}\label{rmk:sep_fun_can_be_avoid}
	We're only concerned with separable Frobenius functors to $\vect$ in this work. It is shown in \cite[Lemma 6.2]{Szlachanyi_2004} that if $\Fff\:\Cc\to\vect$ is a separable Frobenius functor, then $\Fff(\mone)$ is a separable Frobenius algebra, and $\Fff$ factors as 
	\[\diagram{\Cc\ar[r]^-{F} & \BiMod(F(\mone)|F(\mone))\ar[r]^-{\Uuu} &\vect}\,,\]
	where $F$ is a strong monoidal functor and $\Uuu$ is the separable Frobenius functor defined in \Cref{expl:bimod_forget} . This factorization shows that the language of separable Frobenius functor can be avoided by working solely with the notions of strong monoidal functors and separable Frobenius algebras. However, following the practice of \cite{Szlachanyi_2000, Szlachanyi_2004}, we choose to stick to this language, as it provides convenience for both the statement and the proof of the Reconstruction Theorem.
\end{rmk}

For an object $X$ in a monoidal category, we use $X^L$ and $X^R$ to denote the left dual and right dual of $X$, respectively. We end this section with the following basic observation that Frobenius functors preserve duals.
\begin{lem}[\textnormal{\cite[Theorem 2]{Day_Pastro_2008}}]\label{lem:frob_preserve_dual}
	Let $F=(F,F_2,F_0,F_{-2},F_{-0})\:\Cc\to\Dd$ be a Frobenius functor. Suppose $(X^L,\ev\:X^L\os X\to\mone,\coev\:\mone\to X\os X^L)$ is a left dual of $X\in\Cc$. Then $(F(X^L),\Ev,\Coev)$ is a left dual of $F(X)$ in $\Dd$ with
\eqnn{
	\begin{split}
		\Ev=&(\diagram@C=3.35pc{F(X^L)\os'F(X) \ar[r]^-{{F_2}_{X^L,X}} & F(X^L\os X) \ar[r]^-{F(\ev)} & F\mone \ar[r]^-{F_{-0}} & \mone'}) \\ 
		\Coev=&(\diagram@C=3.35pc{\mone' \ar[r]^-{F_0} & F(\mone) \ar[r]^-{F(\coev)} & F(X\os X^L) \ar[r]^-{{F_{-2}}_{X,X^L}} & F(X)\os' F(X^L)})\,.  
	\end{split}
	}
\end{lem}

\subsection{From weak Hopf algebras to weak fiber functors}\label{subsec:wha_to_wff}
\begin{dfn}\label{dfn:wff}
	A \newdef{weak fiber functor} on a finite multi-tensor category $\Cc$ is a faithful and exact separable Frobenius functor from $\Cc$ to $\vect$.
\end{dfn}
In this subsection, we recall basic definition, examples and properties of weak Hopf algebras. Then we show how to construct a finite multi-tensor category together with a weak fiber functor on it from a weak Hopf algebra.

In \Cref{dfn:wff}, we adopt the terminologies from \cite{Etingof_Gelaki_Nikshych_Ostrik_2015}: a \newdef{finite multi-tensor category} is a finite $\kb$-linear rigid monoidal category such that the tensor product is bi-$\kb$-linear; a $\kb$-linear category is \newdef{finite} if it is equivalent to the category $\rep(B)$ of finite-dimensional left modules over a finite-dimensional algebra $B$.

\begin{dfn}\label{dfn:wba}
	A \newdef{weak bialgebra} over $\kb$ is an algebra $(A,\mu\:A\os A\to A,\eta\:\kb\to A)$ equipped with a coalgebra structure $(A,\Delta\:A\to A\os A,\epsilon\:A\to\kb)$ satisfying the following constraints:
	\bnu[(\textbf{Axiom }1)]
		\item \label{item:Axiom_1}The comultiplication is multiplicative: 
		\[\Delta(x)\Delta(y)  =\Delta(xy),\;\forall x,y\in A\,.\]
		\item \label{item:Axiom_2}The counit satisfies
		\eqn{\label{eq:Axiom_2}
		\epsilon(xy_{(1)})\epsilon(y_{(2)}z) =\epsilon(xyz)=\epsilon(xy_{(2)})\epsilon(y_{(1)}z),\;\forall x,y,z\in A\,.
		}
		\item \label{item:Axiom_3} The unit $1\defdtobe 
		\eta(1)$ satisfies
		\eqn{\label{eq:Axiom_3}
		1_{(1)}\os 1_{(2)}1_{(1')}\os 1_{(2')}  =1_{(1)}\os 1_{(2)}\os 1_{(3)}=1_{(1)}\os 1_{(1')}1_{(2)}\os 1_{(2')}\,.
		}
	\enu
	Here we use the Sweedler's notation $\Delta(x)=x_{(1)}\os x_{(2)}$.

	A \newdef{weak Hopf algebra over $\kb$} is a weak bialgebra $(A,\mu,\eta,\Delta,\epsilon)$ equipped with a linear map $S\:A\to A$ satisfying the following condition:
	\bnu[(\textbf{Axiom }1)]
		\setcounter{enumi}{3}
		\item For any $x\in A$, we have 
		\agnn{\label{eq1:antipode}
			x_{(1)}S(x_{(2)})&=\epsilon(1_{(1)}x)1_{(2)}\,;
		\\
		\label{eq2:antipode}
			S(x_{(1)})x_{(2)}&=1_{(1)}\epsilon(x1_{(2)})\,;
		\\
		\label{eq3:antipode}
			S(x_{(1)})x_{(2)}&S(x_{(3)})=S(x)\,.
		}
	\enu
	The map $S$ is called an \newdef{antipode}.
	
	Given two weak bialgebras $(A,\mu,\eta,\Delta,\epsilon)$ and $(B,\mu',\eta',\Delta',\epsilon')$, a \newdef{homomorphism of weak bialgebras $A\to B$} is a linear map $\phi\:A\to B$ such that 
	\agn{
		\mu'\circ (\phi\os \phi) &=\phi\circ \mu\qquad \manualformatting \hspace{3.9pc} \eta'=\phi\circ \eta\,;
	\\
		\Delta'\circ \phi &=(\phi\os \phi)\circ \Delta\qquad \epsilon'\circ \phi=\epsilon\,.
	}
	A \newdef{homomorphism of weak Hopf algebras} is a homomorphism of the underlying weak bialgebras.
\end{dfn}
The basic theory of weak Hopf algebras needed in this article is developed in \cite{Nill_1998, Bohm_Nill_Szlachanyi_1999, Szlachanyi_2000}.
\begin{rmk}\label{rmk:antipode}
	\bnu
		\item An antipode on a weak bialgebra, if exists, is unique. Moreover, a homomorphism of weak Hopf algebras necessarily preserves the antipode. 
		\item The antipode $S\:A\to A$ of a weak Hopf algebra $A$ must be an algebra anti-homomorphism and a coalgebra anti-homomorphism.
		\item The antipode of a finite-dimensional weak Hopf algebra is always invertible.
		\item \label{item4:rmk:antipode} Let $(A,S)$ be a finite-dimensional weak Hopf algebra. Then both $(A^\op,S\inverse)$ and $(A^\cop,S\inverse)$ are weak Hopf algebras, where $A^\op$ and $A^\cop$ are respectively the weak bialgebra obtained by reversing the multiplication and comultiplication of $A$.
		\enu
\end{rmk}
Only finite-dimensional weak Hopf algebras are considered in this work. The following examples of weak Hopf algebras, while not directly related to our main example in \Cref{sec:2}, are presented for pedagogical purposes.
\begin{expl}[Groupoid algebra]\label{expl:grpd}
	Let $\Gg$ be a groupoid with a finite set of morphisms $\Mor(\Gg)$. Then there exists a weak Hopf algebra structure on the vector space $\kb[\Gg]\defdtobe\mathrm{span}\{g|g\in\Mor(\Gg)\}$ defined as follows:
	\[
		\mu(h\os g)=
		\begin{cases}
			h\circ g,\phantom{0}\quad\text{if $b=c$;} \\
			0,\phantom{hg}\quad\text{otherwise,}
		\end{cases}
		\quad\forall (\diagram{c\ar[r]^-h & d}),(\diagram{a\ar[r]^-g& b}) \in\Mor(\Gg)\,;
	\]
	\[
		\eta(1)=\sum_{a\in\Gg}\id_a\,;
	\]
	\[
		\Delta(g)=g\os g,\quad \epsilon(g)=1,\quad S(g)=g\inverse, \quad\forall (\diagram{a\ar[r]^-g& b})\in\Mor(\Gg)\,.
	\]
\end{expl}
\begin{expl}[\textnormal{\cite[Appendix A]{Bohm_Nill_Szlachanyi_1999}}]\label{expl:BBop}
	Let $B=(B,s\:B\to B\os B,\delta\:B\to \kb)$ be a separable Frobenius algebra with separability idempotent $p=s(1)$. Then, there exists a weak Hopf algebra structure on the algebra $B\os B^\op$ defined as follows:
	\agn{
		\Delta\: a\os b & \longmapsto a\os p^{(1)}\os p^{(2)}\os b \\
		\epsilon\:a\os b& \longmapsto \delta(ab) \\
		S\: a\os b& \longmapsto b\os \tau(a)\,,
	}
	where $\tau\:B\to B,\;a\longmapsto \delta(ap^{(2)})p^{(1)}$ is the Nakayama automorphism of $B$.
\end{expl}
Given a weak bialgebra, it is instructive to define two idempotent maps:
\agn{
	\epsilon^{lr}\:& A\to A,\quad x\longmapsto \epsilon(1_{(1)}x)1_{(2)}\,; \\
	\epsilon^{rr}\:&A\to A,\quad x\longmapsto 1_{(1)}\epsilon(1_{(2)}x)\,. 
}
We denote $A^l\defdtobe\epsilon^{lr}(A)$ and $A^r\defdtobe\epsilon^{rr}(A)$. 

\begin{thm}[\textnormal{\cite{Nill_1998, Bohm_Nill_Szlachanyi_1999}}]\label{thm:comm_subalg} 
Let $A$ be a weak bialgebra.
	\bnu
			\item \label{item1:thm:comm_subalg} $A^l$ and $A^r$ are unital subalgebras of $A$.
			\item \label{item2:thm:comm_subalg} $A^l$ and $A^r$ mutually commute in $A$, i.e., for $x\in A^l$ and $y\in A^r$, we have $xy=yx$.
			\item \label{item3:thm:comm_subalg} $\epsilon^{lr}\mid_{A^r}\:A^r\to A^l$ and $\epsilon^{rr}\mid_{A^l}\:A^l\to A^r$ are mutually inverse algebra anti-isomorphisms.
			\item \label{item4:thm:comm_subalg} The algebra $A^l$ has a structure of separable Frobenius algebra, with the separability idempotent given by $p\defdtobe (\epsilon^{lr}\os\id)\Delta(1)\in A^l\os A^l$ and the counit given by $\epsilon\mid_{A^l}\:A^l\to \kb$.
		\enu
	\pf
		1. It is established by \cite[Proposition 2.6]{Nill_1998}. 2. It follows from an easy application of (\textbf{Axiom }\ref{item:Axiom_3}). 3. It follows from \cite[Corollary 3.6]{Nill_1998}. 4. It is shown in \cite[Proposition 5.2]{Nill_1998} that $A^l$ is a separable algebra with separability idempotent $p$. It is easy to conclude that $A^l$ is in fact a separable Frobenius algebra, with separability idempotent $p$ and the counit defined above; see also \cite[Proposition 4.2]{Schauenburg_2003}.
	\epf
\end{thm}
The separable Frobenius algebras $A^l,A^r$ are called \newdef{base algebras} of $A$ and are crucial in the theory of weak Hopf algebras.
\begin{prp}[\textnormal{\cite{Nill_1998, Szlachanyi_2000, Bohm_Szlachanyi_2000}}]\label{prp:repA}
	The category $\rep(A)$ of left $A$-modules over a weak Hopf algebra $A$ is a finite multi-tensor category.
\end{prp}
We do not provide a full proof of \Cref{prp:repA} here; instead, we focus on presenting the finite multi-tensor category structure on $\rep(A)$. We first give the rigid monoidal category structure on $\rep(A)$ in the following steps:
\bnu
	\item Given $V,W\in\rep(A)$, we need to define their tensor product $V\bos W$. We define the underlying vector space of $V\bos W$ as the retract of the idempotent 
\eqn{\label{eq:act_by_1_1_and_1_2}
	e_{V,W}\:V\os W\to V\os W,\quad v\os w\longmapsto 1_{(1)}.v\os 1_{(2)}.w\,.
}
It is convenient to identify this retract as $\im(e_{V,W})\defdtobe\{u\in V\os W\mid e_{V,W}(u)=u\}\subset V\os W$. The action of $x\in A$ on $V\bos W$ can then be given by the restriction of the map
\[
	V\os W\to V\os W,\quad v\os w\longmapsto x_{(1)}.v\os x_{(2)}.w
\]
on $V\bos W$. 
	\item The tensor unit $\mone\in\rep(A)$ is given by the space $A^l$ endowed with the following left $A$-module action:
\[
	A\os A^l\to A^l,\quad x\os y\longmapsto \epsilon^{lr}(xy)\,.
\]
Equivalently, it can be given by $A^r$ endowed with the action $A\os A^r\to A^r,\; x\os y\longmapsto \epsilon^{rr}(xy)$; the two representations are isomorphic via the two linear maps in \Cref{thm:comm_subalg}.\ref{item3:thm:comm_subalg}.
	\item To define the associator, notice that for $V,W,U\in\rep(A)$, both the spaces $(V\bos W)\bos U$ and $V\bos (W\bos U)$ are given by the retract of the idempotent 
\[
	V\os W\os U\to V\os W\os U,\quad v\os w\os u\longmapsto 1_{(1)}.v\os 1_{(2)}.w\os 1_{(3)}.u\,,
\]
which can be verified using (\textbf{Axiom} \ref{item:Axiom_1}).
We then define the associator $a_{V,W,U}$ as the canonical map between the two retracts. One can check that $a_{V,W,U}$ is an $A$-module map, and satisfies the pentagon equation.
	\item Given $V\in\rep(A)$, we set the left unitor $l_V\:A^l\bos V\to V$ to be the restriction of the map
\[
	A^l\os V\to V,\quad x\os v\longmapsto x.v
\]
on $A^l\bos V$; we set the right unitor $r_V\:V\bos A^l\to V$ to be the restriction of the map
\[
	V\os A^l \to V,\quad v\os y\longmapsto \epsilon^{rr}(y).v
\]
on $V\bos A^l$. One can check that $l_V$ and $r_V$ are indeed invertible $A$-module maps, and satisfy the triangle equations.
	\item The left dual $V^L$ of an object $V\in\rep(A)$ is given by the dual vector space $V^\ast\defdtobe\Hom(V,\kb)$ endowed with the $A$-action
\[
	x.\omega=\omega(S(x).-),\; \forall\omega\in V^\ast, x\in A\,.
\]
Similarly, the right dual $V^R$ of $V$ is given by $V^\ast$ endowed with $A$-action
\[
	x.\omega=\omega(S\inverse(x).-),\; \forall\omega\in V^\ast, x\in A\,.
\]
\enu

Secondly, note that $\rep(A)$ is clearly a finite $\kb$-linear category, with the tensor product $\bos$ being bi-$\kb$-linear. This concludes our construction of $\rep(A)$.
\begin{expl}\label{expl:illust}
	We illustrate the above construction of $\rep(A)$ when $A$ is the weak Hopf algebra $B\os B^\op$ defined in \Cref{expl:BBop}. Since a left $A$-module is precisely a $B$-$B$-bimodule, $\rep(A)$ is equivalent to $\BiMod(B|B)$ as categories. It remains to find the monoidal structure on $\rep(A)$. Given left $A$-modules $V$ and $W$, which we identify as $B$-$B$-bimodules, the underlying vector space of $V\bos W$ is the retract of the idempotent
	\[
		V\os W\to V\os W,\quad v\os w\longmapsto v.p^{(1)}\os p^{(2)}.w\,.
	\]
	The action of $a\os b\in B\os B^\op$ on $V\bos W$ is given by the restriction of the map
	\[
		V\os W\to V\os W,\quad v\os w\longmapsto a.v.p^{(1)}\os p^{(2)}.w.b
	\]
	on $V\bos W$.	Using \Cref{cor:rel_over_sep}, it can be shown that the $B$-$B$-bimodule $V\bos W$ is precisely $V\os_BW$.
	
	To find the tensor unit of $\rep(B\os B^\op)$, one first computes 
	\[\epsilon^{lr}\:B\os B^\op\to B\os B^\op,\quad a\os b\longmapsto ab\os 1\,;
	\]
	\[\epsilon^{rr}\: B\os B^\op\to B\os B^\op,\quad a\os b\longmapsto 1\os ab\,.\]
	Therefore, $(B\os B^\op)^l=B\os 1$ and $(B\os B^\op)^r=1\os B^\op$. The tensor unit is hence $B\os 1$ with the $B\os B^\op$-action 
	\[(B\os B^\op)\os (B\os 1)\to (B\os 1),\quad a\os b\os c\os 1\longmapsto acb\os 1\,,\]
	or equivalently $1\os B$ with action $(B\os B^\op)\os (1\os B)\to (1\os B),\; a\os b\os 1\os c\longmapsto 1\os acb$.
	This shows that the tensor unit is isomorphic to the regular $B$-$B$-bimodule $B$. 
	
	With some additional efforts, one can show that $\rep(A)$ is equivalent, as a finite multi-tensor category, to the monoidal category $\BiMod(B|B)$, whose tensor product is given by the relative tensor product.
\end{expl}
\begin{rmk}
	\Cref{expl:illust} in particular shows that there exists infinitely-many weak Hopf algebras $A$ such that $\rep(A)\cong\vect$ as finite multi-tensor categories. Namely, for any $n\geq 1$, one can take $A=M_n(\kb)\os M_n(\kb)^\op$, where $M_n(\kb)$ is the algebra of $n\times n$-matrices equipped with the canonical symmetric separable Frobenius algebra structure. This echos with the fact that there are infinitely many weak Hopf algebras that be ``reconstructed'' from a fusion category, as will be clear in \Cref{rmk:to_funMM}.
\end{rmk}
Our next step is to construct a weak fiber functor $\rep(A)\to\vect$.

By \ref{item1:thm:comm_subalg}-\ref{item3:thm:comm_subalg} of \Cref{thm:comm_subalg}, there is an algebra homomorphism
\[
	\kappa\:A^l\os(A^l)^\op\to A, \quad x\os y\longmapsto x\epsilon^{rr}(y)\,,
\]
which induces a ``change of scalars'' functor:
\[
	F^A\:\rep(A)\to \BiMod(A^l|A^l),\quad {}_AV\longmapsto {}_\kappa V\,.
\]
Explicitly, for a left $A$-module $V$, the left $A^l$-action on ${}_\kappa V$ is given by $x\os v\longmapsto x.v$ while the right $A^l$-action is given by $v\os y\longmapsto \epsilon^{rr}(y).v$. 
\begin{thm}[\cite{Szlachanyi_2000}]\label{thm:forget_to_bimod}
	The functor $F^A$ is a faithful and exact strong monoidal functor.
	\pf
		The faithfulness and exactness come from the fact that $F^A$ is induced by an algebra homomorphism. It remains to show that $F^A$ is a strong monoidal functor. Let $V,W\in\rep(A)$. By \Cref{thm:comm_subalg}.\ref{item4:thm:comm_subalg} and \Cref{cor:rel_over_sep}, the space $F^A(V)\os_{A^l}F^A(W)$ is given by the retract of $g\:V\os W\to V\os W,\;v\os w\longmapsto \epsilon^{rr}\epsilon^{lr}(1_{(1)}).v\os 1_{(2)}.w$. However, one can check that $\epsilon^{rr}\epsilon^{lr}(1_{(1)})\os 1_{(2)}=1_{(1)}\os 1_{(2)}$ using (\textbf{Axiom }\ref{item:Axiom_2}) and (\textbf{Axiom }\ref{item:Axiom_3}), hence $g=e_{V,W}$, with $e_{V,W}$ defined by \eqref{eq:act_by_1_1_and_1_2}. Thus, by definition of $V\bos W$, we have a canonical isomorphism ${F^A_2}_{V,W}\:F^A(V)\os_{A^l}F^A(W)\isom F^A(V\bos W)$ of vector spaces. We leave it to the reader to check that ${F^A_2}_{V,W}$ is an isomorphism of bimodules. One can also show that $F^A(\mone)$ is precisely the regular $A^l$-$A^l$-bimodule ${}_{A^l}{A^l}_{A^l}$, for which we can take $F^A_0\:{}_{A^l}{A^l}_{A^l}\to F^A(\mone)$ to be the identity map. We also leave it to the reader to verify that $(F^A,F^A_2,F^A_0)$ form a strong monoidal functor. 
	\epf
\end{thm}
\begin{cor}\label{cor:wha_to_wff}
	Let $A$ be a weak Hopf algebra. Then the forgetful functor $\Fff^A\:\rep(A)\to\vect$ has a structure of weak fiber functor. Its separable Frobenius functor structure is given as follows:
	\bnu
		\item For $V,W\in\rep(A)$, ${\Fff^A_2}_{V,W}\:V\os W\to V\bos W$ is given by a retraction of $e_{V,W}\:V\os W\to V\os W,\;v\os w\longmapsto 1_{(1)}.v\os 1_{(2)}.w$.
		\item $\Fff^A_0\:\kb\to A^l$ is given by $1\longmapsto 1_A$.
		\item For $V,W\in\rep(A)$, ${\Fff^A_{-2}}_{V,W}\:V\bos W\to V\os W$ is given by the section of $e_{V,W}$ associated with ${\Fff^A_2}_{V,W}$.
		\item $\Fff^A_{-0}\:A^l\to\kb$ is given by $\epsilon\mid_{A^l}$.
	\enu
	\pf
	It is clear that $\Fff^A$ is faithful and exact. The given separable Frobenius structure on $\Fff^A$ comes from viewing $\Fff^A$ as the composition of two separable Frobenius functors below:
	\[\diagram{
		\rep(A)\ar[r]^-{F^A} & \BiMod(A^l|A^l) \ar[r]^-{\Uuu} & \vect
	}\,,
	\]
	where $\Uuu$ is the separable Frobenius functor associated with $A^l$ defined in \Cref{expl:bimod_forget}.
	\epf
\end{cor}
\begin{rmk}
	In the case $A=B\os B^\op$, building on \Cref{expl:illust}, one can show that $F^A$ is nothing but the identity strong monoidal functor. Thus 
	$\Fff^A\:\BiMod(B|B)\to\vect$ 
	coincides with the separable Frobenius functor associated with $B$ defined in \Cref{expl:bimod_forget}.
\end{rmk}

\subsection{From weak fiber functors to weak Hopf algebras}\label{subsec:wff_to_wha}
In this section, we construct a weak Hopf algebra $A^\Fff$ from a pair
\[
	(\Dd,\diagram{\Dd\ar[r]^-{\Fff} & \vect})\,,
\]
where $\Dd$ is a finite multi-tensor category and $\Fff$ is a weak fiber functor on $\Dd$. Then we show that the construction $(\Dd,\Fff)\longmapsto A^\Fff$ is the inverse of the construction
\[
	A\longmapsto (\rep(A),\diagram{\rep(A)\ar[r]^-{\Fff^A} &\vect})
\]
we introduced in \Cref{subsec:wha_to_wff}. As we have stated, this result is due to \cite{Szlachanyi_2000, Szlachanyi_2004}.

Our first task is to see how to reconstruct an algebra from a (not necessarily monoidal) functor.

Let $\Aa$ be a finite $\kb$-linear category and $F\:\Aa\to\vect$ be a functor. Then the space $\End(F)$ of endo-natural transformations on $F$ is naturally a $k$-algebra, with multiplication given by the composition of natural transformations. Moreover, one can define a \newdef{comparison functor}
\eqn{\label{eq:comparison}
\begin{split}
	\widetilde{F}\:\Aa & \to\rep(\End(F)) \\
			X & \longmapsto F(X)\,,
\end{split}
}
where $F(X)$ is equipped with the following evident left $\End(F)$-action: 
\[\End(F)\os F(X) \to F(X),\quad \alpha\os v\longmapsto \alpha_X(v)\,.\]

 Then one immediately has the following strictly commutative diagram of functors, where $U^{\End(F)}$ is the forgetful functor forgetting the $\End(F)$-action:
\[
	\diagram{
		\Aa \ar[r]^-F \ar[d]_{\widetilde{F}} & \vect \\
		\rep(\End(F)) \ar[ru]_-{\hskip 0.75em U^{\End(F)}} 
	}\,.
\]
We now present the conditions on $F$ such that $\widetilde{F}$ is an equivalence of categories. Of course, when $\widetilde{F}$ is an equivalence, the functor $F$ shares all properties with the forgetful functor $U^{\End(F)}$, hence is faithful and exact. The following well-known fact states that the converse is also true:
\begin{thm}\label{thm:monadicity}
	$\widetilde{F}$ is an equivalence if and only if $F$ is faithful and exact. 
\end{thm}
\Cref{thm:monadicity} can be proved by Beck's monadicity theorem. In this work, we do not prove it but instead treat it as a technical condition for our purposes. Note that if $F$ is exact, then $\End(F)$ is finite-dimensional\footnote{Since $\Fff$ is left exact, there exists an object $X\in\Aa$ which represents $F$. Then by Yoneda lemma, $\End(\Fff)\cong\End(X)^\op$ as vector spaces, which is finite-dimensional because $\Aa$ is finite.}.

The following lemma will also be of later use:
\begin{lem}[See for e.g. \textnormal{\cite[Proposition 1.8.15]{Etingof_Gelaki_Nikshych_Ostrik_2015}}]\label{lem:EndF_times_EndG}
	Let $\Aa$, $\Bb$ be finite $\kb$-linear categories, and let $F\:\Aa\to\vect, G\:\Bb\to\vect$ be faithful and exact functors. For $\alpha\in\End(F),\beta\in\End(G)$, define a natural transformation $J_{F,G}(\alpha\os\beta)\:\os(F\times G)\funto \os(F\times G)$, componentwise, by setting
	\[
		J_{F,G}(\alpha\os\beta)_{X,Y}\defdtobe\alpha_X\os\beta_Y, \;\forall X\in\Aa,Y\in\Bb\,.
	\]
	Then the map
	\[
		J_{F,G}\:\End(F)\os \End(G)\to \End(\os(F\times G)),\quad  \alpha\os\beta\longmapsto J_{F,G}(\alpha\os\beta)
	\]
	is a linear isomorphism.
	\pf
		The proof is given in \Cref{subsec:pf:lem:EndF_times_EndG}.
	\epf
\end{lem}

Now we are prepared to demonstrate the way to obtain a weak Hopf algebra $A^\Fff$ from a weak fiber functor $(\Dd,\Fff)$. We define $A^\Fff\defdtobe\End(\Fff)$ as an algebra. It remains to define a comultiplication, a counit, and an antipode on $\End(\Fff)$.

\paragraph{The comultiplication} We define the comultiplication $\Delta\:\End(\Fff)\to\End(\Fff)\os\End(\Fff)$ in two steps. First, we define a map
\[
	\underline{\Delta}\:\End(\Fff)\to\End(\os(\Fff\times\Fff))
\]
by setting $\underline{\Delta}(\alpha)_{X,Y}$ to be the map
\[\diagram@C=2.5pc{\Fff(X)\os \Fff(Y) \ar[r]^-{{\Fff_2}_{X,Y}} & \Fff(X\os Y) \ar[r]^-{\alpha_{X\os Y}} & \Fff(X\os Y) \ar[r]^-{{\Fff_{-2}}_{X,Y}} & \Fff(X)\os \Fff(Y)}\]
for $X,Y\in\Dd$ and $\alpha\in\End(\Fff)$. Secondly, we define
\[
	\Delta\defdtobe J_{\Fff,\Fff}\inverse\underline{\Delta}\:\End(\Fff)\to \End(\Fff)\os\End(\Fff)\,,
\]
where $J_{\Fff,\Fff}$ is defined in \Cref{lem:EndF_times_EndG}.

\paragraph{The counit} We define the counit $\epsilon\:\End(\Fff)\to\kb$ by setting $\epsilon(\alpha)$ to be the image of $1\in\kb$ under the map
\[
	\diagram@C=2.5pc{
		\kb \ar[r]^-{\Fff_0} & \Fff(\mone) \ar[r]^-{\alpha_\mone} & \Fff(\mone) \ar[r]^-{\Fff_{-0}} & \kb
	}
\]
for $\alpha\in\End(\Fff)$.

\paragraph{The antipode} Let $X^L$ be the left dual of $X\in\Dd$. Then by \Cref{lem:frob_preserve_dual}, $\Fff(X^L)$ is the left dual of $\Fff(X)$. We define the antipode $S\:\End(\Fff)\to\End(\Fff)$ by setting 
\eqn{
\label{eq:reconstruct_antipode}
S(\alpha)_X\defdtobe 
\begin{array}{l}
(\diagram@C=4pc{\Fff(X) \ar[r]^-{\Coev\os 1} & \Fff(X)\os \Fff(X^L)\os \Fff(X) } \\
\hskip -0.5em \diagram@C=4pc{\ar[r]^-{1\os\alpha_{X^L}\os 1} & \Fff(X)\os \Fff(X^L)\os \Fff(X) \ar[r]^-{1\os\Ev} & \Fff(X)})
\end{array}
}
for $X\in\Dd$ and $\alpha\in\End(\Fff)$, where $\Coev$ and $\Ev$ are respectively the unit and the counit witnessing $\Fff(X^L)$ as the left dual of $\Fff(X)$, given also in \Cref{lem:frob_preserve_dual}.
\begin{rmk}
	The weak Hopf algebra structure reconstructed from a weak fiber functor $\Fff\:\Dd\to\vect$ can also be defined on the ``end'' $\eend(\Fff)\defdtobe\int_{X\in\Dd}\Hom(\Fff(X),\Fff(X))$, which is well-known to be isomorphic to $\End(\Fff)$ as algebras. The comultiplication, the counit and the antipode on $\eend(\Fff)$ can be obtained using the structure maps of $\Fff$ and the Fubini theorem for ends. For an introduction to ends, we refer the reader to \cite[\S IX.5 and \S IX.8]{MacLane_1978}\cite{Loregian_2021}.
\end{rmk}

For a weak Hopf algebra $A$, recall from \Cref{cor:wha_to_wff} that the forgetful functor $\Fff^A\:\rep(A)\to\vect$ has a canonical weak fiber functor structure. 
\begin{thm}[\textnormal{\cite{Szlachanyi_2000, Szlachanyi_2004}} Reconstruction theorem for weak Hopf algebras, part I]\manualformatting\ 
\label{thm:reconstruct}
\bnu
	\item \label{item1:thm:reconstruct} $(\End(\Fff),\Delta,\epsilon,S)$ is a weak Hopf algebra\footnote{In the exposition of the reconstruction theorem for weak Hopf algebras in the textbook \cite{Etingof_Gelaki_Nikshych_Ostrik_2015}, the weak Hopf algebra structure appears to be constructed on $\End(F)$ for a faithful exact strong monoidal functor $F\:\Dd\to\BiMod(R|R)$, where $R$ is a separable Frobenius algebra \cite[Proposition 7.23.11]{Etingof_Gelaki_Nikshych_Ostrik_2015}. This seems inconsistent with the reconstruction theorem presented here, as we are essentially working with $\End(\Uuu F)$ rather than $\End(F)$, where $\Uuu$ denotes the forgeful functor from $\BiMod(R|R)$ to $\vect$. We believe that there is no meaningful weak Hopf algebra structure on $\End(F)$. On the other hand, $\End(F)$ in the sense of \cite[Proposition 7.23.11]{Etingof_Gelaki_Nikshych_Ostrik_2015} actually refers to $\End(\Uuu F)$, as suggested by the following explanation written  elsewhere by three of the authors of \cite{Etingof_Gelaki_Nikshych_Ostrik_2015}: ``Let $A =\End_\kb(F)$ (i.e. the algebra of endomorphisms of the composition of $F$ with the forgetful functor to vector spaces).''\cite[\S 2.5]{Etingof_Nikshych_Ostrik_2005}.}.
	\item \label{item2:thm:reconstruct} Let $A=(A,\mu',\eta',\Delta',\epsilon')$ be a weak Hopf algebra. Then $A\cong \End(\Fff^A)$ as weak Hopf algebras. 
	\item  \label{item3:thm:reconstruct}  Let $(\Dd,\Fff)$ be a weak fiber functor, and let $\End(\Fff)$ be the reconstructed weak Hopf algebra in \ref{item1:thm:reconstruct}. Then the comparison functor 
	\[
		\widetilde{\Fff}\:\Dd\to\rep(\End(\Fff))
	\]
	is a monoidal equivalence such that the following diagram of separable Frobenius functors strictly commutes:
	\[
	\diagram{
		\Dd \ar[r]^-\Fff \ar[d]_{\widetilde{\Fff}} & \vect \\
		\rep(\End(\Fff)) \ar[ru]_-{\hskip 0.75em \Fff^{\End(\Fff)}} 
	}\,.
	\]
\enu
\pf
	\bnu
		\item The proof is given in \Cref{subsec:pf:thm:reconstruct}.
		\item For $x\in A$, let $j(x)$ denote the natural transformation defined by $(j(x))_V\defdtobe x.-\:V\to V$ for $V\in\rep(A)$. Then the map
		\[
			j\:A\to \End(\Fff^A), \quad x\longmapsto j(x)\,,
		\]
		is clearly an algebra isomorphism. 
		
		To show that $j$ preserves the comultiplication, it is enough to verify that $j(x_{(1)})\os j(x_{(2)})=\Delta j(x)$ for $x\in A$, which is equivalent to
		\[
			J_{\Fff^A,\Fff^A}(j(x_{(1)})\os j(x_{(2)}))=J_{\Fff^A,\Fff^A}\Delta j(x)\equiv \underline{\Delta}(j(x))\,,
		\]
		with $J_{\Fff^A,\Fff^A}$ given by \Cref{lem:EndF_times_EndG}. This indeed holds true, since for any $V,W\in\rep(A)$, we have that the map $J_{\Fff^A,\Fff^A}(j(x_{(1)})\os j(x_{(2)}))_{V,W}$ reads
		\[
			\diagram@C=5.25pc{V\os W\ar[r]^-{x_{(1)}.- \os x_{(2)}.-} &  V\os W}
		\]
		while the map $\underline{\Delta}(j(x))_{V,W}$ reads
		\[
			\diagram@C=5.25pc{V\os W \ar[r]^-{1_{(1)}.- \os 1_{(2)}.-} & V\os W \ar[r]^-{x_{(1)}.- \os x_{(2)}.-} & V\os W \ar[r]^-{1_{(1)}.- \os 1_{(2)}.-} & V\os W}\,.
		\]
		The two maps are equal by (\textbf{Axiom }\ref{item:Axiom_1}) of weak bialgebras.
		
		Finally, $j$ preserves the counit since for any $x\in A$, we have that $\epsilon j(x)$ is the image of $1\in\kb$ under the map 
		\[
			\diagram@C=2.25pc{\kb\ar[r]^-{\eta'} & A^l \ar[r]^-{x.-} & A^l \ar[r]^{\epsilon'\mid_{A^l}} & \kb}\,,
		\]
		which reads
		\[
			\epsilon\epsilon^{lr}(x)=\epsilon(x)\,.
		\]
		
		\item Note that $\widetilde{\Fff}$ is an equivalence by \Cref{thm:monadicity}. To check that $\widetilde{\Fff}$ is a monoidal equivalence, we need only verify that $\widetilde{\Fff}$ is a strong monoidal functor. For $X,Y\in\Dd$, by definition, the underlying vector space of $\widetilde{\Fff}(X)\bos\widetilde{\Fff}(Y)$ is given the retract of the idempotent $\underline{\Delta}(\id_\Fff)$, which reads
		\[
			\diagram@C=2.5pc{\Fff(X)\os\Fff(Y) \ar[r]^-{{\Fff_2}_{X,Y}} & \Fff(X\os Y) \ar[r]^-{{\Fff_{-2}}_{X,Y}} & \Fff(X)\os\Fff(Y)}.
		\]
		This retract is manifestly $\Fff(X\os Y)$ by the separability condition \eqref{eq:s:dfn:sf} satisfied by $\Fff$. Then there exists a canonical isomorphism ${\widetilde{{\Fff}_2}}_{X,Y}\:\widetilde{\Fff}(X)\bos \widetilde{\Fff}(Y)\isom\widetilde{\Fff}(X\os Y)$ of vector spaces. We also define a linear isomorphism
		\[
			\widetilde{\Fff}_0\:\End(\Fff)^l\to \Fff(\mone)
		\]
		by sending $\gamma\in\End(\Fff)^l$ to the image of $1\in\kb$ under the map
		\[
			\diagram@C=2.25pc{\kb\ar[r]^-{\Fff_0} & \Fff(\mone)\ar[r]^-{\gamma_\mone} & \Fff(\mone)}\,.
		\]
		The inverse $\widetilde{\Fff}_0\inverse$ sends $y\in\Fff(\mone)$ to the natural transformation $\Fff\funto\Fff$ given by ``left multiplication by $y$''. To be precise, notice that the element $y$ can be identified as a map $\kb\to\Fff(\mone),\;1\longmapsto y$. Then we define $\widetilde{\Fff}_0\inverse(y)_X$ for $X\in\Dd$ by 
		\[
			\diagram@C=2.25pc{
				\Fff(X)\ar[r]^-1_-\sim & \kb\os\Fff(X) \ar[r]^-{y\os 1}  & \Fff(\mone)\os \Fff(X) \ar[r]^-{{\Fff_2}_{\mone,X}} & \Fff(\mone\os X)\ar[r]^-1_-\sim & \Fff(X)
			}\,.
		\]
		We leave it to the reader to verify that $\widetilde{\Fff}\defdtobe(\widetilde{\Fff},\{{\widetilde{{\Fff}_2}}_{X,Y}\}_{X,Y\in\Dd},\widetilde{\Fff}_0)$ is a strong monoidal functor. 
		
		By our construction of $\widetilde{\Fff}$, it is straightforward to see that $\Fff^{\End(\Fff)}\widetilde{\Fff}=\Fff$ as separable Frobenius functors.
	\enu
\epf
\end{thm}
From now on, we will abuse the notation by using  the term ``weak fiber functor'' to denote a pair consisting of a finite multi-tensor category and a weak fiber functor on it. 
\begin{dfn}\label{dfn:equiv_wff}
	Two weak fiber functors $(\Dd,\diagram{\Dd\ar[r]^-\Fff &\vect})$ and $(\Ee,\diagram{\Ee\ar[r]^-{\Ggg} &\vect})$ are \newdef{equivalent} if there exists a monoidal equivalence $\Phi\:\Dd\to\Ee$ such that $\Ggg\Phi$ and $\Fff$ are isomorphic as separable Frobenius functors.
\end{dfn}
One can show that \Cref{dfn:equiv_wff} indeed defines an equivalence relation on the set of all weak fiber functors.
\begin{thm}[\cite{Szlachanyi_2000, Szlachanyi_2004} Reconstruction theorem for weak Hopf algebras, part II]\label{thm:reconstruct_II}
	The assignment
	\[
		A \longmapsto (\rep(A),\diagram{\rep(A)\ar[r]^-{\Fff^A} & \vect})
	\]
	sends isomorphic weak Hopf algebras to equivalent weak fiber functors. 
	
	The assignment
	\[
		(\Dd,\Fff)\longmapsto \End(\Fff)
	\]
	sends equivalent weak fiber functors to isomorphic weak Hopf algebras.
	
	Consequently, by \ref{item2:thm:reconstruct} and \ref{item3:thm:reconstruct} of \Cref{thm:reconstruct}, these assignments establish mutually inverse bijections between the set of isomorphism classes of weak Hopf algebras and the set of equivalence classes of weak fiber functors. 
	\pf
		To show the first statement, let $A$, $B$ be weak Hopf algebras and $\phi\:A\to B$ be an isomorphism of weak Hopf algebras. Then, it can be verified that the ``change of scalars'' functor 
		\[
			\phi^\ast\:\rep(B)\to\rep(A),\quad {}_BV\longmapsto {}_\phi V
		\]
		is a monoidal equivalence such that $\Fff^A\phi^\ast=\Fff^B$ as separable Frobenius functors.
				
		To show the second statement, let $(\Dd,\Fff)$ and $(\Ee,\Ggg)$ be weak fiber functors, $\Phi\:\Dd\to\Ee$ be a monoidal equivalence, and $\xi\:\Ggg\Phi\funto\Fff$ be an isomorphism of separable Frobenius functors (see \Cref{dfn:iso_spf}) as illustrated in the diagram
		\[
			\diagram@=2.5pc{
				\Dd \ar[d]_\Phi \ar@{}[r]_(0.4){\raisebox{-1em}{$\scriptstyle\Uparrow \xi$}} \ar[r]^-{\Fff} & \vect \\
				\Ee \ar[ru]_-{\Ggg}
			}\,.
		\]
		Then one can verify that the map
		\[
			\End(\Ggg)\to\End(\Fff), \quad \alpha\longmapsto \xi\cdot(\alpha\Phi)\cdot\xi\inverse
		\]
		defines an isomorphism of weak Hopf algebras.
		
		The first two statements show that there are well-defined maps
		\[
			[A]\longmapsto[(\rep(A),\Fff^A)]\quad\text{and}\quad 	
[(\Dd,\Fff)]\longmapsto [\End(\Fff)] 	
		\]
		between the set of isomorphism classes of weak Hopf algebras and the set of equivalence classes of weak fiber functors. Then \ref{item2:thm:reconstruct} and \ref{item3:thm:reconstruct} of \Cref{thm:reconstruct} immediately imply that these two maps are mutually inverse.
	\epf
\end{thm}
\begin{rmk}\label{rmk:to_funMM}
	It has long been known that if a finite multi-tensor category $\Dd$ admits a faithful exact strong monoidal functor
	\eqnn{
		F\:\Dd\to\fun(\Mm,\Mm)\,,
	}
	where $\Mm$ is a finite semisimple category, then $\Dd\cong\rep(A)$ for a weak Hopf algebra $A$ \cite{Hayashi_1999}\cite[\S 4]{Ostrik_2003}. The functor $F$ arises in many circumstances: it appears precisely when $\Mm$ is a \emph{faithful} module category over $\Dd$ \cite[Definition 7.12.9]{Etingof_Gelaki_Nikshych_Ostrik_2015}, and when $\Dd$ is fusion, every non-zero module category is faithful. In this remark, we review the construction of $A$ using the Reconstruction Theorem (\Cref{thm:reconstruct,thm:reconstruct_II}), and discuss the uniqueness of the weak Hopf algebras constructed in this manner. 
	
	By \Cref{thm:reconstruct}, to construct $A$, it suffices to find a faithful exact separable Frobenius functor $\fun(\Mm,\Mm)\to\vect$. This can be done in two steps. First, choose an algebra $B$ such that $\rep(B)\cong\Mm$. Since $\Mm$ is semisimple, $B$ is necessarily semisimple, and it takes the form 
	\eqn{\label{eq:semisimple}
		\Oplus_{x\in\Irr(\Mm)} M_{n_x}(\kb)
	}
	for positive integers $\{n_x\}_{x\in\Irr(\Mm)}$. Along with the choice of $B$ we have a monoidal equivalence
	\[
		\Psi\:\diagram{\fun(\Mm,\Mm)\ar[r]^-\sim & \fun(\rep(B),\rep(B))\ar[r]^-\sim & \BiMod(B|B)}\,,
	\]
	where the latter equivalence follows from the Eilenberg-Watts theorem.
	
	In the second step, we choose a separable Frobenius algebra structure on $B$; every semisimple algebra admits at least one such structure. This gives us a faithful exact separable Frobenius functor
	\[
		\Vvv\:\BiMod(B|B)\to\vect
	\]
	by \Cref{expl:bimod_forget}. Consequently, we obtain a weak fiber functor $\Vvv\Psi F\:\Dd\to\vect$, and hence $A\defdtobe \End(\Vvv\Psi F)$ is a weak Hopf algebra satisfying $\Dd\cong\rep(A)$ as monoidal categories.
	
	We discuss the uniqueness of the weak Hopf algebra $A$. The conclusion is that it is \emph{far from unique}. The weak Hopf algebra $A$ has base algebra $A^l\cong B$. Since base algebras in two isomorphic weak Hopf algebras must be isomorphic, non-isomorphic choices of the separable Frobenius algebra $B$ will necessarily lead to non-isomorphic weak Hopf algebras. As can be seen from above, these non-isomorphic choices of $B$ arise from either (i) non-isomorphic choices of the semisimple algebra $B$, or (ii) non-isomorphic separable Frobenius algebra structures on $B$. Since any semisimple algebra $B$ of the form \eqref{eq:semisimple} serves the purpose, there are infinitely-many non-isomorphic choices of the semisimple algebra $B$. This imply that there exist infinitely-many non-isomorphic choices of $A$.
	
	Nevertheless, we remark that there is arguably a quasi-canonical\footnote{and unarguably the simplest} choice for $A$ \cite{Hayashi_1999}. Namely, we take $B=\kb^{\Oplus\lvert\Irr(\Mm)\rvert}$, with the unique separable Frobenius structure on it. As will be clear in \eqref{eq:a_non_trivial_equivalence}, the underlying semisimple algebra of $B$ is canonical in the sense that it has the interpretation 
	\[
		B=\Oplus_{x\in\Irr(\Mm)}\Mm(x,x)^\op\,.
	\]
	However, this interpretation does not provide any guidance on how to choose the separable Frobenius algebra structure on $B$, except to note that, in this particular case, such a structure is unique. 
	
	In \Cref{sec:2}, we will apply the general paradigm of reconstructing weak Hopf algebras outlined in this remark, where we adopt this ``quasi-canonical'' choice for $B$ \cite{Hayashi_1999}.
\end{rmk}

\section{\texorpdfstring{Reconstruction of the weak Hopf algebra $A_\Mm^\Cc$}{Reconstruction of the weak Hopf algebra \$A\_M\^{}C\$}}\label{sec:2}
Let $\Cc$ be a fusion category and $\Mm=(\Mm,\odot)$ be a finite semisimple left $\Cc$-module. Let $\fun_\Cc(\Mm,\Mm)$ denote the monoidal category of $\Cc$-module endofunctors on $\Mm$. For definitions of and general facts on fusion categories and their modules, we refer the reader to \cite{Etingof_Gelaki_Nikshych_Ostrik_2015}.

In \Cref{subsec:statement}, we directly present the statement of \Cref{mainthm:1} (\Cref{thm:main_equivalence}), which asserts that there exists certain weak Hopf algebra $A_\Mm^\Cc$ satisfying 
\eqn{\label{eq:sec:2}
	\rep(A_\Mm^\Cc)\cong\fun_\Cc(\Mm,\Mm)
}
as monoidal categories. In \Cref{subsec:ih,subsec:the_wff,subsec:the_wha}, we prove this claim using the reconstruction theorem for weak Hopf algebras. In 
\Cref{subsec:right_regular}, we show that given a fusion category $\Cc$, how to use \eqref{eq:sec:2} to obtain a weak Hopf algebra such its representation category is monoidally equivalent to $\Cc$.
\subsection{\texorpdfstring{The weak Hopf algebra $A_\Mm^\Cc$ and the statement of the main theorem}{The weak Hopf algebra \$A\_M\^{}C\$ and the statement of the main theorem}}\label{subsec:statement}
Let $\Cc,\Mm$ be defined at the begining of \Cref{sec:2}.
\begin{ntn}\label{ntn}
	Recall that for an object $a$ in a generic monoidal category, we use $a^L$ and $a^R$ to denote the left dual and right duals of $a$, respectively. The corresponding evaluation and coevaluation maps are denoted respectively by $\ev_a\:a^La\to\mone$ and $\coev_a\:\mone\to aa^L$; sometimes, we omit the subscripts for simplicity. For objects $a,b\in\Cc$ and $x\in\Mm$, we frequently write $a\os b$ as $ab$, and similarly, $a\odot x$ as $ax$. By MacLane's coherence theorem, the expression $a_1a_2\cdots a_n$ for $n\geq 3$ is unambiguous for $a_1,\cdots,a_n\in\Cc$. A similar statement holds for the expression $a_1a_2\cdots a_nx$ when $x\in\Mm$. For a $\Cc$-module functor $F\:\Mm\to\Mm$ and objects $a\in\Cc,x\in\Mm$, we denote the $\Cc$-module structure by ${F_2}_{a,x}\:aF(x)\isom F(ax)$. Lastly, we use $\Irr(\Cc)$ and $\Irr(\Mm)$ to refer to a complete set of pairwise non-isomorphic simple objects in $\Cc$ and $\Mm$, respectively. 
\end{ntn}
We define the structure of the weak Hopf algebra $A_\Mm^\Cc$ in six steps.
\paragraph{The vector space} The underlying vector space is given by 
	\[
		A_\Mm^\Cc\defdtobe\Oplus_{\substack{x,x',y,y'\in\Irr(\Mm)}}\Oplus_{a\in\Irr(\Cc)}\Mm(y',ay)\os \Mm(ax,x')\,.
	\]
\paragraph{The multiplication}
	For simple objects $y,y',\widetilde{y'},y'',x,x',\widetilde{x'},x''\in\Irr(\Mm)$, $a,b\in\Irr(\Cc)$, and elements 
	\[ u\os s\in\Mm(y'',b\widetilde{y'})\os\Mm(b\widetilde{x'},x''),\qquad v\os t\in\Mm(y',ay)\os\Mm(ax,x')
	\]
	in $A_\Mm^\Cc$, the multiplication $\mu$ reads
			\[
				\mu(u\os s\os v\os t)= 
				\begin{array}{l}
				\delta_{\widetilde{y'},y'}\delta_{\widetilde{x'},x'}\sum_{c\in\Irr(\Cc)}\sum_{\alpha=1}^{n_c}(\diagram@C=1.75pc{y''\ar[r]^-u & by'\ar[r]^-{1v} & bay \ar[r]^-{P^\alpha_c1} & cy}) \\
				\hskip -0.01em\os(\diagram@C=1.75pc{cx\ar[r]^-{I^\alpha_c1} & bax \ar[r]^-{1t} & bx' \ar[r]^-s & x''})
				\end{array}\,.
			\]
			Here, $I^\alpha_c$ and $P^\alpha_c$ represent the inclusion and projection maps, respectively, in the direct sum decomposition 
			\[b\os a\cong\oplus_{c\in\Irr(\Cc)}c^{\oplus n_c}\]
			for $\alpha=1,\cdots,n_c$.			
			
\paragraph{The unit} 
The unit $\eta\:\kb\to A_\Mm^\Cc$ is given by $\eta(1)=\sum_{x,y\in\Irr(\Mm)}\id_y\os \id_x$.
\paragraph{The comultiplication}
For $y,y',x,x'\in\Irr(\Mm)$, $a\in\Irr(\Cc)$, and element $u\os s\in\Mm(y',ay)\os\Mm(ax,x')$ in $A_\Mm^\Cc$, the comultiplication $\Delta$ reads
			\[\Delta(u\os s)=\sum_{z,z'\in\Irr(\Mm)}\sum_{\alpha=1}^{n^z_{z'}} u\os P^{z,\alpha}_{z'}\os I^{z,\alpha}_{z'} \os s\,,\]
			where for each $z\in\Irr(\Mm)$, the morphisms $I^{z,\alpha}_{z'}\:z'\to az$ and $P^{z,\alpha}_{z'}\:az\to z'$ denote the inclusions and projections, respectively, in the direct sum decomposition $az\cong \oplus_{z'\in\Irr(\Mm)}z'^{\oplus n^z_{z'}}$ for $\alpha=1,\cdots,n^z_{z'}$.
\paragraph{The counit}
For $y,y',x,x'\in\Irr(\Mm)$, $a\in\Irr(\Cc)$, and element $u\os s\in\Mm(y',ay)\os\Mm(ax,x')$ in $A_\Mm^\Cc$, the counit $\epsilon$ reads
			\[
				\epsilon(u\os s)=\delta_{y,x}\delta_{y',x'}\Lambda_{y'}(s\circ u)\,,
			\]
			where $\Lambda_{y'}\:\Mm(y',y')\to\kb$ is the unique linear map sending $\id_{y'}$ to $1$.
\paragraph{The antipode}
For $y,y',x,x'\in\Irr(\Mm)$, $a\in\Irr(\Cc)$, and element $u\os s\in\Mm(y',ay)\os\Mm(ax,x')$ in $A_\Mm^\Cc$, the antipode $S$ reads
			\eqn{\label{eq:antipode}
				S(u\os s)=s_1\os u_1\in\Mm(x,a^Rx')\os\Mm(a^Ry',y)\,.
			}
Here $s_1=(\diagram{x\ar[r]^-{\coev} & a^Rax \ar[r]^-{1s} & a^Rx'})$
and 
\[
	u_1=\sum_{\alpha=1}^{n_y} (\diagram{a^Ry'\ar[r]^-{1u} & a^Ray \ar[r]^-{1I^\alpha_y} & a^Raa^Ry' \ar[r]^-{1\ev 1} & a^Ry' \ar[r]^-{P^\alpha_y} & y})\,,
\]
where $I^\alpha_{\widetilde{y}}\:\widetilde{y}\to a^Ry'$ and $P^\alpha_{\widetilde{y}}\:a^Ry'\to y$ are the inclusions and projections, respectively, in the direct sum decomposition
\[
	a^Ry'\cong\oplus_{\widetilde{y}\in\Irr(\Mm)}\widetilde{y}^{\oplus n_{\widetilde{y}}}
\]
for $\alpha=1,\cdots,n_{\widetilde{y}}$.

\vskip 1em

Note that the maps $\mu$, $\Delta$ and $S$ do not rely on the direct sum decomposition  we choose. 
\begin{thm}[\cite{Kitaev_Kong_2012}]\label{thm:main_equivalence}
	\bnu
		\item \label{item1:thm:main_equivalence} $(A_\Mm^\Cc,\mu,\eta,\Delta,\epsilon,S)$ is a weak Hopf algebra.
		\item \label{item2:thm:main_equivalence} There exists a monoidal equivalence
		\eqn{\label{eq:thm:main_equivalence}
		\begin{split}
			K\:\fun_\Cc(\Mm,\Mm) & \to\rep(A_\Mm^\Cc) \\ 
				G &\longmapsto \Oplus_{x,y\in\Irr(\Mm)}\Mm(y,G(x))\,,
		\end{split}
		}
		where the action of $A_\Mm^\Cc$ on $K(G)$ is defined as follows: for simple objects $x,x',x_0,y,y',y_0\in\Irr(\Mm)$ and $a\in\Irr(\Cc)$, and morphisms
		\agn{
				\diagram{y'\ar[r]^-u  & ay}, \qquad \diagram{ax\ar[r]^-s & x'}, \qquad \diagram{y_0\ar[r]^-g & G(x_0)}
		}
		in $\Mm$, we have 
		\[
		(u\os s).g=\delta_{y,y_0}\delta_{x,x_0} \diagram@C=2.5pc{
			(y' \ar[r]^-{u} & ay \ar[r]^-{1g} & aG(x) \ar[r]^-{{G_2}_{a,x}}_-\sim & G(ax) \ar[r]^-{Gs} & G(x')
		})\,.
		\]
	\enu
\end{thm}
\begin{rmk}\label{rmk:ih}
	As we will see in \Cref{rmk:rmk_on_presentation}, the algebra $A_\Mm^\Cc$ actually takes a more concise form
	\[
		\Oplus_{x,x',y,y'\in\Irr(\Mm)}\Cc(\mone,[x,x'][y',y])
	\]
	if one employs the language of internal homs introduced in \Cref{subsec:ih}.
\end{rmk}

\begin{rmk}\label{rmk:KK12_pf}
We note that, up to the minor differences that will be discussed in \Cref{rmk:diff_of_alg},  \Cref{thm:main_equivalence} (and in particular \Cref{thm:main_equivalence}.\ref{item2:thm:main_equivalence}) was proposed in \cite[\S 4]{Kitaev_Kong_2012}, and also sketchily proved there. The key point of their proof of \Cref{thm:main_equivalence}.\ref{item2:thm:main_equivalence} is to disclose that the defining data of a $\Cc$-module functor is equivalent to the defining data of a left module over $A_\Mm^\Cc$ \cite[Eqs. (27-30)]{Kitaev_Kong_2012}. We refer the reader to \cite[\S 4]{Kitaev_Kong_2012} for the nice and self-evident graphical intuitions behind the structure maps of $A_\Mm^\Cc$ and the equivalence \eqref{eq:thm:main_equivalence}, which complements the present article.

We wish also to informally comment on other potential proofs of \Cref{thm:main_equivalence}. First of all, we believe a proof of \Cref{thm:main_equivalence} based on \cite[Proposition 10]{Barter_Bridgeman_Jones_2019a} and a well-known equivalence between $\fun_\Cc(\Mm,\Mm)$ and certain relative tensor product of module categories is possible. Secondly, a purely graphical proof, which retains the maximal graphical intuition behind \Cref{thm:main_equivalence} (none of which is preserved in the present article), could potentially be developed \cite{Kitaev_Kong_2012, Morrison_Walker_2012, Hoek_2019, Liu_Ming_Wang_Wu_2023, Jia_Tan_Kaszlikowski_2024}. In fact, a proof of \Cref{thm:main_equivalence}.\ref{item1:thm:main_equivalence} based on bordism categories is already provided in \cite[\S 3.2]{Cordova_Holfester_Ohmori_2024}, drawing from unpublished works by Johnson-Freyd and Reutter. It remains unknown to the authors whether this proof can be extended to a full proof of \Cref{thm:main_equivalence}.
\end{rmk}

\begin{rmk}\label{rmk:diff_of_alg}
The minor differences between the algebra $A_\Mm^\Cc$ and the algebra introduced in \cite[\S 4]{Kitaev_Kong_2012} (denoted by $\Au_\Mu^\Cu$) include the following: (a) The algebra $\Au_\Mu^\Cu$ is a $C^\ast$-weak Hopf algebra, which requires $\Cu$ to be a unitary fusion category and $\Mu$ being a unitary module. In contrast, $A_\Mm^\Cc$ here is only a weak Hopf algebra, without requiring the unitary structures on $\Cc$ and $\Mm$. (b) The algebra $A_\Mm^\Cc$ has reversed comultiplication as $\Au_\Mu^\Cu$. The difference (b) leads to a \emph{warning}: by \Cref{rmk:antipode}.\ref{item4:rmk:antipode}, the antipode $S$ of $A_\Mm^\Cc$ corresponds to the \emph{inverse} of the antipode of $\Au_\Mu^\Cu$ given in \cite[Eq. (24)]{Kitaev_Kong_2012}. Lastly, we comment on another difference between $A_\Mm^\Cc$ and $\Au_\Mu^\Cu$, namely the disparity between the ``perplexed'' form of the antipode of $S$ of $A_\Mm^\Cc$ given in \eqref{eq:antipode}, and the simpler form of the antipode given in \cite[Eq. (24)]{Kitaev_Kong_2012}. We expect that this discrepancy is accounted for by (a): when $\Cc$ and $\Mm$ carry certain additional structures such as unitarity, the antipode $S$ may reduce to the simpler form.
\end{rmk}
\begin{rmk}\label{rmk:phy}
	Let us briefly present the physical application of $A_\Mm^\Cc$ in Levin-Wen models appearing in the original article \cite{Kitaev_Kong_2012}; we refer the reader to the latter and also \cite{Kong_2012, Lan_Wen_2014} for further discussions. In \cite{Kitaev_Kong_2012}, a topological excitation on the $\Mm$-boundary of a $\Cc$-Levin-Wen model is identified with a left $A_\Mm^\Cc$-module. The fusion of two topological excitations is given by the tensor product in $\rep(A_\Mm^\Cc)$, that is, governed by the comultiplication $\Delta$ (cf. \cite[Figure 6]{Kitaev_Kong_2012}). It is also implicit in \cite{Kitaev_Kong_2012} that the vacuum excitation corresponds to the tensor unit of $\rep(A_\Mm^\Cc)$.
	
	For other physical applications of $A_\Mm^\Cc$, see for instance \cite{Cordova_Holfester_Ohmori_2024, Inamura_Ohyama_2024, Choi_Rayhaun_Zheng_2024b} and the references therein.
\end{rmk}

The next three subsections are devoted to the proof of \Cref{thm:main_equivalence}. We believe that the reader focused on applications of this theorem may safely skip them.

\subsection{Recap of internal homs}\label{subsec:ih}
In this subsection, we recall some basic facts about internal homs, a powerful tool in tensor category theory that emerged in the early development of the theory \cite{Ostrik_2003, Etingof_Ostrik_2004}.

Let $\Cc,\Mm$ be as in \Cref{subsec:statement}. 
\begin{dfn}
	For $x\in\Mm$, we denote the right adjoint of the functor $-\odot x\:\Cc\to\Mm, \; a\longmapsto ax$ by $[x,-]\:\Mm\to\Cc$. We denote the image of $y\in\Mm$ under the functor $[x,-]$ by $[x,y]_\Cc$, or simply $[x,y]$. We call $[x,y]$ the \newdef{internal hom} from $x$ to $y$.
\end{dfn}
\begin{rmk}
	The right adjoint functor $[x,-]$ always exists.
\end{rmk}
By definition, we have a natural isomorphism
\eqn{\label{eq:adj}
	\Mm(ax,y)\isom\Cc(a,[x,y])
}
for any $a\in\Cc, y\in\Mm$; let the counit and the unit of the adjunction \eqref{eq:adj} be denoted by 
\[
	\epsilon_{x,y}\:[x,y]x\to y \quad \text{and}\quad \eta_{a,x}\:a\to [x,ax]
\]
respectively. 
\begin{expl}\label{expl:ih}
	\bnu[(1)]
		\item Treat $\Cc$ as a left module over itself. Then for $x,y\in\Cc$, we have $[x,y]_\Cc=yx^L$.
		\item Treat $\Cc$ as a left module over $\Cc^\rev$ with $\Cc^\rev\times\Cc\to\Cc,\; (a,x)\mapsto xa$. Then for $x,y\in\Cc$, we have  $[x,y]_{\Cc^\rev}=x^Ry$.
	\enu
\end{expl}

Note that $[x,-]\:\Mm\to\Cc$ is automatically a left $\Cc$-module functor \cite[Corollary 7.9.5]{Etingof_Gelaki_Nikshych_Ostrik_2015}. Its $\Cc$-module structure 
\[
	{[x,-]_2}_{a,y}\:a[x,y]\isom [x,ay]
\]
for $a\in\Cc,y\in\Mm$ is given by the image of (\(
	\diagram@C=2.5pc{a[x,y]x\ar[r]^-{1\epsilon_{x,y}} & ay}
\))
under the isomorphism
\[
	\Mm(a[x,y]x,ay)\isom\Cc(a[x,y],[x,ay])\,.
\]

We need to introduce more natural maps regarding internal homs: given objects $x,y,z\in\Mm$, we denote by $\mu_{x,y,z}\:[y,z][x,y]\to [x,z]$ the image of
\[
	\diagram@C=2.75pc{[y,z][x,y]x\ar[r]^-{1\epsilon_{x,y}} & [y,z]y\ar[r]^-{\epsilon_{y,z}} & z}	
\]
under the isomorphism
\[
	\Mm([y,z][x,y]x,z)\isom\Cc([y,z][x,y],[x,z])\,.
\]
We also denote $\eta_x\defdtobe\eta_{\mone,x}\:\mone\to[x,x]$. 
\begin{rmk}
One can check that the maps $\{\mu_{x,y,z}\}_{x,y,z\in\Mm}$ and $\{\eta_x\}_{x\in\Mm}$ satisfy a ``generalized'' associativity and unitality conditions:
\bnu
	\item For any $x,y,z,w\in\Mm$, we have 
	\[
		\mu_{x,y,w}\circ(\mu_{y,z,w}\os\id_{[x,y]})=\mu_{x,z,w}\circ(\id_{[z,w]}\os\mu_{x,y,z})\:[z,w][y,z][x,y]\to[x,w]\,.
	\]
	\item For any $x,y\in\Mm$, we have
	\[
		\mu_{x,y,y}\circ(\eta_y\os\id_{[x,y]})=\id_{[x,y]}=\mu_{x,x,y}\circ(\id_{[x,y]}\os\eta_x)\:[x,y]\to[x,y]\,.
	\]
\enu
\end{rmk}

\subsection{\texorpdfstring{The weak fiber functor $\fun_\Cc(\Mm,\Mm)\rightarrow\vect$}{The weak fiber functor \$Fun\_C(M,M)\text{->}Vect\$}}\label{subsec:the_wff}
In this and the following subsection, we continue to prove \Cref{thm:main_equivalence}. In this subsection, we construct a weak fiber functor $\Fff\:\fun_\Cc(\Mm,\Mm)\to\vect$ on $\fun_\Cc(\Mm,\Mm)$, which implies that $\fun_\Cc(\Mm,\Mm)\cong\rep(\End(\Fff))$ by \Cref{thm:reconstruct}. In the next subsection, we establish an isomorphism of weak Hopf algebras $A_\Mm^\Cc\cong\End(\Fff)$ by explicitly specifying the structure maps of $\End(\Fff)$, yielding a proof of \Cref{thm:main_equivalence}.

Notice that $\fun_\Cc(\Mm,\Mm)$ is a finite multi-tensor category by \cite[Proposition 7.11.6 and Exercise 7.12.1]{Etingof_Gelaki_Nikshych_Ostrik_2015}. We define a weak fiber functor $\Fff$ on $\fun_\Cc(\Mm,\Mm)$ as the composition of the following three faithful exact separable Frobenius monoidal functors:
\[
	\diagram@C=1.5pc{\fun_\Cc(\Mm,\Mm)\ar[r]^-{\Gamma} & \fun(\Mm,\Mm) \ar[r]^-{\Psi}_-\cong & \BiMod(\kb^{\Oplus\lvert\Irr(\Mm)\rvert}|\kb^{\Oplus\lvert\Irr(\Mm)\rvert}) \ar[r]^-{\Vvv} & \vect}\,.
\]
Note that this construction of $\Fff$ fits within the general reconstruction paradigm outlined in \Cref{rmk:to_funMM}. We now introduce the functors $\Vvv$, $\Psi$ and $\Gamma$ as follows.
\paragraph{The functor $\Vvv\:\BiMod(\kb^{\Oplus\lvert\Irr(\Mm)\rvert}|\kb^{\Oplus\lvert\Irr(\Mm)\rvert})\rightarrow\vect$}
We endow $\kb^{\Oplus\lvert\Irr(\Mm)\rvert}$ with the unique separable Frobenius algebra structure, i.e., the one given by the direct sum of $\vert\Irr(\Mm)\rvert$ copies of the trivial algebra $\kb$. Let 
\[\Vvv\:\BiMod(\kb^{\Oplus\lvert\Irr(\Mm)\rvert}|\kb^{\Oplus\lvert\Irr(\Mm)\rvert})\to\vect\]
be the separable Frobenius functor associated with $\kb^{\Oplus\lvert\Irr(\Mm)\rvert}$ in \Cref{expl:bimod_forget}. It is faithful and exact.

\paragraph{The functor $\Psi\:\fun(\Mm,\Mm)\rightarrow\BiMod(\kb^{\Oplus\lvert\Irr(\Mm)\rvert}|\kb^{\Oplus\lvert\Irr(\Mm)\rvert})$}
It is clear that any finite semisimple category $\Mm$ is equivalent to the category of left modules over the algebra $\kb^{\Oplus\lvert\Irr(\Mm)\rvert}$. Expressing this fact in a slightly more basis-independent way\footnote{or equivalently, invoking the reconstruction theorem for ordinary algebras (\Cref{thm:monadicity}) on the faithful exact representable functor $\Mm(\oplus_{y\in\Irr(\Mm)}y,-)\:\Mm\to\vect$.}, one can say that the functor
\eqn{\label{eq:a_non_trivial_equivalence}
\begin{split}
	\Mm & \textstyle\isom\rep(\Oplus_{y\in\Irr(\Mm)}\Mm(y,y)^\op)\isom\rep(\kb^{\Oplus\lvert\Irr(\Mm)\rvert}) \\
	w& \longmapsto \textstyle\Oplus_{y\in\Irr(\Mm)}\Mm(y,w)
\end{split}
}
is an equivalence of categories, where the left $\Oplus_{y\in\Irr(\Mm)}\Mm(y,y)^\op$-action on $\Oplus_{y\in\Irr(\Mm)}\Mm(y,w)$ is induced from the evident \emph{right} action of $\Oplus_{y\in\Irr(\Mm)}\Mm(y,y)$.

Now we define the monoidal functor $\Psi$ to be the composition of the two monoidal equivalences
\[
	\diagram{\fun(\Mm,\Mm) \ar[r]^-\sim & \fun(\rep(\kb^{\Oplus\lvert\Irr(\Mm)\rvert}),\rep(\kb^{\Oplus\lvert\Irr(\Mm)\rvert}) \ar[r]^-\sim & \BiMod(\kb^{\Oplus\lvert\Irr(\Mm)\rvert}|\kb^{\Oplus\lvert\Irr(\Mm)\rvert})}\,,
\]
where the second equivalence follows from the Eilenberg-Watts theorem. Then the explicitly form of $\Psi$ is given by 
\agn{
	\Psi\: \fun(\Mm,\Mm) & \to \BiMod(\kb^{\Oplus\lvert\Irr(\Mm)\rvert}|\kb^{\Oplus\lvert\Irr(\Mm)\rvert}) \\
	G & \longmapsto	\Oplus_{x,y\in\Irr(\Mm)}\Mm(y,G(x))\,.
}
Here, importantly, the left $\kb^{\Oplus\lvert\Irr(\Mm)\rvert}\cong\Oplus_{y\in\Irr(\Mm)}\Mm(y,y)^\op$-action on $\Psi(G)$ is given as follows: for $x_0,y_0,y\in\Irr(\Mm)$ and morphisms
\[\diagram{y\ar[r]^-u & y}\quad\text{and}\quad\diagram{y_0\ar[r]^-v & G(x_0)}\,,\]
there is $u.v=\delta_{y_0,y}(v\circ u)$. The right $\kb^{\Oplus\lvert\Irr(\Mm)\rvert}\cong\Oplus_{x\in\Irr(\Mm)}\Mm(x,x)^\op$-action on $\Psi(G)$ is given as follows: for $x_0,x,y_0\in\Irr(\Cc)$ and morphisms
\[\diagram{x\ar[r]^-s & x}\quad\text{and}\quad\diagram{y_0\ar[r]^-v & G(x_0)},\]
there is $v.s=\delta_{x_0,x}(G(s)\circ v)$. 

$\Psi$ is a faithful and exact separable Frobenius monoidal functor since it is a monoidal equivalence.
\paragraph{The functor $\Gamma\:\fun_\Cc(\Mm,\Mm)\rightarrow\fun(\Mm,\Mm)$} We define $\Gamma$ to be the forgetful functor
\[
	\fun_\Cc(\Mm,\Mm)\to\fun(\Mm,\Mm)
\]
sending each $\Cc$-module functor to its underlying functor. It is naturally a strong monoidal functor, hence a separable Frobenius monoidal functor. It is faithful by definition. We still need to prove that $\Gamma$ is exact.
\begin{lem}\label{lem:Gamma_admit_adj}
	The functor 
	\[
		L\:\fun(\Mm,\Mm)\to\fun_\Cc(\Mm,\Mm),\quad F\longmapsto \oplus_{x\in\Irr(\Mm)}[x,-]\odot F(x)
	\]
	is left adjoint to $\Gamma$. Similarly, the functor 
	\[
		R\:\fun(\Mm,\Mm)\to\fun_\Cc(\Mm,\Mm),\quad F\longmapsto \oplus_{x\in\Irr(\Mm)}[-,x]^R\odot F(x)
	\]
	is right adjoint to $\Gamma$. In particular, the functor $\Gamma$ is exact.
	\pf
	We only show the part for $L$, and leave the proof of the other part to the reader. First, $[x,-]\odot F(x)$ is a well-defined $\Cc$-module functor since it is the composition of the $\Cc$-module functor $[x,-]\:\Mm\to\Cc$ in \Cref{subsec:ih} and the $\Cc$-module functor $-\odot F(x)\:\Cc\to\Mm$. 
	
To show $L\ladj \Gamma$, it is enough to construct a natural isomorphism
\eqn{\label{eq:lem:Gamma_admit_adj}
	b_{F,G}\:\fun_\Cc(\Mm,\Mm)(L(F),G)\isom\fun(\Mm,\Mm)(F,G)
}
for a functor $F\:\Mm\to\Mm$ and a $\Cc$-module functor $(G,G_2)\:\Mm\to\Mm$. Note that we have 
\[
	\fun_\Cc(\Mm,\Mm)(L(F),G)\cong\Oplus_{x\in\Irr(\Mm)}\fun_\Cc(\Mm,\Mm)([x,-]F(x),G)\,.
\]
Therefore, to define the isomorphism \eqref{eq:lem:Gamma_admit_adj}, it is enough to realize $\fun(\Mm,\Mm)(F,G)$ as the direct sum of $\{\fun_\Cc(\Mm,\Mm)([x,-]F(x),G)\}_{x\in\Irr(\Cc)}$.

To this end, we define
\agn{
	P_x\:\fun(\Mm,\Mm)(F,G) & \to\fun_\Cc(\Mm,\Mm)([x,-]F(x),G) \\
	\beta & \longmapsto P_x(\beta)\,,
}
where ${P_x(\beta)}_y$ is defined as the composition 
\[\diagram@C=3pc{[x,y]F(x)\ar[r]^-{1\beta_x} & [x,y]G(x) \ar[r]^-{{G_2}_{[x,y],x}} & G([x,y]x) \ar[r]^-{G\epsilon_{x,y}} & G(y)}\]
 for $y\in\Mm$. That $P_x(\beta)$ is indeed a $\Cc$-module natural transformation can be seen by the commutativity of the outermost diagram in 
\[
	\diagram@C=4pc{
		a[x,y]F(x)\ar[d]_{1\beta_x}  \ar[rr]^-{{[x,-]_2}_{a,y}1} & & [x,ay]F(x) \ar[d]^{1\beta_x}  \\
		a[x,y]G(x) \ar[rd]^-{{G_2}_{a[x,y],x}}_-{\spadesuit\hskip 1.5pc} \ar[d]_{1{G_2}_{[x,y],x}}  \ar[rr]^-{{[x,-]_2}_{a,y}1} & & [x,ay]G(x) \ar[d]^{{G_2}_{[x,ay],x}}  \\
		aG([x,y]x) \ar[d]_{1G(\epsilon_{x,y})}  \ar[r]_-{{G_2}_{a,[x,y]x}} & G(a[x,y]x)  \ar[rd]_-{G(1\epsilon_{x,y})}^-{\hskip 1.5pc \clubsuit} \ar[r]^-{G({[x,-]_2}_{a,y}1)}& G([x,ay]x) \ar[d]^{G(\epsilon_{x,ay})} \\
		aG(y) \ar[rr]_-{{G_2}_{a,y}} & & G(ay)
	}
\]
The commutativity of ($\scriptstyle\spadesuit$) follows from that $G$ is a $\Cc$-module functor, and the commutativity of ($\scriptstyle\clubsuit$) follows from the definition of ${[x,-]_2}_{a,y}$ in \Cref{subsec:ih}. 

Next, we define 
\agn{
	I_x\:\fun_\Cc(\Mm,\Mm)([x,-]F(x),G) & \to\fun(\Mm,\Mm)(F,G) \\
			\gamma & \longmapsto I_x(\gamma)\,,
}
where for a simple object $x'\in\Irr(\Cc)$, we set ${I_x(\gamma)}_{x'}$ to be 
\[
	\begin{cases} 
		0,\phantom{\diagram{F(x) \ar[r]^-{\eta_x1} & [x,x]F(x) \ar[r]^-{\gamma_x} & G(x)}}\quad \text{if $x'\neq x$}\,; \\
	\diagram{F(x) \ar[r]^-{\eta_x1} & [x,x]F(x) \ar[r]^-{\gamma_x} & G(x)}\phantom{0},\quad\text{otherwise.}
	\end{cases}
\]
It is not hard to verify that $\sum_{x\in\Irr(\Mm)}I_x\circ P_x=\id$ and $P_{x'}\circ I_x=\delta_{x,x'}\id$. Therefore, the morphisms 
\[\{I_x\}_{x\in\Irr(\Mm)} \quad\text{and}\quad\{P_x\}_{x\in\Irr(\Mm)}
\]
establish $\fun(\Mm,\Mm)(F,G)$ as the direct sum of $\{\fun_\Cc(\Mm,\Mm)([x,-]F(x),G)\}_{x\in\Irr(\Cc)}$, uniquely determining an isomorphism $\fun_\Cc(\Mm,\Mm)(L(F),G)\isom\fun(\Mm,\Mm)(F,G)$. We define $b_{F,G}$ to be this isomorphism.

It suffices to show that $b_{F,G}$ is natural in $F$ and $G$, which we leave to the reader.
\epf
\end{lem}
\Cref{lem:Gamma_admit_adj} concludes our construction of the faithful and exact separable Frobenius monoidal functor $\Gamma$.
\bigskip
\begin{prp}\label{prp:the_main_wff}
	The functor 
	\agn{
		\Fff\:\fun_\Cc(\Mm,\Mm) & \to\vect	 \\
					G & \longmapsto \Oplus_{x,y\in\Irr(\Mm)}\Mm(y,G(x))
	}
	has a structure of weak fiber functor. Its separable Frobenius monoidal structure is given as follows:
	\bnu[1.]
		\item For $G,F\in\fun_\Cc(\Mm,\Mm)$ and $x,y,y',z\in\Irr(\Mm)$, the value of the map
		\[{\Fff_2}_{G,F}\:\Fff(G)\os\Fff(F)\to\Fff(GF)\]
		at $(\diagram{z\ar[r]^-g &  G(y)}) \os (\diagram{y'\ar[r]^-f & F(x)})$ is $\delta_{y,y'} (G(f)\circ g)$.
		\item $\Fff_0\:\kb\to\Fff(\Id_\Mm)\cong\Oplus_{y\in\Irr(\Mm)}\Mm(y,y)$ is given by $1\longmapsto \sum_{y\in\Irr(\Mm)}\id_y$.
		\item For $G,F\in\fun_\Cc(\Mm,\Mm)$ and $x,z\in\Irr(\Mm)$, the value of the map 
		\[{\Fff_{-2}}_{G,F}\:\Fff(GF)\to\Fff(G)\os\Fff(F)\]
		at $(\diagram{ z\ar[r]^-h & GF(x)})$ is given in two steps. First, we fix a decomposition $F(x)\cong\oplus_{y\in\Irr(\Mm)}y^{\oplus n_y}$ with inclusion maps $I^\alpha_y\:y\to F(x)$ and projection maps $P^\alpha_y\:F(x)\to y$ for $\alpha=1,\cdots,n_y$ and $y\in\Irr(\Mm)$. Next, we define
		\[
			{\Fff_{-2}}_{G,F}(h)\defdtobe\sum_{y\in\Irr(\Mm)}\sum_{\alpha=1}^{n_y}(\diagram@C=1.85pc{z\ar[r]^-h & GF(x)\ar[r]^-{GP^\alpha_y} & G(y)})\os(\diagram@C=1.85pc{y \ar[r]^-{I^\alpha_y} & F(x)})\,.
		\]
		\item $\Fff_{-0}\:\Fff(\Id_\Mm)\cong\Oplus_{y\in\Irr(\Mm)}\Mm(y,y)\to\kb$ is induced by the linear maps
		\[
			\Lambda_y\:\Mm(y,y)\to\kb,\quad s\cdot \id_y\longmapsto s\,.
		\]
	\enu
	\pf
		It is enough to take $\Fff\defdtobe\Vvv\Psi\Gamma\:\fun_\Cc(\Mm,\Mm)\to\vect$.
	\epf
\end{prp}
\begin{cor}\label{cor:the_main_wff}
	$\End(\Fff)$ has a natural structure of weak Hopf algebra. Moreover, there exists a monoidal equivalence
	\agn{
		\widetilde{\Fff}\:\fun_\Cc(\Mm,\Mm) & \to  \rep(\End(\Fff)) \\
		G & \longmapsto \Fff(G)=\Oplus_{x,y\in\Irr(\Mm)}\Mm(y,G(x))\,,
	}
	where the action of $\End(\Fff)$ on $\Fff(G)$ is given by 
	\[
		\End(\Fff)\os\Fff(G)\to\Fff(G),\quad \alpha\os g\longmapsto \alpha_G(g)\,.
	\]
	\pf
		The proof directly follows from \Cref{prp:the_main_wff}, and \ref{item1:thm:reconstruct}, \ref{item3:thm:reconstruct} of \Cref{thm:reconstruct}.
	\epf
\end{cor}
By \Cref{cor:the_main_wff}, \Cref{thm:main_equivalence} can be immediately proved once we can show that $A_\Mm^\Cc$ is a weak Hopf algebra isomorphic to $\End(\Fff)$. This is the subject of the next subsection.

\subsection{The reconstruction process}\label{subsec:the_wha}

In this subsection, we finish the proof of \Cref{thm:main_equivalence} by establishing an isomorphism $\End(\Fff)\cong A_\Mm^\Cc$ of weak Hopf algebras.

Let $F_0\:\Mm\to\Mm$ be the functor defined by $F_0(x)\defdtobe\oplus_{y\in\Irr(\Mm)}y$ for all $x\in\Irr(\Mm)$.
\begin{lem}\label{lem:represent}
	 The functor $\Fff\:\fun_\Cc(\Mm,\Mm)\to\vect$ is represented by $L(F_0)$, where $L$ is defined in \Cref{lem:Gamma_admit_adj}.
	\pf
		For any $\Cc$-module functor $G\:\Mm\to\Mm$, we have 
		\[
			\fun_\Cc(\Mm,\Mm)(L(F_0),G)\cong\fun(\Mm,\Mm)(F_0,\Gamma(G))\cong\Oplus_{x,y\in\Irr(\Mm)}\Mm(y,G(x))\,.
		\]
	\epf
\end{lem}
\begin{rmk}
	Explicitly, $L(F_0)(x')$ for $x'\in\Irr(\Mm)$ is given by
	\(
		\oplus_{x,y\in\Irr(\Mm)}[x,x']y
	\).
\end{rmk}
\begin{cor}\label{cor:vs_transport}
As a vector space, $\End(\Fff)$ is canonically isomorphic to
\[{}^\wr A_\Mm^{\Cc}\defdtobe \Oplus_{x,x',y,y'\in\Irr(\Mm)}\Mm(y',[x,x']y)\,.\]

To be precise, there is a canonical linear isomorphism from $\End(\Fff)$ to ${}^\wr A_\Mm^\Cc$ given by \eqn{\label{eq1:cor:vs_transport}
		\begin{split}
		\phi\:\End(\Fff)& \isom {}^\wr A_\Mm^{\Cc}
 \\
		\gamma & \longmapsto\gamma_{L(F_0)}({}^\wr 1)\,.
		\end{split}
		}
Here ${}^\wr 1\in \Fff L(F_0)={}^\wr A_\Mm^\Cc$ is the distinguished element whose component ${}^\wr 1_{x,x';y,y'}$ in $\Mm(y',[x,x']y)$ for $x,x',y,y'\in\Irr(\Mm)$ reads 
		\eqn{\label{eq:wr_1}
		{}^\wr 1_{x,x';y,y'} = \begin{cases} 0,\phantom{\diagram@C=3pc{y  \ar[r]^-{\eta_x1} & [x,x]y}} \text{if $x\neq x'$ or $y\neq y'$}\,; \\
		\diagram@C=3pc{y  \ar[r]^-{\eta_x1} & [x,x]y},\phantom{0} \text{otherwise,}
		\end{cases}
		}	
		where $\eta_x$ is defined in \Cref{subsec:ih}. 
		
		The inverse of $\phi$ is given by 
			\agn{
				\phi\inverse\: {}^\wr A_\Mm^{\Cc} & \isom \End(\Fff)\\
					(\diagram@C=1pc{y'\ar[r]^-f & [x,x']\odot' y}) & \longmapsto \phi\inverse(f)\:\Fff\funto\Fff\,,
			}
			where $\phi\inverse(f)$ reads, componentwise,
			\agn{\phi\inverse(f)_G\:\Fff(G) & \to \Fff(G)\\
			(\diagram@C=1pc{y_0\ar[r]^-{g} & G(x_0)}) & \longmapsto 
			\begin{array}{l}		
			\delta_{x_0,x}\delta_{y_0,y}(\diagram{y' \ar[r]^-f & [x,x']y \ar[r]^-{1g} & [x,x']G(x) }\\
			\hskip -0.5em \diagram@C=4pc{\ar[r]^-{{G_{2}}_{[x,x'],x}} & G([x,x']x) \ar[r]^-{G(\epsilon_{x,x'})} & G(x')})
			\end{array}
			}
			for $G\in\fun_\Cc(\Mm,\Mm)$ and  $x_0,y_0\in\Irr(\Mm)$.
	\pf
		It is enough to take $\phi$ as the composition of the following isomorphisms:
		\[
		\begin{multlined}
			\phi\:\End(\Fff)\isom\fun_\Cc(\Mm,\Mm)(L(F_0),L(F_0))\isom\Fff L(F_0)
			\\
			=\Oplus_{x',y'\in\Irr(\Mm)}\Mm(y',L(F_0)(x'))=\Oplus_{x,x',y,y'\in\Irr(\Mm)}\Mm(y',[x,x']y)={}^\wr A_\Mm^\Cc\,,
		\end{multlined}
		\]
		where the first isomorphism comes from the Yoneda lemma.
	\epf
\end{cor}
Let $(\mu',\eta',\Delta',\epsilon',S')$ represent the weak Hopf algebra structure on $\End(\Fff)$ given by \Cref{cor:the_main_wff}. Using the isomorphism $\phi$ in \eqref{eq1:cor:vs_transport}, we can ``transport'' the weak Hopf algebra structure on $\End(\Fff)$ to ${}^\wr A_\Mm^\Cc$. That is,
define linear maps
\eqn{\label{eq:with_wr}
	\begin{split}
	&\manualformatting \hskip 4em {}^\wr\mu=\phi\circ \mu'\circ (\phi\inverse\os  \phi\inverse)\qquad\quad  {}^\wr \eta=\phi\circ  \eta'
\\
	{}^\wr\Delta&=(\phi\os \phi)\circ \Delta' \circ \phi\inverse\qquad\quad {}^\wr\epsilon= \epsilon' \circ \phi\inverse \qquad\quad {}^\wr S=\phi \circ S'\circ \phi\inverse\,.
	\end{split}
}
Then it is trivial to see that $({}^\wr A_\Mm^\Cc,{}^\wr\mu,{}^\wr\eta,{}^\wr\Delta,{}^\wr\epsilon,{}{}^\wr S)$ is a weak Hopf algebra isomorphic to $\End(\Fff)$, and that by \Cref{cor:the_main_wff}, there exists a monoidal equivalence
\[
	\fun_\Cc(\Mm,\Mm)\to\rep({}^\wr A_\Mm^\Cc)\, .
\]
These data are explicitly computed as follows.
\begin{thm}\label{thm:with_wr}
\bnu
	\item The maps ${}^\wr\mu,{}^\wr\eta,{}^\wr\Delta,{}^\wr\epsilon,{}{}^\wr S$ are given as follows:
	\bnu
		\item \label{item11:thm:with_wr} For simple objects $y,y',\widetilde{y'},y'',x,x',\widetilde{x'},x''\in\Irr(\Mm)$, and elements
		\[\diagram@C=2pc{y''\ar[r]^-g & [\widetilde{x'},x'']\widetilde{y'}},\qquad\quad \diagram@C=2pc{y'\ar[r]^-f & [x,x']y}\] 
		in ${}^\wr A_\Mm^\Cc$, we have
\agn{
	{}^\wr\mu(g\os f)=
	\delta_{x',\widetilde{x'}}\delta_{y',\widetilde{y'}} (\diagram@C=3.2pc{
	y'' \ar[r]^-{g} & [x',x'']y' \ar[r]^-{1f} & [x',x''][x,x']y
	 \ar[r]^-{\mu_{x,x',x''}1} & [x,x'']y
})\,,
}
	where $\mu_{x,x',x''}$ is defined in \Cref{subsec:ih}.
		\item \label{item12:thm:with_wr} The map ${}^\wr\eta$ sends $1\in\kb$ to ${}^\wr 1\in {}^\wr A_\Mm^\Cc$ defined in \eqref{eq:wr_1}.
		\item \label{item13:thm:with_wr} For simple objects $y,y',x,x'\in\Irr(\Mm)$ and  $(\diagram@C=1pc{y'\ar[r]^-f &[x,x']y})\in{}^\wr A_\Mm^\Cc$, the value ${}^\wr\Delta(f)$ is given in two steps. First, we choose a direct sum decomposition
		\[
			[x,x']z \cong\oplus_{z'\in \Irr(\Mm)} z'^{\oplus n^z_{z'}}
		\]
		for each $z\in\Irr(\Mm)$, with inclusion maps and projection maps given respectively by
		\[
			I^{z,\alpha}_{z'}\:z'\to [x,x']z \quad\text{and}\quad P^{z,\alpha}_{z'}\:[x,x']z\to z'
		\]
		for $\alpha=1,\cdots,n^z_{z'}$ and $z'\in\Irr(\Mm)$. Also, let $(P^{z,\alpha}_{z'})^\sharp\:[x,x']\to[z,z']$ be the map induced from $P^{z,\alpha}_{z'}$ via the adjunction $-\odot z\ladj [z,-]$. In the second step, we set
		\[
			{}^\wr\Delta(f)=\sum_{z,z'\in\Irr(\Mm)}\sum_{\alpha=1}^{n^z_{z'}} f^{z,\alpha}_{z',1}\os  f^{z,\alpha}_{z',2}\,,
		\]
		where $f^{z,\alpha}_{z',1}=(\diagram@C=2.25pc{y'\ar[r]^-f & [x,x']y \ar[r]^-{(P^{z,\alpha}_{z'})^\sharp 1} & [z,z']y})$ and $f^{z,\alpha}_{z',2}=I^{z,\alpha}_{z'}$.
		\item \label{item14:thm:with_wr} For simple objects $y,y',x,x'\in\Irr(\Mm)$ and $(\diagram@C=1pc{y'\ar[r]^-f &[x,x']y})\in{}^\wr A_\Mm^\Cc$, we have 
		\[
			{}^\wr\epsilon(f)=\delta_{x,y}\delta_{x',y'}\Lambda_{y'}(\diagram{y'\ar[r]^-f & [y,y']y\ar[r]^-{\epsilon_{y,y'}} & y'})\,,
		\]
		where $\Lambda_{y'}\:\Mm(y',y')\to\kb$ is the unique linear map sending $\id_{y'}$ to $1\in\kb$.
		\item \label{item15:thm:with_wr} For simple objects $y,y',x,x'\in\Irr(\Mm)$ and  $(\diagram@C=1pc{y'\ar[r]^-f &[x,x']y})\in{}^\wr A_\Mm^\Cc$, the value ${}^\wr S(f)$ is given in two steps. First, recall that $[x,x']^R$ denotes the right dual of $[x,x']$. Take a direct sum decomposition
		\[
			[x,x']^Ry'\cong\oplus_{\widetilde{y}\in\Irr(\Mm)}\widetilde{y}^{\oplus n_{\widetilde{y}}}\,,
		\]
		with inclusion maps and projection maps given respectively by
		\[
			I^\alpha_{\widetilde{y}}\:\widetilde{y}\to[x,x']^R y' \quad\text{and}\quad P^\alpha_{\widetilde{y}}\:[x,x']^Ry'\to \widetilde{y}
		\]
		for $\alpha=1,\cdots,n_{\widetilde{y}}$ and $\widetilde{y}\in\Irr(\Mm)$. Let $(P^\alpha_{\widetilde{y}})^\sharp\:[x,x']^R\to[y',\widetilde{y}]$ be the map induced from $P^\alpha_{\widetilde{y}}$ via the adjunction $-\odot y'\ladj [y',-]$. Secondly, there is
		\[
			{}^\wr S(f)=\sum_{\alpha=1}^{n_y}\Lambda_{y'}(f^\alpha_1)f^\alpha_2\,,
		\]
		where 
		\[f_1^\alpha=(\diagram@C=2.5pc{y'\ar[r]^-f & [x,x']y \ar[r]^-{1I^\alpha_y} & [x,x'][x,x']^Ry' \ar[r]^-{\ev 1} & y'})\]
		and \manualformatting
		\[f_2^\alpha=(\diagram@C=2.5pc{x\ar[r]^-{\coev 1} & [x,x']^R[x,x']x\ar[r]^-{1\epsilon_{x,x'}} & [x,x']^Rx'\ar[r]^-{(P^\alpha_y)^\sharp 1} & [y',y]x'})\,.
		\]
	\enu
	\item \label{item2:thm:with_wr} There exists a monoidal equivalence
	\eqn{\label{eq:thm:with_wr}
	\begin{split}
		{}^\wr K\:\fun_\Cc(\Mm,\Mm)& \to \rep({}^\wr A_\Mm^\Cc) \\
		G & \longmapsto \Oplus_{x,y\in\Irr(\Mm)}\Mm(y,G(x))\,,
	\end{split}
	}
	where the action of ${}^\wr A_\Mm^\Cc$ on ${}^\wr K(G)$ is defined as follows: for simple objects $x,x',x_0,y,y',y_0\in\Irr(\Mm)$, and morphisms
		\agn{
				\diagram{y'\ar[r]^-f & [x,x']y} \quad\text{and}\quad \diagram{y_0\ar[r]^-g & G(x_0)}
		}
		in $\Mm$, we have 
		\[
		f.g=\delta_{y,y_0}\delta_{x,x_0} \diagram@C=2.5pc{
			(y' \ar[r]^-{f} &  [x,x']y \ar[r]^-{1g} & [x,x']G(x) \ar[r]^-{{G_2}_{[x,x'],x}}_-\sim & G([x,x']x) \ar[r]^-{G(\epsilon_{x,x'})} & G(x') })\,.
		\]
\enu
	\pf
		\bnu
			\item (\ref{item11:thm:with_wr})-(\ref{item14:thm:with_wr}) can be shown directly by definition of ${}^\wr\mu$, ${}^\wr\eta$, ${}^\wr\Delta$ and ${}^\wr\epsilon$ in \eqref{eq:with_wr}. To prove (\ref{item15:thm:with_wr}), one needs a tedious though direct computation, which can be carried out using the definition of ${}^\wr S$ in \eqref{eq:with_wr}, the definition of the antipode $S'$ on $\End(\Fff)$ given in \eqref{eq:reconstruct_antipode}, and the fact that for $x,y\in\Mm$, the left adjoint to the functor $L(F_0)=\Oplus_{x,y\in\Irr(\Mm)} [x,-]y\:\Mm\to\Mm$ is $\Oplus_{x,y\in\Irr(\Mm)}[-,y]^Rx$.
			\item Take ${}^\wr K$ as the composition of the monoidal equivalence $\widetilde{\Fff}$ in \Cref{cor:the_main_wff} and the functor $(\phi\inverse)^\ast\:\rep(\End(\Fff))\to\rep({}^\wr A_\Mm^\Cc),\; V\longmapsto {}_{\phi\inverse}V$. Then the statement follows from the expression of $\phi\inverse$ given in \Cref{cor:vs_transport}.
		\enu
	\epf
\end{thm}
\pf[Proof of \Cref{thm:main_equivalence}]
	Define a linear isomorphism 
	\eqn{\label{eq:iso_bw_wr_and_nowr}
	\begin{split}
		\psi\: A_\Mm^\Cc & \to {}^\wr A_\Mm^\Cc \\
		(\diagram{y' \ar[r]^-u & ay})\os (\diagram{ax\ar[r]^-s & x'}) & \longmapsto (\diagram{y' \ar[r]^-{u} & ay \ar[r]^-{s^\sharp} & [x,x']y })\,,
	\end{split}
	}
	where $s^\sharp\:a\to [x,x']$ is induced from $s\:ax\to x'$ via the adjunction \eqref{eq:adj}. It is not hard to establish the identities 
	\agn{
		\psi\circ \mu&={}^\wr\mu\circ (\psi\os\psi) \qquad\; \psi\circ \eta={}^\wr\eta\\
		(\psi\os\psi)\circ \Delta&={}^\wr\Delta\circ \psi\qquad\;\epsilon={}^\wr\epsilon\circ \psi\qquad\; \psi \circ S={}^\wr S\circ \psi\,.
	}
	Thus, $(A_\Mm^\Cc,\mu,\eta,\Delta,\epsilon,S)$ is a weak Hopf algebra, and $\psi$ becomes an isomorphism of weak Hopf algebras. The equivalence \eqref{eq:thm:main_equivalence} is then obtained by composing \eqref{eq:thm:with_wr} with 
	\[\psi^\ast\:\rep({}^\wr A_\Mm^\Cc)\isom\rep(A_\Mm^\Cc),\quad {}_{{}^\wr A_\Mm^\Cc}V\longmapsto {}_\psi V\,,
	\]
	where the $A_\Mm^\Cc$-action on ${}_\psi V$ is given by
	\[
		x.v=\psi(x).v,\;\forall x\in A_\Mm^\Cc, v\in V\,.
	\]
\epf
\begin{rmk}\label{rmk:rmk_on_presentation}
	Although \Cref{thm:with_wr} is only used as an intermediate step for proving the main theorem in this article, it provides an alternative presentation of the same weak Hopf algebra $A_\Mm^\Cc$ reconstructed from $\Fff$, namely ${}^\wr A_\Mm^\Cc$. This presentation is quite useful becauses it makes the use of internal homs, which is highly efficient in packaging data. As we will see, the other main result of this article, \Cref{thm:qt_main}, is also proved by first obtaining structures on ${}^\wr A_\Mm^\Cc$, and then ``transporting'' these structures to $A_\Mm^\Cc$ via $\psi$. Besides, we expect that the presentation ${}^\wr A_\Mm^\Cc$ will work well when $\Cc$ is a multi-fusion category, while the analog of the presentation $A_\Mm^\Cc$ in this case is more cumbersome to describe due to the fact that $\mone$ is no longer simple. 
	
	 In fact, the vector space $A_\Mm^\Cc$ is also isomorphic to
	 \[
	 	{}^{\wr\wr}A_\Mm^\Cc\defdtobe\Oplus_{x,x',y,y'\in\Irr(\Mm)}\Cc(\mone,[x,x'][y',y])\,,	\]
	 which looks more symmetric, via the isomorphisms
	\[
		\Mm(y',[x,x']y)\isom\Cc(\mone,[y',[x,x']y])\isom\Cc(\mone,[x,x'][y',y])\,.
	\]
	The presentation of the weak Hopf algebra structure of $A_\Mm^\Cc$ using ${}^{\wr\wr}A_\Mm^\Cc$ is left as an instructive exercise for the reader, and omitted in this article.
\end{rmk}

\begin{rmk}\label{rmk:on_multi_object}
	As anticipated by physicists \cite[\S 6]{Kitaev_Kong_2012}\cite[Eq. (3.13)]{Choi_Rayhaun_Zheng_2024b}, \Cref{thm:main_equivalence} admits a generalization. In \cite{Bai_Zhang}, we aim to prove this generalized result, which we briefly discuss as follows. For finite semisimple left $\Cc$-modules $\Mm,\Nn$, define an algebra 
	\[
		A_{\Mm,\Nn}^\Cc\defdtobe\Oplus_{\substack{x,x'\in\Irr(\Mm) \\y,y'\in\Irr(\Nn)}}\Nn(y',[x,x']y)\cong \Oplus_{\substack{x,x'\in\Irr(\Mm)\\ y,y'\in\Irr(\Nn)}}\Oplus_{a\in\Irr(\Cc)}\Mm(y',ay)\os \Mm(ax,x')\,.
	\]
	whose multiplication is similar to $A^\Cc_\Mm$. Let $\LMOD(\Cc)$ denote the set of finite semisimple left $\Cc$-modules. Then a ``generalized comultiplication''
	\[
		\Delta_{\Mm,\Kk,\Nn}\:A_{\Mm,\Nn}^\Cc\to A_{\Kk,\Nn}^\Cc\os A_{\Mm,\Kk}^\Cc
	\]
	can be defined for $\Kk\in\LMOD(\Cc)$ in a similar way as the comultiplication of $A_{\Mm,\Mm}^\Cc\equiv A_{\Mm}^\Cc$ \cite[\S 6]{Kitaev_Kong_2012}. However, the comultiplication $\Delta_{\Mm,\Kk,\Nn}$ does not fit into the definition of a weak Hopf algebra; to encompass the whole structure 
	\eqn{
		\{\{A_{\Mm,\Nn}^\Cc\}_{\Mm,\Nn\in\LMOD(\Cc)},\{\Delta_{\Mm,\Nn,\Kk}\}_{\Mm,\Nn,\Kk\in\LMOD(\Cc)}\}\,,
	}
	one needs to generalize the notion of weak Hopf algebras to a ``multi-object'' version, which corresponds to the ``weak'' version of \emph{dual $\kb$-linear Hopf category} in the sense of \cite{Batista_Caenepeel_Vercruysse_2016}. Moreover, to fully generalize \Cref{thm:main_equivalence}, the notion of ``representations'' needs to be reintepretated in this context. Such a mathematical theory can be developed completely in parrallel with that of weak Hopf algebras. In \cite{Bai_Zhang}, we develop such a theory by proving a Reconstruction Theorem for multi-object weak Hopf algebras, with which we show that a generalization of \Cref{thm:main_equivalence} to this setting is available.
		
	In \Cref{rmk:phy}, some physical meaning of the comultiplication $\Delta$ of $A_\Mm^\Cc$ in Levin-Wen models is discussed. Likewise, the generalized comultiplication $\Delta_{\Mm,\Kk,\Nn}$ controls the fusion of two topological excitations at the $\Kk$-$\Nn$ junction and the $\Mm$-$\Kk$-junction, respectively (cf. \cite[\S 6]{Kitaev_Kong_2012}), where $\Mm,\Kk,\Nn$ are all boundary labels of the $\Cc$-Levin-Wen model. In the context of conformal field theory, $\Delta_{\Mm,\Kk,\Nn}$ controls the OPE of boundary changing local operators in different Hilbert spaces \cite[\S 3.2.3]{Choi_Rayhaun_Zheng_2024b}. 
	
	Finally, we remark that the weak Hopf algebra $A_\Mm^\Cc$ can also be extended in an orthogonal direction to the direction discussed above. For details, we refer the reader to \cite[Eq. (3.5)]{Kong_2012}, and also to \cite[Eqs. (107)-(110)]{Lan_Wen_2014}, which is based on \cite{Kong, Kong_2012}.
\end{rmk}

\subsection{Example: reconstruction from an arbitrary fusion category}\label{subsec:right_regular}
Let $\Cc$ be a fusion category. We treat $\Cc$ as the regular right $\Cc$-module, or equivalently the left $\Cc^\rev$-module with action
\[
	\Cc^\rev\times\Cc\to\Cc,\quad (a,x)\longmapsto xa\,,
\]
where $\Cc^\rev$ is the fusion category with reversed tensor product. Then the following monoidal equivalence is well-known:
\eqn{\label{eq:C_funCCC}
\begin{split}
	\Cc& \to\fun_{\Cc^\rev}(\Cc,\Cc)\\
	w & \longmapsto w\os-\,.
\end{split}
}
Composing \eqref{eq:C_funCCC} with the monoidal equivalence given in \Cref{thm:main_equivalence}, we obtain
\begin{cor}\label{cor:right_regular}
There is a monoidal equivalence
\eqn{\label{eq:cor:right_regular}
\begin{split}
	\Cc & \to \rep(A_\Cc^{\Cc^\rev}) \\
	w & \longmapsto \Oplus_{y_0,x_0\in\Irr(\Cc)}\Cc(y_0,wx_0)\,,
\end{split}
}
where the action of $A_\Cc^{\Cc^\rev}=\Oplus_{y',y,x',x,a\in\Irr(\Cc)}\Cc(y',ya)\os \Cc(xa,x')$ on $\Oplus_{y_0,x_0\in\Irr(\Cc)}\Cc(y_0,wx_0)$ is given as follows: for $x,x',x_0,y,y',y_0,a\in\Irr(\Cc)$, and morphisms
\[
	\diagram{y'\ar[r]^-u & ya}, \qquad \diagram{xa\ar[r]^-s &  x'},\qquad \diagram{y_0\ar[r]^-g &  wx_0}
\]
in $\Cc$, the action of $u\os s$ on $g$ is given by
\[
	(u\os s).g=\delta_{x,x_0}\delta_{y,y_0}(\diagram{y' \ar[r]^-{u} & ya \ar[r]^-{g1} &wxa \ar[r]^-{1s} & wx'})\,.
\]
\end{cor}

Note that $A_\Cc^{\Cc^\rev}$ is reconstructed from the weak fiber functor 
$\diagram{\Cc\ar[r]^-\Omega & \fun(\Cc,\Cc)\ar[r]^-{\Vvv\Psi} &\vect}$,
where $\Omega$ is given by $x\mapsto x \os -$. Thus it is an explicit weak Hopf algebra reconstructed from $\Cc$ using the paradigm developed in \cite{Hayashi_1999} and introduced in \cite[Section 7.23]{Etingof_Gelaki_Nikshych_Ostrik_2015} (see also \Cref{rmk:to_funMM}). To be more specific, this algebra $A_\Cc^{\Cc^\rev}$ is precisely the one reconstructed in \cite{Hayashi_1999} (up to dual opposite), although the explicit form was not given there. 

The following is an example of $A_\Cc^{\Cc^\rev}$ and the equivalence \eqref{eq:cor:right_regular} when $\Cc=\vectG^\omega$ \cite[\S 4]{Hayashi_1999}\cite[\S 2.3.3]{Cordova_Holfester_Ohmori_2024}.
\begin{expl}\label{expl:regular_right_omega}
Let $G$ be a finite group and $\omega\in H^3(G,\kb^\times)$ be a 3-cocycle. Let $\vectG^\omega$ be the category of $G$-graded vector spaces whose associators are given by $\omega$; see for example \cite[Example 2.3.8]{Etingof_Gelaki_Nikshych_Ostrik_2015}. Without loss of generality, we assume that $\omega$ is \emph{normalized}, i.e., satisfies the condition
\[
	\omega(g,1,h)=\omega(1,g,h)=\omega(g,h,1)=1,\;\forall g,h\in G\,.
\]

Let us present the structure of the weak Hopf algebra $B_G^\omega\defdtobe A_{\vectG^\omega}^{(\vectG^\omega)^\rev}$. 
Suppose for a fusion category $\Cc$, we denote the subspace $\Oplus_{y',x'\in\Irr(\Cc)}\Cc(y',ya)\os \Cc(xa,x')\subset A_\Cc^{\Cc^\rev}$ for $a,x,y\in\Irr(\Cc)$ by $W_{a|y|x}$. Then in the case $\Cc=\vectG^\omega$, each $W_{a|y|x}$ is 1-dimensional, therefore the whole algebra $B_G^\omega$ is $|G|^3$-dimensional, where $|G|$ is the order of $G$. For $a,y,x\in G$, we set 
\[
	\ff_{a|y|x}\defdtobe\id_{ya}\os\id_{xa}\in W_{a|y|x}\,,
\]
so that $\{\ff_{a|y|x}\}_{a,y,x\in G}$ form a basis of $B_G^\omega$.

The weak Hopf algebra structure on $B_G^\omega$ is given as follows:
\bit
	\item The multiplication reads
	\[
		\ff_{a'|y'|x'}\cdot\ff_{a|y|x}=\delta_{y',ya}\delta_{x',xa}\frac{\omega(y,a,a')}{\omega(x,a,a')}\ff_{aa',y,x}\,.
	\]
	\item The unit reads
	\[
		1_{B_G^\omega}=\sum_{y,x\in G}\ff_{1|y|x}\,.
	\]
	\item The comultiplication reads
	\[
		\Delta(\ff_{a|y|x})=\sum_{z\in G}\ff_{a|y|z}\os\ff_{a|z|x}\,.
	\]
	\item The counit reads
	\[
		\epsilon(\ff_{a|y|x})=\delta_{x,y}\,.
	\]
	\item The antipode reads
	\[
		S(\ff_{a|y|x})=\frac{\omega(y,a,a\inverse)}{\omega(x,a,a\inverse)}\ff_{a\inverse|xa|ya}\,.
	\]
\eit

The equivalence $\vectG^\omega\isom\rep(B_G^\omega)$ sends a simple object $g\in\vectG^\omega$ to the $|G|$-dimensional vector space $\mathrm{span}\{\fh_{x_0}\}_{x_0\in G}$, on which the $B_G^\omega$-action is given by
\[
	\ff_{a|y|x}.\fh_{x_0}=\delta_{x,x_0}\delta_{y,gx}\omega(g,x,a)\fh_{xa},\;\forall a,y,x,x_0\in G\,.
\]

Finally, we remark that $B_G^\omega$ is not isomorphic to a groupoid algebra when $G$ is non-trivial, as any groupoid algebra must be cocommutative (cf. \Cref{expl:grpd}).
\end{expl}

\section{\texorpdfstring{The quasi-triangular structure on $A_\Cc^{\Cc\boxtimes\Cc^\rev}$}{The quasi-triangular structure on \$A\_C\^{}{CCrev}\$}}\label{sec:qt}
Let $\Cc$ be a fusion category. Then $\Cc$ can be viewed as a left $\Cc\boxtimes\Cc^\rev$-module via the evident action
\[
	\odot\:\Cc\boxtimes\Cc^\rev\times\Cc\to\Cc,\quad (a\boxtimes b,c)\longmapsto acb\,.
\]
Here $\boxtimes$ denotes the Deligne tensor product. By \Cref{thm:main_equivalence}, there is a monoidal equivalence
\eqn{\label{eq:sec:qt}
	\fun_{\Cc\boxtimes\Cc^\rev}(\Cc,\Cc)\isom\rep(A_\Cc^{\Cc\boxtimes\Cc^\rev})\,.
}

On the other hand, there is a monoidal equivalence $\fun_{\Cc\boxtimes\Cc^\rev}(\Cc,\Cc)\isom\Zz(\Cc)$, where $\Zz(\Cc)$ is the Drinfeld center of $\Cc$. Thus, we obtain a monoidal equivalence
\eqn{\label{eq:Z_C_equiv_A_CCC}
	\Zz(\Cc)\isom\rep(A_\Cc^{\Cc\boxtimes\Cc^\rev})\,.
}
The category $\Zz(\Cc)$ is a braided monoidal category. Moreover, it is well-known that braidings on the representation category of a weak Hopf algebra are in 1:1 correspondence with quasi-triangular structures on the algebra. In this section, we use the braiding on $\Zz(\Cc)$ to endow $A_\Cc^{\Cc\boxtimes\Cc^\rev}$ with a quasi-triangular structure, making the equivalence in  \eqref{eq:Z_C_equiv_A_CCC} a braided monoidal equivalence.
\subsection{Braidings and quasi-triangular structures}\label{subsec:qt}
Let $A$ be a weak Hopf algebra. In this subsection, we briefly recall the correspondence between braidings on $\rep(A)$ and quasi-triangular structures on $A$. For definition of braided monoidal categories, we refer the reader to \cite[Chapter 8]{Etingof_Gelaki_Nikshych_Ostrik_2015}.

For vector spaces $V$ and $W$, let $\tau_{V,W}\:V\os W\to W\os V$ denote the canonical braiding defined by $\tau_{V,W}(v\os w)=w\os v$ for $v\in V$ and $w\in W$. 
\begin{dfn}
A \newdef{quasi-triangular structure} on $A$ is an element 
\[
	\Rr\in(A\os A)\Delta(1)\equiv \{u\in A\os A\mid u\Delta(1)=u\}
\]
satisfying the conditions
\agn{
	\Rr\Delta(x) & =\Delta^{\cop}(x)\Rr,\;\forall x\in A\,,\\
	(\Delta\os\id)(\Rr) & =\Rr_{13}\Rr_{23}\,, \\
	(\id\os\Delta)(\Rr) & =\Rr_{13}\Rr_{12}\,,
}
such that there exists an element $\overline{\Rr}\in (A\os A)\Delta^\cop(1)\equiv\{u\in A\os A\mid u\Delta^\cop(1)=u\}$ with
\[
	\Rr\overline{\Rr}=\Delta^{\cop}(1) \qquad \overline{\Rr}\Rr=\Delta(1)\,.
\]
Here $\Delta^\cop$ denotes the comultiplication opposite to $\Delta$, and we use the standard notation $\Rr_{13}=\Rr^{(1)}\os 1\os  \Rr^{(2)}\in A\os A\os A$ and $\Rr_{23}=1\os \Rr^{(1)}\os\Rr^{(2)}\in A\os A\os A$, etc. (cf. \cite[\S VIII.2]{Kassel_1995}). We also call a quasi-triangular structure on $A$ an \newdef{$R$-matrix}.
\end{dfn}
Recall that for left $A$-modules $V$ and $W$, there are canonical maps 
\[
	r_{V,W}\:V\os W\to V\bos W \quad\text{and} \quad i_{V,W}\:V\bos W\to V\os W
\]
given by the retraction and the associated section, respectively, of the map $e_{V,W}$ defined by \eqref{eq:act_by_1_1_and_1_2}.\begin{thm}\label{thm:qt}
	\bnu
		\item \label{item1:thm:qt} Given an $R$-matrix $\Rr$, set 
		\[
		\widetilde{c^\Rr}_{V,W}\:V\os W\to W\os V,\quad v\os w\longmapsto (\Rr^{(2)}.w)\os (\Rr^{(1)}.v)
		\]
		for $V,W\in\rep(A)$. Then $c^\Rr$ is a braiding on $\rep(A)$ with 
		\[
			c^\Rr_{V,W}\defdtobe (\diagram@C=2.25pc{V\bos W \ar[r]^-{i_{V,W}} & V\os W \ar[r]^-{\widetilde{c^\Rr}_{V,W}} & W\os V \ar[r]^-{r_{W,V}} & W\bos V})\,.
		\]

		\item  \label{item2:thm:qt} Conversely, given a braiding $c$ on $\rep(A)$, define a map
	\[
		f^c\defdtobe(\diagram@C=2.25pc{A\os A\ar[r]^-{r_{A,A}} & A\bos A \ar[r]^-{c_{A,A}} & A\bos A \ar[r]^-{i_{A,A}} & A\os A \ar[r]^-{\tau_{A,A}} & A\os A})\,.
	\]
	Then $f^c(1\os 1)$ is an $R$-matrix on $A$.
		\item  \label{item3:thm:qt} Moreover, the assignments $\Rr\longmapsto c^\Rr$ and $c\longmapsto f^c(1\os 1)$ establish a bijection between quasi-triangular structures on $A$ and braidings on $\rep(A)$.
	\enu
	\pf
		The statements in \ref{item1:thm:qt} and \ref{item2:thm:qt} are proved in \cite[Proposition 5.2.2]{Nikshych_Vainerman_2000}. To prove \ref{item3:thm:qt}, it remains to check that $c=c^{f^c(1\os 1)}$ and $\Rr=f^{c^\Rr}(1\os 1)$. The first equality is also contained in \cite[Proposition 5.2.2]{Nikshych_Vainerman_2000}. The second equality is verified by the following calculation:
		\[
			f^{c^\Rr}(1\os 1)=\tau_{A,A}(\Delta(1)(\Rr^{(2)}\os\Rr^{(1)})\Delta^\cop(1))=\Delta^\cop(1)\Rr\Delta(1)=\Rr\,.
		\]
	\epf
\end{thm}
When $\rep(A)$ has a braiding $c$, we also define the \newdef{reduced $R$-matrix}, denoted by $\Rr_r$, to be the image of $1\os 1$ under the map 
\[
	\diagram{A\os A \ar[r]^-{r_{A,A}} & A\bos A \ar[r]^-{c_{A,A}} & A\bos A}\,.
\]
By \Cref{thm:qt}, the $R$-matrix corresponding to $c$ can be given by the reduced $R$-matrix via the formula 
\eqn{\label{eq:from_reduced}
	\Rr=\tau_{A,A}i_{A,A}(\Rr_r)\,.
}
Sometimes it is more convenient to first work out the reduced $R$-matrix, then apply \eqref{eq:from_reduced} to obtain the $R$-matrix.

\subsection{Computation of the quasi-triangular structure}
\label{subsec:qt_comp} 
In this subsection, we present a monoidal equivalence
\[
	\Zz(\Cc)\isom\rep(A_\Cc^{\Cc\boxtimes\Cc^\rev})\,,
\]
and then apply \Cref{thm:qt} to endow $A_\Cc^{\Cc\boxtimes\Cc^\rev}$ with a quasi-triangular structure.
 
Let us first recall some basic facts on the Drinfeld center. We adopt the definition given in \cite[Definition 7.13.1]{Etingof_Gelaki_Nikshych_Ostrik_2015}. In particular, objects in $\Zz(\Cc)$ are pairs $(z,\gamma_{-,z})$, where $z$ is an object of $\Cc$, and 
\[
	\gamma_{-,z}=\{\gamma_{w,z}\:wz\isom zw\}_{w\in\Cc}
\]
is a \emph{half-braiding} on $z$, a family of isomorphisms natural in $w$ and satisfying certain constraints. Morphisms in $\Zz(\Cc)$ are morphisms in $\Cc$ that are compatible with the half-braiding.

We will introduce two aspects of the Drinfeld center, both of which are categorifications of facts related to the center $Z(A)$ of a $\kb$-algebra $A$. The first fact is that $Z(A)$ is a commutative algebra. This is categorified in that the Drinfeld center $\Zz(\Cc)$ is a braided monoidal category. Specifically, for $(z,\gamma_{-,z})$ and $(z',\gamma'_{-,z'})$, the braiding is given by 
\eqn{\label{eq:braiding_on_Z}
	c_{(z,\gamma_{-,z}),(z',\gamma'_{-,z'})}\defdtobe \gamma'_{z,z'}\:z\os z'\isom z'\os z\,.
}

The other fact on $Z(A)$ is that the map
\[
	Z(A)\to\Hom_{A\os A^\op}(A,A), \quad a\mapsto a\cdot-
\]
defines an algebra isomorphism from the center to the algebra $\Hom_{A\os A^\op}(A,A)$, which consists of $A$-$A$-bimodule maps from $A$ to itself. This is categorified by the following well-known lemma:
\begin{lem}
There is an equivalence of monoidal categories
\eqn{\label{eq:drinfeld_center}
\begin{split}
	\Zz(\Cc) & \to\fun_{\Cc\boxtimes\Cc^\rev}(\Cc,\Cc) \\
		(z,\gamma_{-,z}) & \longmapsto (z\os-,(z\os -)_2)\,,
\end{split}
}
where the $\Cc\boxtimes\Cc^\rev$-module structure $(z\os -)_2$ is given by 
\[
	{(z\os-)_2}_{a\boxtimes b,c}\:\diagram{azcb\ar[r]^-{\gamma_{a,z}1} & zacb},\;\forall a,b,c\in\Cc\,.
\]
\end{lem}

Now, composing the equivalence \eqref{eq:drinfeld_center} with \eqref{eq:sec:qt}, we have
\begin{cor}\label{cor:Drinfeld_center}
There is a monoidal equivalence
\eqn{\label{eq:qt_drinfeld_center}
\begin{split}
	\Zz(\Cc) & \to\rep(A_\Cc^{\Cc\boxtimes\Cc^\rev}) \\
	(z,\gamma_{-,z}) & \longmapsto \Oplus_{x_0,y_0\in\Irr(\Cc)}\Cc(y_0,zx_0)\,,
\end{split}
}
where the action of \(A_\Cc^{\Cc\boxtimes\Cc^\rev}=\Oplus_{y,y',x,x',a,b\in\Irr(\Cc)}\Cc(y',ayb)\os\Cc(axb,x')\) on $\Oplus_{x_0,y_0\in\Irr(\Cc)}\Cc(y_0,zx_0)$ is given as follows: for $y,y',y_0,x,x',x_0,a,b\in\Irr(\Cc)$, and morphisms 
\[
\diagram{y'\ar[r]^-u & ayb}, \qquad \diagram{axb\ar[r]^-s & x'}, \qquad \diagram{y_0\ar[r]^-g & zx_0}
\]
in $\Cc$, the action of $u\os s$ on $g$ is given by
\[
	(u\os s).g=\delta_{x,x_0}\delta_{y,y_0}(\diagram{
		y' \ar[r]^-u & ayb \ar[r]^-{1g} & azxb \ar[r]^-{\gamma_{a,z}1} & zaxb \ar[r]^-{1s} & zx'
	})\,.
\]
\end{cor}
By \Cref{thm:qt}, the weak Hopf algebra $A_\Cc^{\Cc\boxtimes\Cc^\rev}$ possesses a quasi-triangular structure which corresponds to the braiding \eqref{eq:braiding_on_Z} on $\Zz(\Cc)$. To present this quasi-triangular structure, several additional preparations are needed. 

First, let us denote 
\(
	U_{a|b|y'|y|x'|x}\defdtobe\Cc(y',ayb)\os\Cc(axb,x')\subset A_\Cc^{\Cc\boxtimes\Cc^\rev}
\)
for $a,b,y',y,x',x\in\Irr(\Cc)$. Then note that there is an obvious inclusion
\[
	\iota_1\:A_\Cc^{\Cc^\rev}\to A_\Cc^{\Cc\boxtimes\Cc^\rev}
\]
sending $(\diagram{y'\ar[r]^-u & yb})\os(\diagram{xb \ar[r]^-s & x'})$ for $b,y',y,x',x\in\Irr(\Cc)$ to $u\os s\in U_{\mone|b|y'|y|x'|x}$. Similarly, there is an inclusion 
\[
	\iota_2\:A_\Cc^\Cc\to A_\Cc^{\Cc\boxtimes\Cc^\rev}\,.
\]
Let $\psi_1\:A_\Cc^{\Cc^\rev}\to {}^\wr A_\Cc^{\Cc^\rev}$ and $\psi_2\:A_\Cc^\Cc\to{}^\wr A_\Cc^\Cc$ defined by \eqref{eq:iso_bw_wr_and_nowr}.

Secondly, define a linear map
\eqn{\label{eq:copairing}
	\Theta\:\kb\to {}^\wr A_\Cc^{\Cc^\rev}\os {}^\wr A_\Cc^{\Cc}
}
by setting $\Theta(1)$ to be 
\[
	\begin{multlined}
		\sum_{w,x,y,z\in\Irr(\Cc)}\sum_{\alpha=1}^{n^{x,y,z}_w}(\diagram@C=3pc{w\ar[r]^-{I^{x,y,z,\alpha}_w} & yx^Rz\ar[r]^-1_-\sim & [x,z]_{\Cc^\rev}\os^\rev y}) \\
		\os (\diagram@C=3.5pc{y\ar[r]^-{1\coev} & yx^Rzz^Lx \ar[r]^-{P^{x,y,z,\alpha}_w 1} & wz^Lx \ar[r]^-1_-\sim & [z,w]_\Cc\os x})\,,
	\end{multlined}
\]
where $[x,z]_{\Cc^\rev}$ and $[z,w]_\Cc$ are computed in \Cref{expl:ih}, and for each $x,y,z\in\Irr(\Cc)$, the maps $I^{x,y,z,\alpha}_{w}\:w\to yx^Rz$ and $P^{x,y,z,\alpha}_{w}\:yx^Rz\to w$ are the inclusions and projections, respectively, in the direct sum decomposition
\[
	yx^Rz\cong \oplus_{w\in\Irr(\Cc)}w^{\oplus n^{x,y,z}_{w}}
\]
for $w\in\Irr(\Cc)$ and $\alpha=1,\cdots,n^{x,y,z}_{w}$.

Lastly, we need the following computation of the internal homs in $\Cc\boxtimes\Cc^\rev$:
\begin{lem}\label{lem:reg_ih}
	For $x,y\in\Cc$, let $[x,y]_\sboxtimes$ denote the internal hom from $x$ to $y$ in $\Cc\boxtimes\Cc^\rev$. Then we have
	\[
		[x,y]_\sboxtimes=\oplus_{i\in\Irr(\Cc)}yi^Lx^L\boxtimes i=\oplus_{i\in\Irr(\Cc)}i^Lx^L\boxtimes iy\,.
	\]
	\pf
		To show the first equality, we observe that there are natural isomorphisms
		\[
			\begin{multlined}\Cc(axb,y)\cong\Cc(b,x^Ra^Ry)\cong\Oplus_{i\in\Irr(\Cc)}\Cc(i,x^Ra^Ry)\os\Cc(b,i)\cong\Oplus_{i\in\Irr(\Cc)}\Cc(y^Lax,i^L)\os\Cc(b,i)\\
			\cong\Oplus_{i\in\Irr(\Cc)}\Cc(a,yi^Lx^L)\os\Cc(b,i)\cong\Cc\boxtimes\Cc^\rev(a\boxtimes b,\oplus_{i\in\Irr(\Cc)}yi^Lx^L\boxtimes i)
			\end{multlined}
		\]
		for $a,b\in\Cc$. The second equality is proved similarly.
	\epf
\end{lem}
\begin{thm}\label{thm:qt_main}
	The quasi-triangular structure $\Rr\in A_\Cc^{\Cc\boxtimes\Cc^\rev}\os A_\Cc^{\Cc\boxtimes\Cc^\rev}$ corresponding to the braiding on $\Zz(\Cc)$ via the equivalence \eqref{eq:qt_drinfeld_center} is given by 
	\eqn{\label{eq:qt:1}
		\Rr=(\iota_1\psi_1\inverse\os\iota_2\psi_2\inverse)\Theta(1)\,,
	}
	where $\Theta$ is defined in \eqref{eq:copairing}. Explicitly, we have
	\eqn{\label{eq:qt:2}
		\Rr=\sum_{a,b,w,x,y,z\in\Irr(\Cc)}\Rr_{a,b,w,x,y,z}\,,
	}
	where $\Rr_{a,b,w,x,y,z}$ is given in the following steps:
	
	\bnu[$1^\circ$]
		\item Choose a direct sum decomposition
	\[
		yx^Rz\cong \oplus_{w'\in\Irr(\Cc)}w'^{n_{w'}}
	\]
	with inclusions $I^\alpha_{w'}\:w'\to yx^Rz$ and projections $P^\alpha_{w'}\:yx^Rz\to w'$ for $\alpha=1,\cdots,n_{w'}$ and $w'\in\Irr(\Cc)$. Choose a direct sum decomposition 
	\[
		x^Rz\cong \oplus_{b'\in\Irr(\Cc)}b'^{n'_{b'}}
	\]
	with inclusions $I'^\beta_{b'}\:b'\to x^Rz$ and projections $P'^\beta_{b'}\:x^Rz\to b'$ for $\beta=1,\cdots,n'_{b'}$ and $b'\in\Irr(\Cc)$. Choose a direct sum decomposition
	\[
		wz^L\cong \oplus_{a'\in\Irr(\Cc)}a'^{n''_{a'}}
	\]
	with inclusions $I''^\gamma_{a'}\:a'\to wz^L$ and projections $P''^\gamma_{a'}\:wz^L\to a'$ for $\gamma=1,\cdots,n''_{a'}$ and $a'\in\Irr(\Cc)$.
		\item Then there is 
		\[
			\Rr_{a,b,w,x,y,z}=\sum_{\alpha=1}^{n_w}\sum_{\beta=1}^{n'_b}\sum_{\gamma=1}^{n''_a}g^{\alpha,\beta,\gamma}_1\os g^{\alpha,\beta,\gamma}_2\os g^{\alpha,\beta,\gamma}_3\os g^{\alpha,\beta,\gamma}_4\in U_{\mone|b|w|y|z|x}\os U_{a|\mone|y|x|w|z}\,, 
		\]
		where
		\agn{
			g^{\alpha,\beta,\gamma}_1 & = (\diagram{w\ar[r]^-{I^\alpha_w} & yx^Rz \ar[r]^-{1P'^\beta_b} & yb}) \\
			g^{\alpha,\beta,\gamma}_2 & = (\diagram{xb \ar[r]^-{1I'^\beta_b} & xx^Rz \ar[r]^-{\ev1} & z})\\
			g^{\alpha,\beta,\gamma}_3 & = (\diagram{y \ar[r]^-{1\coev} & yx^Rzz^Lx \ar[r]^-{P^\alpha_w1} & wz^Lx \ar[r]^-{P''^\gamma_a1} &  ax} )\\
			g^{\alpha,\beta,\gamma}_4 & = (\diagram{az \ar[r]^-{I''^\gamma_a1} & wz^Lz \ar[r]^-{1\ev} & w})\,.
		}
	\enu
	\pf
		We only show the explicit form \eqref{eq:qt:2} is true; it is easy to deduce \eqref{eq:qt:1} from \eqref{eq:qt:2}. 
		
		Let $A\defdtobe A_\Cc^{\Cc\boxtimes\Cc^\rev}$. We use the notations from \Cref{sec:2}. Specifically, let $\Fff\:\fun_{\Cc\boxtimes\Cc^\rev}(\Cc,\Cc)\to\vect$ represent the weak fiber functor we constructed in \Cref{subsec:the_wff}. Let $F_0\:\Cc\to\Cc$ refer to the functor sending each simple object $x\in\Irr(\Cc)$ to $\oplus_{y\in\Irr(\Cc)}y$. Let $L\:\fun(\Cc,\Cc)\to\fun_{\Cc\boxtimes\Cc^\rev}(\Cc,\Cc)$ denote the left adjoint of the forgetful functor $\fun_{\Cc\boxtimes\Cc^\rev}(\Cc,\Cc)\to\fun(\Cc,\Cc)$ given in \Cref{lem:Gamma_admit_adj}, so that $G_0\defdtobe L(F_0)$ reads
\agn{
	G_0\:\Cc & \to \Cc \\
	x'\in\Irr(\Cc) & \mapsto \oplus_{x,y\in\Irr(\Cc)}[x,x']_\sboxtimes\odot y\,.
}
Let ${}^\wr 1\in \Fff(G_0)$ be defined in \eqref{eq:wr_1}. Finally, we use $\widetilde{\Fff}\:\fun_{\Cc\boxtimes\Cc^\rev}(\Cc,\Cc)\isom\rep(A)$ to represent the comparison functor \eqref{eq:thm:main_equivalence} in \Cref{thm:main_equivalence}.
		
		Let us first obtain the reduced $R$-matrix $\Rr_r$ of $A$ (see \Cref{subsec:qt}). Note that $\widetilde{\Fff}(G_0)=A$. Therefore $\Rr_r$ can be viewed as an element in 
\[
\begin{multlined}
\widetilde{\Fff}(G_0)\bos \widetilde{\Fff}(G_0)\cong\widetilde{\Fff}(G_0G_0)
=\Oplus_{y,z\in\Irr(\Cc)}\Cc(y,G_0G_0(z))\,.
\end{multlined}
\]

By definition, $\Rr_r$ is the image at ${}^\wr 1\os {}^\wr 1\in\Fff(G_0)\os\Fff(G_0)$ under the map
\[
\begin{multlined}
	\diagram@=3.75pc{\Fff(G_0)\os\Fff(G_0) \ar[r]^-{{\Fff_2}_{G_0,G_0}} & \Fff(G_0G_0) \ar[r]^-1_-\sim & \Oplus_{y,z\in\Irr(\Cc)} \Cc(y,G_0G_0(z))} \\
	\diagram@C=7.5pc{\ar[r]^-{\Oplus_{z\in\Irr(\Cc)}((c'_{G_0,G_0})_{z})_\ast} & \Oplus_{y,z\in\Irr(\Cc)} \Cc(y,G_0G_0(z))}\,,
\end{multlined}
\]
where $c'$ is the braiding on $\fun_{\Cc\boxtimes\Cc^\rev}(\Cc,\Cc)$ induced from the braiding on $\Zz(\Cc)$. After a tedious though straightforward computation, one finds that $\Rr_r=\sum_{y,z,x\in\Irr(\Cc)}{\Rr_r}_{x,y,z}$, where ${\Rr_r}_{x,y,z}\in\Cc(y,G_0G_0(z))$ is given by
\agn{
	y&\diagram@C=3.5pc{\ar[r]^-{\coev 1\coev} & xx^Lyx^Rzz^Lx} \\
	&\diagram@C=3.5pc{\ar@{^(->}[r]^-{(\ast)} & \oplus_{i\in\Irr(\Cc)}xx^Lyx^Rzi^Lz^Lxi}\\
	&\diagram@C=3.5pc{\ar[r]^-\sim_-{\text{\Cref{lem:reg_ih}}} & [z,xx^Lyx^Rz]_\sboxtimes\odot x} \\
	&\diagram@C=3.5pc{\ar@{^(->}[r] & G_0(xx^Lyx^Rz)}\\
	&\diagram@C=3.5pc{\ar@{^(->}[r]^-{(\ast\ast)} & G_0(\oplus_{j\in\Irr(\Cc)}j^Lx^Lyjz)}\\
	&\diagram@C=3.5pc{\ar[r]^-\sim_-{\text{\Cref{lem:reg_ih}}} & G_0([x,z]_\sboxtimes\odot y)} \\
	&\diagram@C=3.5pc{\ar@{^(->}[r] & G_0G_0(z)}\,.
}
	 Here, $(\ast)$ refers to the inclusion into the component $i=\mone$, while $(\ast\ast)$ refers to the inclusion into the component $j=x^R$.
	
	From \Cref{subsec:qt}, to obtain the $R$-matrix $\Rr\in A\os A$, one needs to apply the map
	\[
		\begin{multlined}
		\diagram@C=4.5pc{\displaystyle\Oplus_{y,z\in\Irr(\Cc)}\manualformatting\hskip -0.5em\Cc(y,G_0G_0(z))\ar[r]^-{{\Fff_{-2}}_{G_0,G_0}} &\displaystyle \Oplus_{y,w\in\Irr(\Cc)}\hskip -0.5em\Cc(y,G_0(w))\os\Oplus_{w',z\in\Irr(\Cc)}\hskip -0.5em\Cc(w',G_0(z)) } \\
		\diagram@C=5pc{\ar[r]^-{\tau_{\Fff(G_0),\Fff(G_0)}} &\displaystyle \Oplus_{w',z\in\Irr(\Cc)}\hskip -0.5em\Cc(w',G_0(z))\os\Oplus_{y,w\in\Irr(\Cc)}\hskip -0.5em\Cc(y,G_0(w)) }\\
		\diagram@C=4.5pc{\ar[r]^-1_-\sim & {}^\wr A_\Cc^{\Cc\boxtimes\Cc^\rev}\os {}^\wr A_\Cc^{\Cc\boxtimes\Cc^\rev} \ar[r]^-{\psi\inverse\os\psi\inverse} & A\os A}
		\end{multlined}
	\]
	to $\Rr_r$, where $\psi$ is given in \eqref{eq:iso_bw_wr_and_nowr}. As an intermediate step, we obtain that
	\[
		{}^\wr\Rr\defdtobe\tau_{\Fff(G_0),\Fff(G_0)}{\Fff_{-2}}_{G_0,G_0}(\Rr_r)=\sum_{x,y,z\in\Irr(\Cc)}{}^\wr\Rr_{x,y,z}\,.
	\]
	Here
	\[
		{}^\wr\Rr_{x,y,z}=\sum_{w\in\Irr(\Cc)}\sum_{\alpha=1}^{n_w} f_{w,1}^\alpha\os f_{w,2}^\alpha\,,
	\]
	where 
	\agn{
	f_{w,1}^\alpha & = (\diagram{w \ar[r]^-{I^\alpha_w} & xx^Lyx^Rz \ar@{^(->}[r] & \oplus_{j\in\Irr(\Cc)}j^Lx^Lyjz \ar[r]^-1_-\sim & [x,z]_\sboxtimes\odot y}) \\
		f_{w,2}^\alpha & =(\diagram@C=2pc{y \ar[rr]^-{\coev 1\coev} && xx^Lyx^Rzz^Lx \ar[r]^-{P^\alpha_w1} & wz^Lx \ar@{^(->}[r] & \oplus_{i\in\Irr(\Cc)}wi^Lz^Lxi \ar[r]^-1_-\sim & [z,w]_\sboxtimes\odot x })\,,
	}
	and $I_w^\alpha$ and $P_w^\alpha$ are respectively the inclusions and the projections of a direct sum decomposition
	\[
		xx^Lyx^Rz\cong\oplus_{w\in\Irr(\Cc)}w^{\oplus n_w}
	\]
	for $\alpha=1,\cdots,n_w$.

	Finally, one can check that 
	\[
		(\psi\inverse\os\psi\inverse)({}^\wr\Rr)
	\]
	is precisely the element given in the \rhs \manualformatting\hspace{-0.1ex} of \eqref{eq:qt:2}.
	\epf
\end{thm}
\begin{rmk}
	We give a long remark on how $A^{\Cc\boxtimes\Cc^\rev}_\Cc$ can be viewed as the Drinfeld double of $A^{\Cc^\rev}_\Cc$, inspired by \cite{Jia_Tan_Kaszlikowski_2024}. In view of the Reconstruction Theorem, the \newdef{Drinfeld double} (also called \newdef{quantum double}) $D(A)$ of a finite-dimensional weak Hopf algebra $A$ is the weak Hopf algebra reconstructed from the weak fiber functor 
	\[
		\diagram{\Zz(\rep(A))\ar[r]^-G & \rep(A) \ar[r]^-{\Fff^A} & \vect}\,,
	\]
	where $G$ is the forgetful functor. One can then see that $A^{\Cc\boxtimes\Cc^\rev}_\Cc$ is the Drinfeld double of $A^{\Cc^\rev}_\Cc$ since $A^{\Cc\boxtimes\Cc^\rev}_\Cc$ is reconstructed from the weak fiber functor
	\[
		\diagram@C=2.25pc{\fun_{\Cc\boxtimes\Cc^\rev}(\Cc,\Cc)\ar[r]^-\cong_-{\eqref{eq:drinfeld_center}} & \Zz(\Cc)\ar[r]^-G & \Cc \ar[r]^-\cong_-{\eqref{eq:C_funCCC}} & \fun_{\Cc^\rev}(\Cc,\Cc) \ar[r]^-{\Gamma} & \fun(\Cc,\Cc) \ar[r]^-{\Vvv\Psi} & \vect}\,.
	\] 
	
	It is interesting to explicitly construct the isomorphism $A^{\Cc\boxtimes\Cc^\rev}_\Cc\cong D(A^{\Cc^\rev}_\Cc)$. This involves a pairing between $A^\Cc_\Cc$ and $A^{\Cc^\rev}_\Cc$. Let $B,A$ be weak Hopf algebras and 
	\[
		\<,\>\:B\os A\to\kb
	\]
	be a non-degenerate pairing satisfying
	\[
		\<b,a_{(1)}\>\<b',a_{(2)}\>=\<bb',a\>\qquad \<1_B,a\>=\epsilon_A(a)
	\]
	\[
		\<b_{(1)},a\>\<b_{(2)},a'\>=\<b,a'a\> \qquad \<b,1_A\>=\epsilon_B(b)
	\]
	for any $a,a'\in A,b,b'\in B$. Note that these conditions equivalently say that the pairing induce a weak Hopf algebra isomorphism $B\isom (A^\ast)^\cop$; in particular, for any weak Hopf algebra $A$, such a pairing for $A$ exists. Given a pairing satsifying the above conditions, the explicit form of $D(A)$ can be defined as follows. As a vector space, $D(A)\defdtobe B\os A/I$, where $I$ is the subspace generated by 
	\[
		b\os xa-b\<{1_B}_{(1)},x\>{1_B}_{(2)}\os a,\quad x\in A^l\,;
	\]
	\[
		b\os ya-b{1_B}_{(1)}\<{1_B}_{(2)},y\>\os a,\quad y\in A^r\,.
	\]
	The multiplication is given by
	\[
		[b'\os a']\cdot[b\os a]=\<b_{(1)},a_{(1)}\>\<b_{(3)},S\inverse(a_{(3)})\>[b'b_{(2)}\os a'_{(2)}a]\,.
	\]
	The unit is given by $[1_B\os 1_A]$. The comultiplication reads
	\[
		\Delta([b\os a])=[b_{(1)}\os a_{(1)}]\os[b_{(2)}\os a_{(2)}]\,.
	\]
	The counit reads 
	\[
		\epsilon([b\os a])=\<b,\epsilon^{rr}(a)\>\,.
	\]
	The $R$-matrix is given by 
	\[
		\Rr=[1_B\os a_i]\os [b_i\os 1_A]\,,	
	\]
	where $1\longmapsto \sum_i a_i\os b_i$ is the copairing associated with the pairing $\<,\>$.
	Note that there are different conventions regarding the definition of Drinfeld double. Our definition adheres to the one given in \cite[\S IX.4]{Kassel_1995} when $A$ is a Hopf algebra, and is different from \cite[\S 5.3]{Nikshych_Vainerman_2000}.\footnote{When $B$ is identified with $(A^\ast)^\cop$, the map
	\[
		S^\ast\os\id\:D_{\mathrm{NV}}(A)\to D(A)
	\]
	provides an isomorphism of weak Hopf algebras from the Drinfeld double in \cite{Nikshych_Vainerman_2000} to $D(A)$ presented here. Identified with this isomorphism, the relation between our $R$-matrix $\Rr$ and the $R$-matrix $\Rr_{\mathrm{NV}}$ given in \cite{Nikshych_Vainerman_2000} is given by
	\[
		\Rr=\tau_{D(A),D(A)}(\overline{\Rr}_{\mathrm{NV}})\,,
	\]
	where $\overline{\Rr}_{\mathrm{NV}}$ is the unique element satisfying $\overline{\Rr}_{\mathrm{NV}}\Rr_{\mathrm{NV}}=\Delta(1)$ and $\Rr_{\mathrm{NV}} \overline{\Rr}_{\mathrm{NV}}=\Delta^\cop(1)$. }
	
	In the case $A=A^{\Cc^\rev}_\Cc$, one can take $B=A^\Cc_\Cc$ with the non-degenerate pairing
	\eqn{\label{eq:pairing}
		A^\Cc_\Cc\os A^{\Cc^\rev}_\Cc \to\kb
	}
	sending $(\diagram{y_1'\ar[r]^-{u_1} & a_1y_1})\os(\diagram{a_1x_1\ar[r]^-{s_1} & x_1'})\os(\diagram{y_2'\ar[r]^-{u_2} & y_2a_2})\os(\diagram{x_2a_2\ar[r]^-{s_2} & x_2'})$ to 
	\[
		\delta_{y_2,y_1'}\delta_{x_2,y_1}\delta_{x_2',x_1}\delta_{y_2',x_1'}\Lambda_{y_2'}(\diagram{y_2'\ar[r]^-{u_2} & y_2a_2\ar[r]^-{u_11} & a_1y_1a_2 \ar[r]^-{1s_2} & a_1x_2' \ar[r]^-{s_1} & x_1'})\,.
	\]
	The associated copairing is $\Theta$ given in \eqref{eq:copairing}. 
		
	Then, there is an isomorphism $D(A_\Cc^{\Cc^\rev})\equiv A^\Cc_\Cc\os A^{\Cc^\rev}_\Cc/I\isom A^{\Cc\boxtimes\Cc^\rev}_\Cc$ induced from the map 
	\[
		\sharp\:A^\Cc_\Cc\os A^{\Cc^\rev}_\Cc\to A^{\Cc\boxtimes\Cc^\rev}_\Cc\,,
	\]
	where $\sharp$ sends $(\diagram{y_1'\ar[r]^-{u_1} & a_1y_1}) \os (\diagram{a_1x_1\ar[r]^-{s_1} & x_1'}) \os (\diagram{y_2'\ar[r]^-{u_2} & y_2a_2}) \os (\diagram{x_2a_2\ar[r]^-{s_2} & x_2'})$ to 
	\[
		\delta_{y_2',y_1}\delta_{x_2',x_1}(\diagram{y_1'\ar[r]^-{u_1} & a_1y_1 \ar[r]^-{1u_2} & a_1y_2a_2})\os (\diagram{a_1x_2a_2\ar[r]^-{1s_2} & a_1x_2' \ar[r]^-{s_1} & x_1'})\,.
	\]
	Note that the non-degenerate pairing \eqref{eq:pairing} in particular shows that $A_\Cc^\Cc\cong((A_\Cc^{\Cc^\rev})^\ast)^\cop$ as weak Hopf algebras.
	
	More generally, if $\Mm$ is a Morita equivalence between fusion categories $\Cc$ and $\Dd$, i.e., $\Dd=\fun_\Cc(\Mm,\Mm)^\rev$ (cf. \cite[Definition 7.12.17]{Etingof_Gelaki_Nikshych_Ostrik_2015}), then there is a similar pairing 
	\[
		A_\Mm^\Cc\os A_\Mm^{\Dd^\rev}\to\kb
	\]
	as \eqref{eq:pairing} exhibiting $A^\Cc_\Mm$ as the coopposite of the dual of $A_\Mm^{\Dd^\rev}$, and $A_\Mm^{\Cc\boxtimes\Dd^\rev}$ is the Drinfeld double of $A_\Mm^{\Dd^\rev}$. This pairing is inspired by \cite[Remark 5.2]{Jia_Tan_Kaszlikowski_2024}, although our pairings are in disagreement.
	
	For a more thorough discussion of some of the content in this remark, we refer the reader to \cite{Jia_Tan}, which independently work out the pairing \eqref{eq:pairing} and the other observations in this remark.
\end{rmk}
The following is an example of the quasi-triangular weak Hopf algebra $A_\Cc^{\Cc\boxtimes\Cc^\rev}$ when $\Cc=\vectG^\omega$.

\begin{expl}\label{expl:qt_vectG}
	We assume $\omega$ is normalized as in \Cref{expl:regular_right_omega}. For any fusion category $\Cc$, we denote the subspace $\Oplus_{y',x'\in\Irr(\Cc)}\Cc(y',(ay)b)\os\Cc((ax)b,x')\subset A_\Cc^{\Cc\boxtimes\Cc^\rev}$ by $V_{a|b|y|x}$. Then in the case $\Cc=\vectG^\omega$, each $V_{a|b|y|x}$ is 1-dimensional, thus the whole algebra $A_G^\omega\defdtobe A_{\vectG^\omega}^{\vectG^\omega\boxtimes(\vectG^\omega)^\rev}$ is $|G|^4$-dimensional, where $|G|$ is the order of $G$. For $a,b,y,x\in G$, we set
\[
	\fe_{a|b|y|x}\defdtobe \id_{ayb}\os\id_{axb}\in V_{a|b|y|x}\,,
\]
so that $\{\fe_{a|b|y|x}\}_{a,b,y,x\in G}$ form a basis of $A_G^\omega$.

The quasi-triangular weak Hopf algebra structure on $A_G^\omega$ is given as follows:
\bit
	\item The multiplication reads
	\[
		\fe_{a'|b'|y'|x'}\cdot \fe_{a|b|y|x} = \delta_{y',ayb}\delta_{x',axb}\frac{\omega(a',a,x)}{\omega(a',a,y)}\frac{\omega(a',ax,b)}{\omega(a',ay,b)}\frac{\omega(a'ay,b,b')}{\omega(a'ax,b,b')}\fe_{a'a|bb'|y|x}\,.
	\]
	\item The unit reads
	\[
		1_{A_G^\omega}=\sum_{y,x\in G}\fe_{1|1|y|x}\,.
	\]
	\item The comultiplication reads
	\[
		\Delta(\fe_{a|b|y|x})=\sum_{z\in G}\fe_{a|b|y|z}\os\fe_{a|b|z|x}\,.
	\]
	\item The counit reads
	\[
		\epsilon(\fe_{a|b|y|x})=\delta_{x,y}\,.
	\]
	\item The antipode reads
	\[
		S(\fe_{a|b|y|x})=\frac{\omega(y,b,b\inverse)}{\omega(x,b,b\inverse)}\frac{\omega(a,y,b)}{\omega(a,x,b)}\frac{\omega(a,a\inverse,axb)}{\omega(a,a\inverse,ayb)}\fe_{a\inverse|b\inverse|axb|ayb}\,.
	\]
	\item The quasi-triangular structure reads 
	\[
		\Rr=\sum_{a,b,z\in G}\omega(a,z,b)\inverse\fe_{1|b|az|z}\os\fe_{a|1|z|zb}\,.
	\]
\eit
Finally, we remark that $A_G^\omega$ is not isomorphic to a groupoid algebra when $G$ is non-trivial, as any groupoid algebra must be cocommutative (cf. \Cref{expl:grpd}).
\end{expl}

\appendix
\section{Supplementary proofs}
\subsection{\texorpdfstring{Proof of \Cref{lem:EndF_times_EndG}}{Proof of Lemma \ref{lem:EndF_times_EndG}}}\label{subsec:pf:lem:EndF_times_EndG}
\pf
		Let $U^{\End(i)}\:\rep(\End(i))\to\vect$ denote the forgetful functor for $i=F,G$. Let 
	\[\widetilde{F}\:\Aa\to\rep(\End(F))\quad\text{and}\quad\widetilde{G}\:\Bb\to\rep(\End(G))\]
be the comparison functors as in \eqref{eq:comparison}, which are equivalences by \Cref{thm:monadicity}. Then, we have a strictly commutative diagram of functors
		\[
			\diagram@R=3pc{
				\Aa\times\Bb\ar[r]^-{F\times G}  \ar[d]_{\widetilde{F}\times\widetilde{G}} & \vect\times\vect \ar[r]^-{\os} & \vect \\
				\rep(\End(F))\times\rep(\End(G)) \ar[ru]_-{\hskip 2em U^{\End(F)}\times U^{\End(G)}}
			}\,.
		\]
		
		We define a map 
		\[J_1\:\End(F)\os\End(G)\to  \End(\os(U^{\End(F)}\times U^{\End(G)}))\]
		by setting $J_1(\alpha\os\beta)_{V,W}$ to be the map defined by 
		\[V\os W\to V\os W,\quad v\os w\longmapsto \alpha.v\os \beta.w\]
		for $V\in\rep(\End(F)),W\in\rep(\End(G))$ and $\alpha\in\End(F),\beta\in\End(G)$. We also set
		\[
			J_2\:\End(\os(U^{\End(F)}\times U^{\End(G)}))\to \End(\os(F\times G))
		\]
		to be the isomorphism induced by the equivalence $\widetilde{F}\times\widetilde{G}$. Then one can verify that 
		\[
			J_{F,G}=J_2J_1\,.
		\] 
		It suffices to show $J_1$ is invertible, which is indeed the case since the inverse can be given by 
\agn{
		K_1\:\End(\os(U^{\End(F)}\times U^{\End(G)})) & \to \End(F)\os\End(G) \\
			\gamma & \longmapsto \gamma_{\End(F),\End(G)}(\id_F\os\id_G)\,.
	}
\epf
\subsection{\texorpdfstring{Proof of \Cref{thm:reconstruct}.\ref{item1:thm:reconstruct}}{Proof of Theorem \ref{thm:reconstruct}.\ref{item1:thm:reconstruct}}}\label{subsec:pf:thm:reconstruct}
\pf[Proof of \Cref{thm:reconstruct}.\ref{item1:thm:reconstruct}]
	For simplicity, let us assume that the monoidal structure on $\Dd$ is strict. 
	
	First, we show that $(\End(\Fff),\Delta,\epsilon)$ form a weak bialgebra. It is not hard to conclude that $(\End(\Fff),\Delta,\epsilon)$ form a coalgebra. (\textbf{Axiom }\ref{item:Axiom_1}) follows from the separability condition of $\Fff$. The first equality in \eqref{eq:Axiom_2} of (\textbf{Axiom }\ref{item:Axiom_2}) holds since for any $\alpha,\beta,\gamma\in\End(\Fff)$, the outermost diagram of 
	\[
		\diagram{
			\kb\ar[r]^-{\Fff_0} & \Fff(\mone)  \ar[r]^-{\gamma_\mone} & \Fff(\mone) \ar[d]_1 \ar[r]^-{\Fff_0\os 1} & \Fff(\mone)\os \Fff(\mone) \ar[d]^{{\Fff_2}_{\mone,\mone}} \\
			& & \Fff(\mone) \ar[d]_{\beta_\mone} \ar[r]^-1 & \Fff(\mone\os\mone) \ar[d]^{\beta_{\mone\os\mone}}  \\
			& & \Fff(\mone) \ar[d]_1 & \Fff(\mone\os\mone) \ar[d]^{{\Fff_{-2}}_{\mone,\mone}} \ar[l]^-1 \\
			\kb & \Fff(\mone) \ar[l]^-{\Fff_{-0}} & \Fff(\mone) \ar[l]^-{\alpha_\mone} & \Fff(\mone)\os\Fff(\mone) \ar[l]^-{1\os\Fff_{-0}}
		}
	\]
	is commutative. The second equality in \eqref{eq:Axiom_2} can be proved similarly. That the first equality in \eqref{eq:Axiom_3} of (\textbf{Axiom }\ref{item:Axiom_3}) holds is equivalent to that for any $X,Y,Z\in\Dd$, the outermost diagram of
	\[
		\diagram{
			\Fff(X)\os\Fff(Y)\os\Fff(Z) \ar[d]_{1\os{\Fff_2}_{Y,Z}} \ar[r]^-1 & \Fff(X)\os\Fff(Y)\os\Fff(Z) \ar[d]^{1\os{\Fff_2}_{Y,Z}} \\
			\Fff(X)\os\Fff(Y\os Z) \ar[r]^-1 \ar[d]_{1\os{\Fff_{-2}}_{Y,Z}}  & \Fff(X)\os\Fff(Y\os Z) \ar[d]^{{\Fff_2}_{X,Y\os Z}} \\
			\Fff(X)\os\Fff(Y)\os\Fff(Z)  \ar@{}[r]|(0.5){\spadesuit} \ar[d]_{{\Fff_2}_{X,Y}\os 1} \ar[d] & \Fff(X\os Y\os Z) \ar[d]^{{\Fff_{-2}}_{X\os Y,Z}} \\
			\Fff(X\os Y)\os\Fff(Z) \ar[d]_{{\Fff_{-2}}_{X,Y}\os 1} \ar[r]^-1 & \Fff(X\os Y)\os\Fff(Z) \ar[d]^{{\Fff_{-2}}_{X,Y}\os 1} \\
			\Fff(X)\os \Fff(Y)\os\Fff(Z) \ar[r]^-1 & \Fff(X)\os \Fff(Y)\os\Fff(Z) \\
		}
	\]
	commutes. The diagram indeed commutes, where the commutativity of ($\scriptstyle\spadesuit$) comes from the Frobenius condition \eqref{diagram:frob_2} obeyed by $\Fff$. Similarly, the second equality of \eqref{eq:Axiom_3} can be derived using the the Frobenius condition \eqref{diagram:frob_1}.
	
	Secondly, let us show that $(\End(\Fff),\Delta,\epsilon,S)$ is a weak Hopf algebra. Note that checking \eqref{eq1:antipode} amounts to check that for any $\gamma\in\End(\Fff)$ and $X\in\Dd$, the outermost diagram of 
	\[
		\diagram@C=2.5pc{
			\Fff(X) \ar[rr]^-1 \ar[d]_{\Fff_0\os 1} & & \Fff(X) \ar[d]^{\Fff_0\os 1} \\
			\Fff(\mone)\os\Fff(X) \ar[d]_{\Fff(\coev)\os 1} \ar[rr]^-1 && \Fff(\mone) \os \Fff(X) \ar[dd]^{\gamma_\mone\os 1} \\
			\Fff(X\os X^L) \os \Fff(X) \ar[d]_{\gamma_{X\os X^L}\os 1}  \\
			\Fff(X\os X^L) \os \Fff(X) \ar[rd]^-{{\Fff_2}_{X\os X^L,X}} \ar[d]_{{\Fff_{-2}}_{X,X^L}\os 1}  && \Fff(\mone) \os \Fff(X) \ar[d]^{{\Fff_{2}}_{\mone,X}} \ar[ll]_-{\Fff(\coev)\os 1}\\
			\Fff(X)\os \Fff(X^L)\os\Fff(X) \ar@{}[r]|(0.5){\clubsuit} \ar[d]_{1\os {\Fff_2}_{X,X^L}} &\Fff(X\os X^L\os X) \ar[d]^{\Fff(1\os\ev)} \ar[ld]^-{\hskip 1.35em {\Fff_{-2}}_{X,X^L\os X}} & \Fff(X) \ar[l]_-{\Fff(\coev\os 1)}\ar[dd]^-1 \\
			\Fff(X)\os\Fff(X^L\os X) \ar[d]_{1\os \Fff(\ev)} & \Fff(X) \ar[rd]_-1 \ar[ld]^-{\hskip 0.35em {\Fff_{-2}}_{X,\mone}} \\
			\Fff(X)\os\Fff(\mone) \ar[rr]_-{1\os \Fff_{-0}} && \Fff(X)
		}
	\]
	is commutative. However, the diagram indeed commutes, where the commutativity of ($\scriptstyle\clubsuit$) follows from the Frobenius condition \eqref{diagram:frob_1}. Similarly, one can verify \eqref{eq2:antipode} using the Frobenius condition \eqref{diagram:frob_2}. Finally, note that $S$ is an algebra anti-homomorphism by construction. Thus \eqref{eq3:antipode} holds by \cite[Lemma 7.4]{Nill_1998}.
\epf

\section{\texorpdfstring{A comparison between $A_\Cc^{\Cc\boxtimes\Cc^\rev}$ and $\Tube_\Cc$}{A comparison between \$A\_C\^{}{CCrev}\$ and \$Tube\_C\$}}\label{sec:comparison}
It is well-known that for a pivotal fusion category $\Cc$, Ocneanu's tube algebra $\Tube_\Cc$ \cite{Ocneanu_1994} provides us with an equivalence \cite{Izumi_2000, Mueger_2003}
\eqn{
\label{eq:Tube_alg_eq}
\Zz(\Cc)\isom\rep(\Tube_\Cc)\,.
}
The combination of this equivalence with \Cref{cor:Drinfeld_center} establishes a Morita equivalence between $A_\Cc^{\Cc\boxtimes\Cc^\rev}$ and $\Tube_\Cc$, motivating a comparison of the two algebras, which we undertake in this appendix.

In \Cref{subsec:Tube_alg}, we review the tube algebra and the equivalence \eqref{eq:Tube_alg_eq}. In \Cref{subsec:Morita_eq}, we construct the explicit Morita equivalence between $A_\Cc^{\Cc\boxtimes\Cc^\rev}$ and $\Tube_\Cc$. In \Cref{subsec:not_wha}, we show that, in contrast to $A_\Cc^{\Cc\boxtimes\Cc^\rev}$, the tube algebra $\Tube_\Cc$ does not generally possess a weak Hopf algebra structure that would make the equivalence \eqref{eq:Tube_alg_eq} a monoidal equivalence.
\Cref{subsec:Morita_eq} and \Cref{subsec:not_wha} can be read independently. 

\subsection{\texorpdfstring{The equivalence between Drinfeld center and the representation category of $\Tube_\Cc$}{The equivalence between Drinfeld center and the representation category of \$Tube\_C\$}}\label{subsec:Tube_alg} 
\begin{dfn}[due to \cite{Ocneanu_1994}]
	Let $\Cc$ be a fusion category. The \newdef{(Ocneanu's) tube algebra associated with $\Cc$}, denoted by $\Tube_\Cc$, is an algebra over $\kb$ defined as follows:
	\bit
		\item As a vector space, we have
		\[
			\Tube_\Cc\defdtobe\Oplus_{x,y,w\in\Irr(\Cc)}\Cc(x\os w,w\os y)\,.
		\]
		We denote the subspace $\Cc(x\os w,w\os y)$ by $Y_{w|x|y}$ for $w,x,y\in\Irr(\Cc)$.
		\item For $(\diagram{x'\os w'\ar[r]^-{h} & w'\os y'})\in Y_{w'|x'|y'}$ and $(\diagram{x\os w\ar[r]^-g & w\os y})\in Y_{w|x|y}$, the multiplication $h\cdot g$ is given in two steps.
		\bnu[$1^\circ$]
			\item First, we choose a direct sum decomposition 
	\[
		w'\os w\cong\oplus_{t\in\Irr(\Cc)}t^{\oplus n_t}
	\]
	with inclusions $I^\alpha_t\:t\to w'\os w$ and projections $P^\alpha_t\:w'\os w\to t$ for $t\in\Irr(\Cc)$ and $\alpha=1,\cdots,n_t$. 
	
			\item Secondly, we have 
	\[
		h\cdot g=\delta_{x,y'}\sum_{t\in\Irr(\Cc)}(h\cdot g)_{t|x'|y}\,,
	\]
	where $(h\cdot g)_{t|x'|y}\in Y_{t|x'|y}$ reads
	\[
	\begin{multlined}
		\sum_{\alpha=1}^{n_t}(\diagram{x'\os t\ar[r]^-{1\os I^\alpha_t} & x'\os w'\os w \ar[r]^-{h\os 1} & w'\os y'\os w \ar[r]^-{1\os g} & w'\os w\os y\ar[r]^-{P^\alpha_t\os 1} & t\os y})\,.
	\end{multlined}
	\]
		\enu
		\item The unit reads
		\[
			1=\sum_{x\in\Irr(\Cc)}(\diagram{x\ar[r]^-{\id_x} & x})\,,
		\]
		where $\id_x\in Y_{\mone|x|x}$.
	\eit
\end{dfn}
Recall that a pivotal structure on $\Cc$ is a monoidal natural isomorphism $\fa\:\Id_\Cc\funto (-)^{RR}$ (see for e.g. \cite[Section 4.7]{Etingof_Gelaki_Nikshych_Ostrik_2015}). In this subsection, we show
\begin{thm}[\cite{Izumi_2000, Mueger_2003}]\label{thm:rep_of_Tube}
	Let $\Cc$ be a fusion category with pivotal structure $\fa$. There exists an equivalence of categories
	\eqn{\label{eq:Tube_eq_of_cat_3}
	\begin{split}
		J\:\Zz(\Cc) & \to\rep(\Tube_\Cc) \\
		(z,\gamma_{-,z}) & \longmapsto \Oplus_{x\in\Irr(\Cc)}\Cc(x,z)\,,
	\end{split}
	}
	where the $\Tube_\Cc$-action on $J((z,\gamma_{-,z}))$ is given as follows: for morphisms $(\diagram{x\os w\ar[r]^-g & w\os y})\in Y_{w|x|y}$ and $\diagram{x_0\ar[r]^-s & z}$ such that $x_0$ is simple, the action $g.s$ reads
	\eqnn{
	\begin{split}
		\delta_{x_0,y}(\diagram@C=3pc{
			x \ar[r]^-{1\os\coev} & x\os w\os w^L \ar[r]^-{g\os 1} & w\os y\os w^L\ar[r]^-{1\os s\os 1} & w\os z\os w^L} \\
			\diagram@C=3pc{ \ar[r]^-{\gamma_{w,z}\os 1} & z\os w\os w^L \ar[r]^-{1\os \fa_{w^L}} & z\os w\os w^R \ar[r]^-{1\os \ev} & z})\,.
	\end{split}
	}
\end{thm}

Most proofs of \Cref{thm:rep_of_Tube} are developed within the context of operator algebras or topological quantum field theory, where it is common to assume that $\Cc$ satisfies the additional conditions of being unitary or spherical. We will present a purely algebraic proof of \Cref{thm:rep_of_Tube}, showing that a pivotal structure is sufficient.

Along the way, we present a ``reconstruction viewpoint'' on the tube algebra. This viewpoint can be used to reproduce the family of algebras whose representation categories are all $\Zz(\Cc)$ in \cite[Remark 5.1.2]{Mueger_2003}, and it will be exploited in the next subsection to show the Morita equivalence result.

Let $\Cc$ be a fusion category. Denote $a_0\defdtobe\oplus_{x\in\Irr(\Cc)}x\in\Cc$. The representable functor
\agn{
	H\:\Cc & \to \vect \\
	a & \longmapsto \Cc(a_0,a)\equiv\Oplus_{x\in\Irr(\Cc)}\Cc(x,a)
}
is faithful and exact, where we use the fact that an exact functor between abelian categories is faithful if and only if it preseves non-zero objects. It is also known that the forgetful functor
\agn{
	G\:\Zz(\Cc) & \to\Cc\\
	(z,\gamma_{-,z}) & \longmapsto z
}
is faithful and exact \cite[Proposition 1]{Street_1998}. Thus we obtain a faithful and exact functor
\agn{
	HG\:\Zz(\Cc) & \to \vect \\
	(z,\gamma_{-,z}) & \longmapsto \Oplus_{x\in\Irr(\Cc)}\Cc(x,z)\,.
}
By \Cref{thm:monadicity}, we have an equivalence of categories
\eqn{\label{eq:Tube_eq_of_cat_1}
\begin{split}
	\Zz(\Cc) & \to\rep(\End(HG)) \\
	(z,\gamma_{-,z}) & \longmapsto \Oplus_{x\in\Irr(\Cc)}\Cc(x,z)\,.
\end{split}
}
\Cref{thm:rep_of_Tube} can now be immediately proved once we have an isomorphism $\End(HG)\isom\Tube_\Cc$ of algebras. It is hence worthwhile to have a presentation of the algebra $\End(HG)$.
\begin{lem}[\cite{Day_Street_2007, Bruguieres_Virelizier_2012}]\label{lem:Z_C_coend}
	The functor $G$ admits a left adjoint. If $F$ denotes this adjoint, then the underlying object of $F(a)$ for $a\in\Cc$ is given by $\oplus_{x\in\Irr(\Cc)}x\os a\os x^R$.
\end{lem}
\begin{rmk}
	The key observation in \cite{Day_Street_2007, Bruguieres_Virelizier_2012} is that when $\Cc$ has certain nice properties, the forgetful functor $G\:\Zz(\Cc)\to\Cc$ is monadic. The monad sends $a\in\Cc$ to the coend $\int^{x\in\Cc}x\os a\os x^R$, which reduces to $\oplus_{x\in\Irr(\Cc)}x\os a\os x^R$ when $\Cc$ is semisimple. 
\end{rmk}
\begin{cor}\label{cor:Tube_alg}
	\bnu
		\item The functor $HG$ is represented by $F(a_0)$.
		\item We have 
		\eqn{
		\label{eq:Tube_alg_wo_pivotal}
		\End(HG)\cong\Cc(a_0,GF(a_0))=\Oplus_{x,y,w\in\Irr(\Cc)}\Cc(x,w\os y\os w^R)
		}
		as vector spaces.
	\enu
	\pf
		1. For any $(z,\gamma_{-,z})\in\Zz(\Cc)$, we have $\Zz(\Cc)(F(a_0),(z,\gamma_{-,z}))\cong\Cc(a_0,z)\cong HG((z,\gamma_{-,z}))$. 2. By Yoneda lemma, we have $\End(HG)\cong\Zz(\Cc)(F(a_0),F(a_0))\cong\Cc(a_0,GF(a_0))$ as vector spaces.
	\epf
\end{cor}
Using the isomorphism \eqref{eq:Tube_alg_wo_pivotal}, we can transport the algebra structure on $\End(HG)$ to the space $\Tube'_\Cc\defdtobe\Oplus_{x,y,w\in\Irr(\Cc)}\Cc(x,w\os y\os w^R)$. Moreover, the equivalence \eqref{eq:Tube_eq_of_cat_1} extends to an equivalence between $\Zz(\Cc)$ and $\rep(\Tube'_\Cc)$. To present these data, it is helpful to denote $X_{w|x|y}\defdtobe\Cc(x,w\os y\os w^R)\subset \Tube'_\Cc$ for $x,y,w\in\Irr(\Cc)$. By computation, we obtain the following proposition:
\begin{prp}\label{prp:rep_of_pre_Tube}
\bnu
	\item There is an algebra structure on $\Tube'_\Cc$ defined as follows.
	\bit
		\item For $(\diagram{x'\ar[r]^-h& w'\os y'\os w'^R})\in X_{w'|x'|y'}$ and $(\diagram{x\ar[r]^-g& w\os y\os w^R})\in X_{w|x|y}$, the multiplication $h\cdot g$ is given in two steps:
	\bnu[$1^\circ$]
		\item First, we choose a direct sum decomposition 
	\[
		w'\os w\cong\oplus_{t\in\Irr(\Cc)}t^{\oplus n_t}
	\]
	with inclusions $I^\alpha_t\:t\to w'\os w$ and projections $P^\alpha_t\:w'\os w\to t$ for $\alpha=1,\cdots,n_t$ and $t\in\Irr(\Cc)$. 
	
		\item Then 
	\[
		h\cdot g=\delta_{x,y'}\sum_{t\in\Irr(\Cc)}(h\cdot g)_{t|x'|y}\,,
	\]
	where $(h\cdot g)_{t|x'|y}\in X_{t|x'|y}$ reads
	\[
	\begin{multlined}
		\sum_{\alpha=1}^{n_t}(\diagram@C=2.1pc{x' \ar[r]^-h & w'\os y'\os w'^R \ar[r]^-{1\os g\os 1} & w'\os w\os y\os w^R\os w'^R \ar[rr]^-{P^\alpha_t\os 1\os (I^\alpha_t)^R} & & t\os y\os t^R})\,.
	\end{multlined}
	\]
	\enu
		\item The unit reads
	\[
		1=\sum_{x}(\diagram{x\ar[r]^-{\id_x} & x})\,,
	\]
	where $\id_x\in X_{\mone|x|x}$. 
	\eit
	\item There is an equivalence of categories
	\eqn{\label{eq:Tube_eq_of_cat_2}
	\begin{split}
		\Zz(\Cc) &  \to \rep(\Tube'_\Cc) \\
		(z,\gamma_{-,z}) & \longmapsto \Oplus_{x\in\Irr(\Cc)}\Cc(x,z)\,,
	\end{split}
	}
	where the $\Tube'_\Cc$-action on $\Oplus_{x\in\Irr(\Cc)}\Cc(x,z)$ is given as follows: for morphisms 
	\[(\diagram@C=1.5pc{x\ar[r]^-g& w\os y\os w^R})\in X_{w|x|y} \qquad \text{and}\qquad \diagram{x_0\ar[r]^-s & z}
	\]
	such that $x_0$ is simple, the action $g.s$ reads
	\[
		\delta_{x_0,y}(\diagram@C=2.5pc{x \ar[r]^-g  & w\os y\os w^R \ar[r]^-{1\os s\os 1} & w\os z\os w^R \ar[r]^-{\gamma_{w,z}\os 1} & z\os w\os w^R \ar[r]^-{1\os \ev} & z})\,.
	\]
\enu
\end{prp}

Now we're ready to prove \Cref{thm:rep_of_Tube}.
\pf[Proof of \Cref{thm:rep_of_Tube}]
 When $\Cc$ is equipped with a pivotal structure $\fa$, it is clear that the linear isomorphism
\eqn{\label{eq:pre_tube_to_tube}
\begin{multlined}
	\Tube'_\Cc=\Oplus_{x,y,w\in\Irr(\Cc)}\Cc(x,w\os y\os w^R)\isom\Oplus_{x,y,w\in\Irr(\Cc)}\Cc(x\os w^{RR},w\os y) \\\isom\Oplus_{x,y,w\in\Irr(\Cc)}\Cc(x\os w,w\os y)=\Tube_\Cc
\end{multlined}
}
induced by $\fa_w\:w\isom w^{RR}$ is an algebra isomorphism. Using this isomorphism and the equivalence \eqref{eq:Tube_eq_of_cat_2} in \Cref{prp:rep_of_pre_Tube}, the equivalence \eqref{eq:Tube_eq_of_cat_3} is immediately obtained.
\epf
\subsection{\texorpdfstring{Morita equivalence between $\Tube_\Cc$ and $A_\Cc^{\Cc\boxtimes\Cc^\rev}$}{Morita equivalence between \$Tube\_C\$ and \$A\_C\^{}{CCrev}\$}}\label{subsec:Morita_eq}
In this subsection, we sketchily prove the following
\begin{thm}\label{thm:Mor}
\bnu
	\item \label{item1:thm:Mor} Let $\Cc$ be a fusion category. Then there is a sequence of mutually Morita equivalent algebras
	\[
		\Tube'^{(1)}_\Cc,\Tube'^{(2)}_\Cc,\cdots,\Tube'^{(n)}_\Cc,\cdots
	\]
	such that $\Tube'^{(1)}_\Cc=\Tube'_\Cc$ and the underlying vector space of $\Tube'^{(n)}_\Cc$ is 
	\[
		\manualformatting \hskip -1.25em\Oplus_{\substack{x_1,\cdots,x_n,\\y_1,\cdots,y_n,w\in\Irr(\Cc)}}\hskip -1.25em\Cc(x_1\os\cdots\os x_n,w\os y_1\os\cdots\os y_n\os w^R)\,.
	\]
	For any $n,m\geq 1$, there exists an invertible $\Tube'^{(n)}_\Cc$-$\Tube'^{(m)}_\Cc$-bimodule whose underlying vector space is given by 
	\[\manualformatting\hskip -1.25em\Oplus_{\substack{x_1,\cdots,x_n,\\y_1,\cdots,y_m,w\in\Irr(\Cc)}}\hskip -1.25em\Cc(x_1\os\cdots\os x_n,w\os y_1\os\cdots\os y_m\os w^R)\,.
	\]
	
	On the other hand, there exists a sequence of mutually Morita equivalent algebras
	\[
		\Tube_\Cc^{(1)},\Tube_\Cc^{(2)},\cdots,\Tube_\Cc^{(n)},\cdots
	\]
	such that the $\Tube_\Cc^{(1)}=\Tube_\Cc$ and the underlying vector space of $\Tube_\Cc^{(n)}$ is 
	\eqn{\label{eq:thm:Mor}
		\manualformatting\hskip -1.25em\Oplus_{\substack{x_1,\cdots,x_n,\\y_1,\cdots,y_n,w\in\Irr(\Cc)}}\hskip -1.25em\Cc(x_1\os\cdots\os x_n\os w,w\os y_1\os\cdots\os y_n)\,.
	}
	For any $n,m\geq 1$, there exists an invertible $\Tube_\Cc^{(n)}$-$\Tube_\Cc^{(m)}$-bimodule whose underlying vector space is given by 
	\eqn{\label{eq2:thm:Mor}
	\manualformatting\hskip -1.25em\Oplus_{\substack{x_1,\cdots,x_n,\\y_1,\cdots,y_m,w\in\Irr(\Cc)}}\hskip -1.25em\Cc(x_1\os\cdots\os x_n\os w,w\os y_1\os\cdots\os y_m)\,.
	}
	
	In addition, we have $\Tube'^{(2)}_\Cc\cong A_\Cc^{\Cc\boxtimes\Cc^\rev}$ as algebras.
	
	\item \label{item2:thm:Mor} If $\Cc$ is pivotal, then for any $n\geq 1$ there is an algebra isomorphism
	\[
		\Tube'^{(n)}_\Cc\isom\Tube_\Cc^{(n)}\,.
	\]
	In particular, $A_\Cc^{\Cc\boxtimes\Cc^\rev}$ and $\Tube_\Cc$ are Morita equivalent.
\enu
\end{thm}
Note that the two sequences of Morita equivalent algebras are not new. When $\Cc$ is pivotal, they form a subfamily of the algebras appearing in \cite[Remark 5.1.2]{Mueger_2003}. A variant of these algebras can also be seen in, for example, \cite[Lemma 2]{Kong_2012}.

To prove \Cref{thm:Mor}, let us first construct the family of algebras $\{\Tube'^{(n)}_\Cc\}_{n\geq 1}$.
We construct each member in this family in a way similarly to our construction of $\Tube'_\Cc$ in \Cref{subsec:Tube_alg}. Recall the object $a_0=\oplus_{x\in\Irr(\Cc)}x$. It can be easily verified that the representable functor
\[
	H^{(n)}\defdtobe\Cc(a_0^{\os n},-)\:\Cc\to\vect
\]
is faithful and exact, as is $H^{(1)}=H$ in \Cref{subsec:Tube_alg}. Then $H^{(n)}G\:\Zz(\Cc)\to\vect$ is faithful and exact, hence we have 
\[\Zz(\Cc)\cong\rep(\End(H^{(n)}G))\]
as categories. Similarly to \Cref{cor:Tube_alg}, we can show that $H^{(n)}G$ is represented by $F(a_0^{\os n})$, where $F$ is the left adjoint of $G$. Moreover, we have an isomorphism
\[\End(H^{(n)}G)\cong\Cc(a_0^{\os n},GF(a_0^{\os n}))=\manualformatting\hskip -1.25em\Oplus_{\substack{x_1,\cdots,x_n,\\y_1,\cdots,y_n,w\in\Irr(\Cc)}}\hskip -1.25em\Cc(x_1\os\cdots\os x_n,w\os y_1\os\cdots\os y_n\os w^R)
\]
of vector spaces. Using this isomorphism, we can transport the algebra structure on $\End(H^{(n)}G)$ to the space $\Tube'^{(n)}_\Cc\defdtobe\Oplus_{x_1,\cdots,x_n,y_1,\cdots,y_n,w\in\Irr(\Cc)}\Cc(x_1\os\cdots\os x_n,w\os y_1\os\cdots\os y_n\os w^R)$. The explicit expression of the algebra structure on $\Tube'^{(n)}_\Cc$ is similar to that on $\Tube'^{(1)}_\Cc=\Tube'_\Cc$, and is not given here. The algebras $\Tube'^{(n)}_\Cc$ are all Morita equivalent since their representation categories are all $\Zz(\Cc)$. An invertible $\Tube'^{(n)}_\Cc$-$\Tube'^{(m)}_\Cc$-bimodule can be given by  
\eqnn{
\begin{multlined}
\nat(H^{(m)}G,H^{(n)}G)\cong\Zz(\Cc)(F(a_0^{\os n}),F(a_0^{\os m}))\cong\Cc(a_0^{\os n}, GF(a_0^{\os m}))
\\\cong\manualformatting\hskip -1.25em\Oplus_{\substack{x_1,\cdots,x_n,\\y_1,\cdots,y_m,w\in\Irr(\Cc)}}\hskip -1.25em\Cc(x_1\os\cdots\os x_n,w\os y_1\os\cdots\os y_m\os w^R)
\end{multlined}
}
with the evident $\End(H^{(n)}G)$-$\End(H^{(m)}G)$-action, where $\nat(K,K')$ denotes the vector space of natural transformations $K\funto K'$ for $\kb$-linear functors $K$ and $K'$.

We now turn to the existence of the family of algebras $\{\Tube_\Cc^{(n)}\}_{n\geq 1}$ in \Cref{thm:Mor}. For each $n\geq 1$, we define the algebra $\Tube_\Cc^{(n)}$ in a manner similar to that of $\Tube_\Cc$ as follows:
\bit
	\item The underlying vector space of $\Tube_\Cc^{(n)}$ is given by \eqref{eq:thm:Mor}. For $x_1,\cdots,x_n,y_1,\cdots,y_n,w\in\Irr(\Cc)$, we denote the subspace 
\[
	\Cc(x_1\os\cdots\os x_n\os w,w\os y_1\os\cdots\os y_n)\subset\Tube_\Cc^{(n)}
\]
as $Y_{w|x_1|\cdots|x_n|y_1|\cdots|y_n}$.
	\item For \[(\diagram{x_1'\os \cdots\os x_n'\os w'\ar[r]^-{h} & w'\os y_1'\os \cdots \os y_n'})\in Y_{w'|x_1'|\cdots|x_n'|y_1'|\cdots|y_n'}\] and \manualformatting \[(\diagram{x_1\os \cdots\os x_n\os w\ar[r]^-{g} & w\os y_1\os \cdots \os y_n})\in Y_{w|x_1|\cdots|x_n|y_1|\cdots|y_n}\,,\]
	the multiplication $h\cdot g$ is given in two steps. 
	\bnu[$1^\circ$]
		\item First, choose a direct sum decomposition 
	\[
		w'\os w\cong\oplus_{t\in\Irr(\Cc)}t^{\oplus n_t}
	\]
	with inclusion maps $I^\alpha_t\:t\to w'\os w$ and projection maps $P^\alpha_t\:w'\os w\to t$ for $\alpha=1,\cdots,n_t$ and $t\in\Irr(\Cc)$.
	
		\item Then we have 
	\[
		h\cdot g=(\prod_{i=1}^n\delta_{x_i,y_i'})\sum_{t\in\Irr(\Cc)}(h\cdot g)_{t|x_1'|\cdots|x_n'|y_1|\cdots|y_n}\,,
	\]
	where $(h\cdot g)_{t|x_1'|\cdots|x_n'|y_1|\cdots|y_n}\in Y_{t|x_1'|\cdots|x_n'|y_1|\cdots|y_n}$ reads
	\[
	\begin{multlined}
		\sum_{\alpha=1}^{n_t}(\diagram@C=2.25pc{x_1'\os \cdots \os x_n'\os t\ar[r]^-{1\os I^\alpha_t} & x_1'\os \cdots \os x_n'\os w'\os w \ar[r]^-{h\os 1} & w'\os y_1'\os \cdots\os y_n' \os w }\\
		\diagram@C=2.25pc{\ar[r]^-{1\os g} & w'\os w\os y_1\os \cdots\os y_n \ar[r]^-{P^\alpha_t\os 1} & t\os y_1\os\cdots\os y_n})\,.
	\end{multlined}
	\]
	\enu
	\item The unit is given by $\sum_{x_1,\cdots,x_n\in\Irr(\Cc)}\id_{x_1\os\cdots\os x_n}$, where $\id_{x_1\os\cdots\os x_n}\in Y_{\mone|x_1|\cdots|x_n|x_1|\cdots|x_n}$.
\eit

Let us show that for $n,m\geq 1$, there exists a $\Tube_\Cc^{(n)}$-$\Tube_\Cc^{(m)}$ bimodule with the underlying vector space given by
\[
		\Tube_\Cc^{(m,n)}\defdtobe\manualformatting\hskip -1.25em\Oplus_{\substack{x_1,\cdots,x_n,\\y_1,\cdots,y_m,w\in\Irr(\Cc)}}\hskip -1.25em\Cc(x_1\os\cdots\os x_n\os w,w\os y_1\os\cdots\os y_m)
\]
as in \eqref{eq2:thm:Mor}. 

Our proof of this fact is modified from a proof of \cite[Lemma 2]{Kong_2012} given in \cite{Kong}. First, note that a map 
\[
	\circ^{mnk}\:\Tube_\Cc^{(n,k)}\os\Tube_\Cc^{(m,n)}\to\Tube_\Cc^{(m,k)}\,.
\]
can be defined in a similiar way as the multiplication of $\Tube_\Cc^{(n)}$, so that $(\Tube_\Cc^{(n,n)},\circ^{nnn})$ is precisely the algebra $\Tube_\Cc^{(n)}$. One can check that $\{\circ^{mnk}\}_{m,n,k\geq 1}$ satisfy the generalized associativity constraints
\[
	\circ^{mnl}\circ (\circ^{nkl}\os\id)=\circ^{mkl}\circ(\id\os\circ^{mnk}),\;\forall m,k,n,l\geq 1\,,
\]
and certain generalized unitality constraints. In particular, the vector space $\Tube_\Cc^{(m,n)}$ carries a $\Tube_\Cc^{(n)}$-$\Tube_\Cc^{(m)}$-bimodule action, and 
\[
	\circ^{nmn}\:\Tube_\Cc^{(m,n)}\os\Tube_\Cc^{(n,m)}\to\Tube_\Cc^{(n,n)}
\]
is a $\Tube_\Cc^{(m,m)}$-balanced map. Now we briefly show that $\Tube_\Cc^{(m,n)}$ is an invertible bimodule. It suffices to show that $\circ^{nmn}$ exhibits $\Tube_\Cc^{(n,n)}$ as the relative tensor product
\[
	\Tube_\Cc^{(m,n)}\os_{\Tube_\Cc^{(m,m)}} \Tube_\Cc^{(n,m)}\,.
\]
To this end, define a map
\[
	s\:\Tube_\Cc^{(n,n)}\to\Tube_\Cc^{(m,n)}\os\Tube_\Cc^{(n,m)}
\]
by setting $s(g)$ for $(\diagram{x_1\os \cdots\os x_n\os w\ar[r]^-{g} & w\os y_1\os \cdots \os y_n})\in Y_{w|x_1|\cdots|x_n|y_1|\cdots|y_n}$ to be
\[
	\sum_{d\in\Irr(\Cc)}\sum_{\alpha=1}^{n_d}g^\alpha_{d,1}\os g^\alpha_{d,2}\,.
\]
Here 
\agn{
	g^\alpha_{d,1} & =(\diagram@C=2.5pc{x_1\os\cdots \os x_n\os w\ar[r]^-g & w\os y_1\os \cdots \os y_m \ar[r]^-{1\os P_d^\alpha} & w\os d\os \overbrace{\mone\os\cdots\os\mone}^{m-1}})\in\Tube_\Cc^{(m,n)} \\
	g^\alpha_{d,2} & =(\diagram@C=2.5pc{d\os\overbrace{\mone\os\cdots\os\mone}^{m-1}\os\mone \ar[r]^-{1\os I^\alpha_d} & \mone\os y_1\os\cdots\os y_n})\in\Tube_\Cc^{(n,m)}\,,
}
where $I^\alpha_d\:d\to y_1\os\cdots\os y_m$ and $P^\alpha_d\:y_1\os\cdots\os y_m\to d$ are respectively the inclusions and the projections in the direct sum decomposition $y_1\os\cdots\os y_m\cong\oplus_{d\in\Irr(\Cc)}d^{\oplus n_d}$ for $\alpha=1,\cdots,n_d$.

One can check that $s$ is a section of $\circ^{nmn}$. Moreover, for any vector space $Q$ and a $\Tube_\Cc^{(m,m)}$-balanced map $q\:\Tube_\Cc^{(m,n)}\os\Tube_\Cc^{(n,m)}\to Q$, we have that $\underline{q}\defdtobe q\circ s$ satisfies
\eqn{\label{eq3:thm:Mor}
	\underline{q}\circ \circ^{nmn}=q\,.
}
On the other hand, by that $s$ is a section of $\circ^{nmn}$, one can verify that $\underline{q}=q\circ s$ is the unique map satsifying \eqref{eq3:thm:Mor}. This establishes the proof of $\Tube_\Cc^{(n,n)}\cong\Tube_\Cc^{(m,n)}\os_{\Tube_\Cc^{(m,m)}} \Tube_\Cc^{(n,m)}$.

Let us first finish the proof of \Cref{thm:Mor}.\ref{item2:thm:Mor}. Observe that when $\Cc$ is equipped with a pivotal structure $\fa$, a linear isomorphism
\[
	\Tube'^{(n)}_\Cc\isom\Tube_\Cc^{(n)}
\]
can be constructed for all $n\geq 1$, utilizing $\fa$ in a way similiar to the isomorphism \eqref{eq:pre_tube_to_tube}. It is easy to conclude that this is an algebra isomorphism. This proves \Cref{thm:Mor}.\ref{item2:thm:Mor}.

To finish the proof of \Cref{thm:Mor}.\ref{item1:thm:Mor}, it suffices to show that there exists an algebra isomorphism $A_\Cc^{\Cc\boxtimes\Cc^\rev}\isom \Tube'^{(2)}_\Cc$ when $\Cc$ is a fusion category. To this end, let us denote the subspace 
\[
\Cc(x_1\os x_2,w\os y_1\os y_2\os w^R)\subset\Tube'^{(2)}_\Cc
\]
for simple objects $x_1,x_2,y_1,y_2,w\in\Irr(\Cc)$ by $X_{w|x_1|x_2|y_1|y_2}$. Then a linear map $\chi\:A_\Cc^{\Cc\boxtimes\Cc^\rev}\to \Tube'^{(2)}_\Cc$ can be defined as follows. For simple objects $a,b,y',y,x',x\in\Irr(\Cc)$ and morphisms
\[
	\diagram{y'\ar[r]^-u & a\os y\os b} \quad\text{and}\quad \diagram{a\os x\os b\ar[r]^-s & x'}\,,
\]
we set $\chi(u\os s)$ to be the following element in $X_{a|y'|x'^R|y|x^R}$:
\[
	\begin{multlined}
		\diagram@C=2.25pc{
	y'\os x'^R\ar[r]^-{u\os 1} & a\os y\os b\os x'^R \ar[rr]^-{1\os\coev\os 1} & & a\os y\os x^R\os a^R\os a\os x\os b \os x'^R} \\
	\diagram@C=3.25pc{ \ar[r]^-{1\os s\os 1} & a\os y\os x^R\os a^R\os x'\os x'^R \ar[r]^-{1\os\ev} & a\os y\os x^R\os a^R}.
	\end{multlined}
\]
\Cref{thm:Mor} is then proved once the following easily verifiable observation is made:
\begin{prp}
	The map $\chi\:A_\Cc^{\Cc\boxtimes\Cc^\rev}\to \Tube'^{(2)}_\Cc$ is an algebra isomorphism.
\end{prp}

\subsection{\texorpdfstring{$\Tube_\Cc$ is in general not a weak Hopf algebra}{\$Tube\_C\$  is in general not a weak Hopf algebra}}\label{subsec:not_wha}

In this subsection, we observe that $\Tube_\Cc$ in general does not possess a weak Hopf algebra structure rendering \eqref{eq:Tube_eq_of_cat_3} a monoidal equivalence. This is in contrast with the scenario in \Cref{cor:Drinfeld_center}.
\begin{prp}
	There exists a pivotal fusion category $\Cc$ satisfying the following property: there is no weak bialgebra structure on $\Tube_\Cc$ such that the induced monoidal structure on $\rep(\Tube_\Cc)$ renders \eqref{eq:Tube_eq_of_cat_3} a monoidal equivalence.
	\pf
		Let $\rep(\Tube_\Cc)$ carry a monoidal structure $(\bos,\mathbb{1})$ such that the equivalence $J$ in \eqref{eq:Tube_eq_of_cat_3} is a monoidal equivalence. In order that this monoidal structure is induced by a weak bialgebra structure on $\Tube_\Cc$, for any $z=(z,\gamma_{-,z}),z'=(z',\gamma'_{-,z'})\in\Zz(\Cc)$, there should be 
		\[
		\dim J(z\os z')=\dim (J(z)\bos J(z'))\leq \dim (J(z)\os J(z'))=\dim J(z) \cdot \dim 	J(z')\,.
		\]
Here, the inequality follows from the definition of the tensor product of two representations over a weak bialgebra recalled in \Cref{subsec:wha_to_wff}. Take $\Cc$ to be the Fibonacci modular tensor category $\fib$ defined in \cite[\S 5.3.2]{Rowell_Stong_Wang_2009}, which is in particular a pivotal fusion category. It has two simple object $\mone$ and $\nu$, with the fusion rule $\nu\os\nu=\mone\oplus\nu$. Let $c$ denote the braiding of $\fib$. Then $z\defdtobe(\nu,c_{-,\nu})$ defines an object in $\Zz(\Cc)$. Now we have $\dim J(z\os z)=\dim J(1\oplus \nu)=2>1\cdot 1=\dim J(z)\cdot \dim J(z)$. Therefore, there is no weak bialgebra structure on $\Tube_\Cc$ rendering \eqref{eq:Tube_eq_of_cat_3} a monoidal equivalence when $\Cc=\fib$.
	\epf
\end{prp}
\begin{rmk}
Nonetheless, there exist pivotal fusion categories $\Cc$ such that $\Tube_\Cc$ has a weak Hopf algebra structure rendering \eqref{eq:Tube_eq_of_cat_3} a monoidal equivalence. One well known example is given by $\Cc=\vectG$ for a finite group $G$, in which case $\Tube_\Cc$ is the Drinfeld double $D[G]$ of the group algebra $\kb[G]$.
\end{rmk}

\addcontentsline{toc}{section}{References}

\bibliographystyle{alpha}
\bibliography{ref_LWKK}

\end{document}